\documentclass[12pt]{article}

\usepackage[a4paper,text={17.5cm,24.0cm},centering]{geometry}
\usepackage{latexsym}
\usepackage{graphicx}
\usepackage{amsmath}
\usepackage{amssymb}
\usepackage{mathrsfs}
\usepackage{bm}
\usepackage{fancybox}
\usepackage{color}
\usepackage{wrapfig}
\usepackage{framed}
\usepackage{mathtools}
\usepackage{framed}
\DeclareMathAlphabet{\mathpzc}{OT1}{pzc}{m}{it}
\def\intp{\phantom{al}\setlength{\unitlength}{0.4mm}
                 \begin{picture}(0, 0)(5, 5)%
                 \put(0, 4){\line(1,0){5}} \put(5,4){\line(0,1){7}}
                 \end{picture}\,\,}

\def\lprod{[\![}
\def\rprod{]\!]}

\def\delt{\partial_t}

\def\delx{\partial_x}

\def\xd{\mathrm{d}}

\def\phat{\hat{\phi}}

\def\Zhat{\widehat{Z}}
\def\util{\tilde{u}}

\def\sech{\mathrm{sech\,}}

\def\J{\mathbf{J}}

\def\C{\mathbb{C}}

\def\lbk{\lbrack\!\lbrack}
\def\rbk{\rbrack\!\rbrack}

\def\q'{q\;\!'\!\!\;}

\def\D{{\mathcal D}}

\def\bM{{\bf M}}
\def\bK{{\bf K}}
\def\ri{\texttt{i}}
\def\rd{{\rm d}}
\def\re{{\rm e}}
\def\fr{\mbox{$\frac{1}{2}$}}
\def\frr{\mbox{$\frac{1}{4}$}}
\def\R{\mathbb R}
\def\C{\mathbb C}

\def\vol{{\textsf{vol}}}

\def\qand{\quad\mbox{and}\quad}

\def\Rg{{\mathcal R}}
\def\bO{{\bm\Omega}}
\def\Dform{{\mathcal A}}

\def\Gt{{\rm G}_\theta}
\def\Gtt{{\rm G}_{\theta(t)}}
\def\g{\mathpzc{g}}
\def\Tq{S^1}
\def\LUw{{\bf L}({\cal U},c)}
\def\Zh{{\widehat Z}}

\def\bu{{\mathbf{u}}}
\def\u{{\mathbf{u}}}
\def\bv{{\mathbf{v}}}

\def\bw{{\mathbf{w}}}
\def\ba{{\mathbf{a}}}

\def\bp{{\mathbf{p}}}
\def\bq{{\mathbf{q}}}

\def\bwedge{{\textstyle\bigwedge}}

\def\B{{\mathbf B}}
\def\M{{\mathbf M}}
\def\K{{\mathbf K}}

\def\w{{\bf w}}
\def\a{{\bf a}}

\def\ptil{\tilde{\phi}}
\def\util{\tilde{u}}
\def\vtil{\tilde{v}}

\renewcommand{\theequation}{\arabic{section}.\arabic{equation}}

\begin{document}

\begin{center}
  \textbf{\Large Symplectic transversality and the Pego-Weinstein theory}
\vspace{.70cm}

\textsf{\large Timothy J. Burchell\footnote{Email: \textcolor{blue}{\texttt{T.Burchell@surrey.ac.uk}}\hfill\textcolor{blue}{\today}}, {\small and} Thomas J. Bridges\footnote{Email: \textcolor{blue}{\texttt{T.Bridges@surrey.ac.uk}}}
  }
\vspace{.25cm}

\textit{Department of Mathematics, University of Surrey, Guildford GU2 7XH, UK}
\vspace{0.75cm}

\setlength{\fboxsep}{10pt}
\doublebox{\parbox{15cm}{
{\bf Abstract.} {\small
This paper studies the linear stability problem for solitary wave solutions of
Hamiltonian PDEs.  The linear stability problem is
formulated in terms of the Evans function, a complex analytic function
denoted by $D(\lambda)$, where $\lambda$ is the spectral parameter.
The main result is the introduction of a new factor, denoted $\Pi$,
in the Pego \& Weinstein (1992) derivative formula
\[
D''(0) = \chi \Pi \frac{dI}{dc}\,,
\]
where $I$ is the momentum of the solitary wave and $c$ is the speed.
Moreover this factor turns out to be
related to transversality of the solitary wave,
modelled as a homoclinic orbit: the homoclinic orbit is transversely constructed
if and only if $\Pi\neq 0$.  The sign of $\Pi$ is a symplectic
invariant, an intrinsic property of the solitary wave, and is a
key new factor affecting the linear stability.  The factor $\chi$ was
already introduced by Bridges \& Derks (1999) and is based on the asymptotics of
the solitary wave.
A supporting result is
the introduction of a new abstract class of Hamiltonian PDEs
built on a nonlinear Dirac-type equation, which model a wide range of PDEs
in applications.  Examples
where the theory applies, other than Dirac operators, are the coupled mode equation
in fluid mechanics and optics, the massive Thirring model, and coupled
nonlinear wave equations. A calculation of $D''(0)$ for
solitary wave solutions of the latter class is included to illustrate the theory.}
}}
\end{center}


\section{Introduction}
\setcounter{equation}{0}
\label{sec-intro}

The stability of solitary waves of Hamiltonian partial differential equations
(PDEs) can be approached many ways.  One strategy is to use the
calculus of variations, since solitary waves can be characterized as critical points of the energy restricted to level sets of the momentum, and show that
the solitary wave is a minimizer on the constraint set, concluding,
with some additional analysis, Lyapunov (nonlinear, orbital) stability.
For dispersive nonlinear wave equations, like the Korteweg de-Vries and
nonlinear Schr\"odinger equations, this approach goes back to
\textsc{Benjamin}~\cite{tbb72}, \textsc{Bona}~\cite{bona75},
\textsc{Weinstein}~\cite{weinstein} and many others (see the book
by \textsc{Pava}~\cite{pava} for history and many references).
This strategy was developed into a general and
powerful approach for a class of abstract Hamiltonian PDEs
with one or more constraints by \textsc{Grillakis, Shatah, \& Strauss}~\cite{gss87,gss90} (hereafter GSS).
There has been a vast amount of work in this direction since, too much
to review here (e.g.\  see Chapter 5 in \textsc{Kapitula \& Promislow}~\cite{kp13} and references therein).

From the perspective of this paper, a key part of the GSS theory
is the connection between the sign
of the derivative of a scalar-valued function and minimization.
When there is a single constraint set, say the momentum denoted by $I$, and a single
Lagrange multiplier, the speed of the solitary wave denoted by $c$, the condition is
\begin{equation}\label{dIdc-1}
\frac{dI}{dc}>0 \quad\Rightarrow\quad \mbox{solitary wave is a minimizer}.
\end{equation}
This condition is useful since $I(c)$ is a property of the basic state
and so is, in principle, available and easy to calculate. 

On the other hand a central hypothesis in the GSS theory, required for
(\ref{dIdc-1}), is that
the second variation of the functional, in the case of one
constraint, should have at most one negative eigenvalue, one zero eigenvalue,
and the remainder of the spectrum strictly positive.
It is this GSS spectral hypothesis that is most
difficult to satisfy, and indeed may not be satisfied, especially for coupled PDEs.

Another approach is to study the spectral problem associated with
the linearization about solitary waves
while incorporating the Hamiltonian structure.
The seminal paper in this direction is \textsc{Pego \& Weinstein}~\cite{pw92}
(hereafter PW). They looked to retain the
derivative of $I$ in (\ref{dIdc-1}) and find its role in the linear stability
problem, but work around the GSS
spectral hypothesis.  By combining the Evans function, $D(\lambda)$, with
the Hamiltonian structure, and the energy-momentum
characterization of solitary waves, PW proved that it
has the following properties at $\lambda=0$,
\begin{equation}\label{Dpp-1}
D(0)=0\,,\quad D'(0)=0\,,\quad D''(0) = \frac{dI}{dc}\,.
\end{equation}
The derivative of $I$ in (\ref{dIdc-1}) appears in this formula
in a natural way,
but the GSS spectral condition is not used in the proof.
This result is useful as it is straightforward, when the evolution equation
is well-posed, to normalize the Evans function so
that it satisfies $D(\lambda)\to 1$ as $\lambda\to+\infty$ along the real axis.
Hence when $dI/dc<0$ the existence of an unstable stability exponent is assured
by the intermediate value theorem.  The theory was applied to
scalar-valued PDEs such as generalized KdV, BBM equation, and Boussinesq equation, and in all cases the formula (\ref{Dpp-1}) was applicable.

\textsc{Bridges \& Derks}~\cite{bd99,bd01,bd03} extended the PW theory and
showed that, in the presence of symmetry, there may be an additional
factor in the second derivative in (\ref{Dpp-1})
\begin{equation}\label{Dpp-BD}
D''(0) = \chi\,\frac{dI}{dc}\,.
\end{equation}
When $D(\lambda)\to 1$ as $\lambda\to+\infty$ along the real axis
it is the negativity of the full product that gives existence
of an unstable eigenvalue.
The factor $\chi$ is calculated independently of the derivative of $I$
and is not just a scale factor.
An explicit formula was found for the factor $\chi$ but
the presence of symmetry,
other than translation invariance in space, was an essential part of
the proof in \cite{bd99,bd01,bd03}. 
Moreover the
theory relied on the ``system at infinity'' having only one positive and one
negative real (spatial) eigenvalue when $\lambda=0$
(see \S\ref{sec-canonicalform} and \S\ref{sec-Evansfunction}
for the definition of ``spatial eigenvalue'' and ``system
at infinity''). Several examples were given with $\chi$ taking both positive and
negative values, showing that an additional factor is essential in general.

In this paper the assumptions of additional symmetry and one-dimensional
stable manifold in the system at infinity are removed.  A new formula
for the second derivative is found in the form
\begin{equation}\label{Dpp-transversality}
D''(0) = \chi \Pi\,\frac{dI}{dc}\,.
\end{equation}
The factor $\Pi$ is associated with the transversal
intersection of the stable and unstable manifolds which form the solitary wave,
characterized as a homoclinic orbit.  
In the case of a two-dimensional stable and unstable manifold,
the new factor is
\begin{equation}\label{Pi-def-a}
  \Pi = \bO({\bf a}^-,{\bf a}^+)\,,
\end{equation}
where $\ba^-$ and $\ba^+$ are $\xi-$dependent
tangent vectors to the stable and unstable manifolds respectively
and $\bO$ is a symplectic form associated with a
$c-$dependent spatial symplectic structure
(defined in \S\ref{sec-canonicalform}).
When
the dimension of the stable and unstable manifolds is greater than two the
formula (\ref{Pi-def-a}) expands to be the determinant of a matrix of
symplectic forms (see \S\ref{sec-cr} and \cite{treschev,cb15}).
The importance of $\bO({\bf a}^-,{\bf a}^+)$
as a symplectic invariant of homoclinic orbits was discovered by
\textsc{Lazutkin}~\cite{gl01},
and hence we call it the Lazutkin invariant, and its properties and connection with the parity of the Maslov index are proved by \textsc{Chardard \& Bridges}~\cite{cb15}. 

There are a number of hypotheses that go into
the result (\ref{Dpp-transversality})
but the most important are firstly that no symmetry (other than translation
invariance) is assumed, and
secondly the system at infinity is not restricted to one (spatial) eigenvalue with positive real part
in the limit $\lambda\to0$.
In the body of the paper we will restrict
attention to the case where the homoclinic orbit is the intersection between
a two-dimensional stable and two-dimensional unstable manifold,
with the generalization to arbitrary dimension discussed in \S\ref{sec-cr}.

The role of transversality in the Evans function formulation of
the linear stability problem for solitary waves here is new but not that
surprising.  In the case of dissipative PDEs, \textsc{Alexander \& Jones}~\cite{aj94} prove that the first derivative of the Evans function can be characterized in terms
of a coefficient of transversality.  This theory is abstracted to an orientation
index in \S9.4 of \textsc{Kapitula \& Promislow}~\cite{kp13}.
\textsc{Chardard \& Bridges}~\cite{cb15} prove that in
gradient systems, with Hamiltonian steady part, the first derivative of the Evans function can be expressed in terms
of transversality. However, in all previous work
transverality shows up in the first derivative, and here symplecticity
in the time direction comes into play, and it shows up in the second derivative and
includes the additional factors $\chi$ and the momentum function $dI/dc$.

In order to give the result  (\ref{Dpp-transversality})
some generality we need an abstract class of Hamiltonian
PDEs.  By way of comparison, the class of Hamiltonian PDEs in
PW \cite{pw92} is
\begin{equation}\label{pw-hamiltonian}
u_t = \mathcal{J}\nabla_u H(u,u_x)\,,\quad u\in \mathbb{X}\,,
\end{equation}
for some function space $\mathbb{X}$,
where $\mathcal{J}:\mathbb{X}^*\to\mathbb{X}$ is the co-symplectic (or Poisson)
operator, $u$ is scalar-valued, and $H:\mathbb{X}\to\R$ is the Hamiltonian function.
However, reduction of (\ref{pw-hamiltonian}) to a steady problem is an ODE,
\[
c u_\xi = \nabla_u H(u,u_\xi)\,,\quad \xi=x+ct\,.
\]
With modest hypotheses on $\mathcal{J}$ this equation is generated by a
Lagrangian with density $L(u,u_\xi)$.  Invoking a Legendre transform
then brings in a second hidden symplectic structure
and a {\it finite-dimensional}  Hamiltonian system. 
  Denote the second ``spatial symplectic operator'' by ${\bf K}$.  This
  spatial symplectic structure is essential for both defining symplectic
  transversality and for the proof of the formula (\ref{Dpp-BD}).

  It is clear that the interplay between two symplectic structures is an
  essential part of the analysis: the time evolution and
  the energy-momentum characterization of
  the solitary wave use the temporal symplectic structure, whereas
  transversality of the homoclinic orbit representation of the solitary
  wave is defined using the spatial symplectic structure.  The Evans function
  is defined using \emph{both} symplectic structures.
  Hence, introducing a finite-dimensional
  representation of $\mathcal{J}$, and new coordinates,
  leads to a formulation of the Hamiltonian
  PDE in terms of multisymplectic structure (e.g.\ \cite{bd99,bd01};
  canonical structures and history are given
  in \cite{tjb06,bhl10}).
  The canonical form for a multisymplectic Hamiltonian PDE of interest
  here is
\begin{equation}\label{MKS-1}
\bM Z_t + \bK Z_x = \nabla S(Z)\,,\quad Z\in \R^{2n}\,,
\end{equation}
where $\bM$ and $\bK$ are symplectic operators which are taken to be
constant and $S$ is a generalized Hamiltonian function with $\bM$
a finite dimensional representation on the phase space $\R^{2n}$ of the
infinite-dimensional operator $\mathcal{J}$ in (\ref{pw-hamiltonian}).
Steady solutions $Z(x,t)=\Zh(\xi)$, $\xi=x+ct$, are orbits of the
finite-dimensional Hamiltonian system
\begin{equation}\label{MKS-steady}
(\bK+c\bM) \Zh_\xi = \nabla S(\Zh)\,,\quad \Zh\in \R^{2n}\,.
\end{equation}
In this system a solitary wave is represented by a homoclinic orbit.
The theory will be developed for the case $n=2$ which is the lowest
dimension of interest, and limits the proliferation of indices, with
comments on the general $n>2$ case in the concluding remarks.
The abstract form (\ref{MKS-1}) is quite 
satisfactory for the theory and represents a wide range of Hamiltonian PDEs
\cite{tjb06,bd99,bhl10,bd01,bd03}.

However, we go one step further in this paper and introduce an abstract
class of multisymplectic Hamiltonian PDEs.
Given a smooth pseudo-Riemannian
manifold there is a natural form on the total exterior algebra
bundle whose variation produces a coordinate-free version of the left-hand side
of (\ref{MKS-1}).  This construction generalizes the symplectic
structure on the cotangent bundle of a Riemannian manifold in
classical mechanics.
With this strategy we get a coordinate-free formulation as well
as the canonical form (\ref{MKS-1}).
In fact the partial differential operator generated is a Dirac
operator.  It is made nonlinear by adding a gradient on the right-hand side.  We
call the class of PDEs generated on the total exterior algebra
bundle {\it multisymplectic Dirac operators}.
This class of Hamiltonian PDEs includes as special cases the coupled mode
equation which appears in fluid dynamics \cite{cgd11,grimshaw-2000,gc01,gs02}
and optics \cite{splj09,cp06,bpz98},
the massive-Thirring model \cite{ow75}, and a class of coupled nonlinear wave equations.

Solitary waves are
{\it relative equilibria}; that is, solutions of the
Hamiltonian PDE that are equilibria in a moving frame of reference.  Hence,
in looking for a motivation for the formula (\ref{Dpp-transversality}), we
first consider the spectral problem for relative equilibria of Hamiltonian ODEs and establish that
\begin{equation}\label{morse-exponent}
D''(0) = (-1)^{\mbox{\footnotesize Morse}}\, \frac{dI}{dc}\,,
\end{equation}
where the exponent
is the Morse index of the constrained critical point problem (the number of strictly negative eigenvalues of the
constrained second variation). In the context of ODEs the proof of (\ref{morse-exponent}) uses elementary linear algebra.  This result ties in with the GSS theory because it contains a weak form of the GSS spectral condition.
It is weak in that only the parity of the number of negative
eigenvalues enters the formula.

On the other hand,
solitary waves in the energy-momentum construction, may or may not
have a well-defined Morse index.  So the result (\ref{morse-exponent})
is not expected to generalize to solitary waves.  However,
using Theorem 10.1 in \cite{cb15} we can go one step further and relate the new
characterization of $\Pi$ to the Maslov index of the solitary wave
\begin{equation}\label{maslov-exponent}
{\rm sign}(\Pi) = (-1)^{\mbox{\footnotesize Maslov}}\,,
\end{equation}
where in this case the Maslov index of the solitary wave
is defined using the Souriau characterization
(cf.\ \S9 of \cite{cb15}).
Solitary waves, with exponential decay at
infinity, always have a well-defined Maslov index, but may not have
a well-defined Morse index.  The Maslov index will not feature in this paper,
as the emphasis here is on transversality, but some elaboration of
(\ref{maslov-exponent}) is given in \S\ref{sec-deriv-evans}.

An outline of the paper is as follows. In Section \ref{Dprimeprime-odes}
the special case (\ref{morse-exponent})
of the derivative formula is proved for relative
equilibria of Hamiltonian ODEs.  In Section \ref{sec-Diracoperators} an
abstract class of multisymplectic Hamiltonian PDEs is introduced.  
Section \ref{sec-canonicalform} is the starting point for proving the
main results on stability of solitary waves.  Here the abstract
class of solitary waves is introduced as well as the properties of the
linearization about these waves.  Section \ref{sec-Evansfunction}
constructs the Evans function and develops the interplay with
symplecticity.
Section \ref{sec-deriv-evans} proves the main result on $D''(0)$ confirming
(\ref{Dpp-transversality}).  Section \ref{sec_coupledwave-eqn} gives
an example where all the details are worked out explicitly.
In \S\ref{sec-complex-eigs} the difficulties with the
case of complex $\mu-$eigenvalues are discussed.
Finally in the concluding remarks Section \ref{sec-cr} some generalizations
are pointed out.

\section{Instability of relative equilibria of ODEs}
\label{Dprimeprime-odes}
\setcounter{equation}{0}

A solitary wave solution of a Hamiltonian PDE
is a relative equilibrium in the following sense.  Focussing
on the form (\ref{pw-hamiltonian}) for description, suppose
the Hamiltonian function and symplectic structure do not depend
explicitly on the spatial coordinate, $x$.  Then $u(x+s,t)$
is a solution for any
$s$ whenever $u(x,t)$ is, and we say that the Hamiltonian PDE is equivariant
with respect to the group $G=\R$, the group of real numbers. When
$s = ct$ with $c$ a constant, and $u$ is otherwise independent of
$t$, the solution is a relative equilibrium of the form
$u(x,t):=\widehat u(x+ct)$; that is, it is an equilibrium when viewed from
a frame of reference moving at constant speed along the group $\R$.
Noether theory then gives the
existence of an invariant associated with the translation
symmetry that is called
the momentum, here denoted by $I$.  It is a functional and depends on
$u$, but when $I$ is evaluated on a family of relative equilibria it becomes a
function of $c$ only, and it is this function that appears in the
derivative formula (\ref{Dpp-BD}).

The main aims of
this section are to show, in the simplest possible setting,
how the product structure
of $D''(0)$ arises, and to show how the geometry of relative
equilibria (RE) enters the analysis.  Although this section is
restricted to ODEs, the abstract structure is the same as that of
the solitary wave stability problem with the Evans function replaced
here by an elementary characteristic function.
The group is simplified to the one-parameter compact group $S^1$.

Consider the finite-dimensional Hamiltonian system on $M:=\R^{2n}$:
\begin{equation}\label{a.1}
\bM Z_t = \nabla H(Z)\,,\quad\quad Z\in M\,,\quad\quad
\bM =
\begin{pmatrix}
  {\bf 0}\phantom{_n} & -{\bf I}_n\\ {\bf I}_n & \hfill{\bf 0}\phantom{_n}
\end{pmatrix}\,,
\end{equation}
where $H: M\to\R$ is a given smooth function.
The system (\ref{a.1}) is assumed to be equivariant with
respect to $S^1$: that is, there is an orthogonal action,
$\Gt$, of
$\Tq$ on $M$ satisfying
\begin{equation}\label{a.2}
  \Gt^T\bM \Gt = \bM\quad
     {\rm and}\quad H(\Gt Z)=H(Z)\quad \forall\,\theta\in\Tq\,.
\end{equation}
In this section $\langle\cdot,\cdot\rangle$ is an inner product on
$\R^{2n}$.
\vspace{.15cm}

\noindent{\bf Proposition 2.1.} {\it
  Let \phantom{\Big|}$\g(Z) = \frac{d\ }{d\theta}\Gt Z\big|_{\theta=0}$ for
  $Z\in M$, then there exists (symplectic Noether theory)
  a functional $I: M\to\R$ satisfying}
$\bM \g(Z)=\nabla I(Z)$.

\noindent{\bf Proof.} In the Hamiltonian setting, Noether theory states
that contraction of the symplectic form by the generator
$\g(Z)$ is the gradient of a function (e.g.\ \S2.7 of \cite{m92}).  In
this case, the existence of $I$ follows from a linear algebra argument.
Since the action of the group is orthogonal, the tangent vector $\g(Z)$ can be expressed as $\g(Z)={\bf S}Z$ with
${\bf S}^T=-{\bf S}$.  Now ${\bf M}$ and ${\bf S}$ are two skew-symmetric
matrices that commute (proved by differentiating
the first of (\ref{a.2}) with respect to the group parameter)
and so their product is symmetric.  Define
$I(Z) = \fr \langle {\bf MS}Z,Z\rangle$, then
\begin{equation}\label{a.3}
  {\bf M}\g(Z) = {\bf MS}Z = \nabla I(Z)\,,
\end{equation}
completing the proof. $\hfill\blacksquare$
\vspace{.15cm}

Now, suppose there exists a family of RE of the above system.
An RE is a solution of (\ref{a.1}) that is aligned with the group
orbit and has constant speed (cf.\ Chapter 4 of \cite{m92}); that is,
a solution of the form
\begin{equation}\label{A.5}
  \widehat Z(t)=\Gtt{\cal U}(c)\quad {\rm with}\quad \theta(t)= c\,t+\theta^o\,.
\end{equation}
Substitution of this form into the governing equation (\ref{a.1})
gives
\[
  \begin{array}{rcl}
    0 &=& {\bf M}\Zh_t - \nabla H(\Zh)\\[2mm]
    &=& \dot\theta{\bf M}\g(\Gt\mathcal{U}) - \nabla H(\Gt\mathcal{U})\\[2mm]
&=& \Gt\big( c \nabla I(\mathcal{U}) - \nabla H(\mathcal U)\big)\,,
  \end{array}
  \]
  and so $\mathcal{U}(c)$ and $c$ are defined by the equations
\begin{equation}\label{A.6}
\nabla H({\cal U})= c\nabla I({\cal U})\quad
         {\rm and}\quad I({\cal U})=I_0\,,
\end{equation}
that is; $\mathcal{U}(c)\in M$ can be characterised as a critical point
of $H$ on level sets, $I=I_0$, of the functional $I$,
with $c$ as a Lagrange multiplier.  This family of RE is assumed to
exist for all $c\in\mathscr{C}$ where $\mathscr{C}$ is
some open subset of $\R$.
A RE in the family is said to be non-degenerate on $\mathscr{C}$ when
$\frac{dI}{dc}\neq 0$ where $I$ is evaluated on the family $\mathcal{U}(c)$.
This family of RE is a finite-dimensional
analogue of a travelling wave.

\subsection{Linear stability of the family of RE}
\label{subsec-lin-stab-RE}

The linear stability equation for the family of
RE is formulated by linearizing (\ref{a.1}) 
about (\ref{A.5})
 \begin{equation}\label{A.8}
   \bM Z_t = D^2H(\widehat Z(t))\,Z\,,
 \end{equation}
 with $\widehat Z(t)$ defined in (\ref{A.5})-(\ref{A.6}). The
 group can be factored out by letting
 $Z(t)=\Gtt W(t)$, then
   \[
 \begin{array}{rcl}
   0 &=& \bM Z_t - D^2H(\widehat Z(t))Z\\[2mm]
   &=& \dot\theta \bM \g(\Gtt W) + \bM \Gtt W_t - D^2H(\Gtt \mathcal{U})Z \\[2mm]
   &=&\displaystyle
   \Gtt\big( c \bM \g(W) + \bM W_t - \Gtt^T D^2H(\Gtt \mathcal{U})\Gtt W\big) \\[2mm]
   &=&\displaystyle
   \Gtt\big( c D^2I(\mathcal{U})W + \bM W_t - D^2H(\mathcal{U}) W\big)\,,
 \end{array}
 \]
 using the first of (\ref{a.2}), noting that $\nabla I(W)=D^2I(\mathcal U)W$,
and using the identity
 \[
 D^2H(\Gt{\cal U})=\Gt D^2H({\cal U})\Gt^T\,,
 \]
 obtained by
 differentiating the second of (\ref{a.2}) with respect to $\theta$,
 twice.
 Hence, the linearized equation (\ref{A.8}) is equivalent to
 the constant coefficient ODE:
\begin{equation}\label{A.9}
  \bM W_t = \LUw\,W\,,\quad{\rm where}\quad
  \LUw=D^2H({\cal U}(c))-c D^2I({\cal U}(c))\,,
\end{equation}
with associated spectral equation $\LUw W=\lambda\bM W$.  Let
\begin{equation}\label{A.10}
  D(\lambda)={\rm det}[\LUw-\lambda\bM ]\,,
\end{equation}
then we have the following sufficient condition for instability:
{\it if there exists a $\lambda\in\C$ with ${\rm Re}(\lambda)>0$ and
$D(\lambda)=0$, the
RE (\ref{A.5}) is linearly (spectrally)
unstable.\/}  $D(\lambda)$ is a finite-dimensional analogue of
the Evans function.

The operator ${\bf L}$ has a zero eigenvalue
with eigenvector $\g(\mathcal U)$.  To show this, act on (\ref{A.6}) with
$\Gt$, use equivariance, differentiate with respect to $\theta$, and set $\theta=0$,
\begin{equation}\label{kernel-of-L}
0 = \frac{d\ }{d\theta}\Big( \nabla H(\Gt\mathcal{U}) - c \nabla(\Gt\mathcal{U}) \Big)\Big|_{\theta=0}\quad\Rightarrow\quad {\bf L}\g(\mathcal{U})=0\,.
\end{equation}
It is assumed that the zero eigenvalue
of ${\bf L}$ is simple for $c\in\mathscr{C}$.
Differentiating the first of (\ref{A.6}) with respect to $c$
and using (\ref{a.3}) gives the generalized eigenvector
  \[
    {\bf L}\mathcal{U}_c={\bf M}\g(\mathcal{U})\,.
    \]
    Normalize the length of these two vectors and define
  \[
  \zeta_1 = \frac{\g(\mathcal{U})}{\|\g(\mathcal{U})\|}\qand
  \zeta_2 = \frac{\mathcal{U}_c}{\|\g(\mathcal{U})\|}\,.
  \]
  Then they generate the Jordan chain of length two
  \begin{equation}\label{zeta12-def}
    {\bf L}\zeta_1=0 \qand {\bf L}\zeta_2 = {\bf M}\zeta_1\,.
  \end{equation}
  With the simple zero eigenvalue (\ref{kernel-of-L}) and $dI/dc\neq0$,
  the Jordan chain
  has length exactly two.  To have length three would require solvability of
  \begin{equation}\label{jc-length-3}
    {\bf L}\zeta_3 = {\bf M}\zeta_2\,.
  \end{equation}
    But solvability of this equation requires
  \begin{equation}\label{jc-length-3a}
    0 = \langle\zeta_1,{\bf M}\zeta_2\rangle = -\frac{1}{\|\g(\mathcal{U})\|^2}\langle {\bf M}\g(\mathcal{U}),\mathcal{U}_c\rangle
=-\frac{1}{\|\g(\mathcal{U})\|^2}\langle \nabla I(\mathcal{U}),\mathcal{U}_c\rangle    =
    -\frac{1}{\|\g(\mathcal{U})\|^2} \frac{dI}{dc}\,.
  \end{equation}
    Hence with $dI/dc\neq0$ this equation is not solvable.

  \subsection{Derivatives of the characteristic function}
  \label{subsec-deriv-of-D}
  
    We are now in a position to prove the following finite-dimensional
    analogue of (\ref{Dpp-1}) with second derivative (\ref{Dpp-transversality}).
    \vspace{.15cm}
    
    \noindent{\bf Theorem 2.2.} {\it Suppose a smooth family of RE exists
      for all $c\in\mathscr{C}$, with $dI/dc\neq0$, and suppose ${\bf L}$ has
      a simple zero eigenvalue.  Then the
      characteristic function $D(\lambda)$ has
      the following derivatives at the origin}
    \[
    D(0)=0\,,\quad D'(0)=0\,,\quad D''(0) = \mu({\bf L}) \frac{dI}{dc}\,,
    \]
    {\it where $\mu({\bf L})$ is the product of the nonzero eigenvalues of
      ${\bf L}$.}
    \vspace{.15cm}

    \noindent{\bf Remark.} The sign of $\mu({\bf L})$ is the parity of
    the number of negative eigenvalues so with a suitable scaling of
    $D(\lambda)$ an equivalent formula for $D''(0)$ is
    \[
    D''(0) = (-1)^{\mbox{\footnotesize Morse}}\,\frac{dI}{dc}\,,
    \]
    where {\small Morse} is the Morse index of
    ${\bf L}$,  confirming (\ref{morse-exponent}) in the introduction.

    \noindent{\bf Proof.} Since ${\bf L}$ has a simple
    zero eigenvalue
    \[
    D(0) = {\rm det}[{\bf L}]=0\,.
    \]
    Differentiate $D(\lambda)$ in (\ref{A.10}) using the formula
  for the derivative of a determinant
  \begin{equation}\label{Dp-formula}
    D'(\lambda) = -{\rm Tr}\Big(({\bf L}-\lambda\bM)^\#\bM\Big)
    \quad\Rightarrow\quad D'(0) = -{\rm Tr}({\bf L}^\#\bM)\,,
  \end{equation}
  where ${\bf L}^\#$ is the adjugate of ${\bf L}$.
  When ${\bf L}$ has only one zero eigenvalue with unit length
  eigenvector $\zeta_1$ then ${\bf L}^\#$ is the rank one matrix
  \begin{equation}\label{Pi-def-ODEs}
    {\bf L}^\# = \mu({\bf L}) \zeta_1 \zeta_1^T\quad\mbox{with}\quad 
  \mu({\bf L})=\prod_{j=2}^{2n}\mu_j\,.
  \end{equation}
  $\mu({\bf L})$ is
  the product of the nonzero eigenvalues, $\mu_j$, of ${\bf L}$
  (taking the zero eigenvalue to be $\mu_1$).
  This formula is stated and proved as Theorem 3 on page 48 of
  \textsc{Magnus \& Neudecker}~\cite{mn-book}.  Substitute
  ${\bf L}^\#$ into (\ref{Dp-formula}),
  \[
  D'(0) = -{\rm Tr}({\bf L}^\#{\bf M}) = -\mu({\bf L}){\rm Tr}(
  \zeta_1\zeta_1^T{\bf M}) = -\mu({\bf L})\langle\zeta_1,{\bf M}\zeta_1\rangle =0\,,
  \]
  since ${\bf M}$ is skew symmetric.  For $D''(0)$, differentiate
  $D'(\lambda)$ in (\ref{Dp-formula})
  \[
  D''(0) = -{\rm Tr}\left( \frac{d\ }{d\lambda}({\bf L}-\lambda\bM)^\#\bigg|_{\lambda=0}\bM\right)\,.
  \]
  The adjugate is defined by
  \begin{equation}\label{adjugate-def}
({\bf L}-\lambda\bM)({\bf L}-\lambda\bM)^\#=  ({\bf L}-\lambda\bM)^\#({\bf L}-\lambda\bM) = D(\lambda){\bf I}\,,
  \end{equation}
  where ${\bf I}$ is the identity on $\R^{2n}$. 
  Now differentiate (\ref{adjugate-def}) with respect to $\lambda$,
  set $\lambda$ to zero, and define
  \[
  \dot{\bf L} := \frac{d\ }{d\lambda}({\bf L}-\lambda\bM)^\#\bigg|_{\lambda=0}\,.
  \]
This gives the following equations for $\dot{\bf L}$
  \[
    {\bf L}\dot{\bf L} = \bM{\bf L}^\#\qand \dot{\bf L}^T=-\dot{\bf L}\,,
    \]
    with skew-symmetry following from commutivity in (\ref{adjugate-def}).
    Combining (\ref{zeta12-def}), (\ref{Pi-def-ODEs}) and skew-symmetry
    of $\dot{\bf L}$ gives
  \[
  \dot{\bf L} = \mu({\bf L})(\zeta_2\zeta_1^T-\zeta_1\zeta_2^T) \,.
  \]
  Substitute into $D''(0)$
  \[
  \begin{array}{rcl}
    D''(0)
    &=& -\mu({\bf L}){\rm Tr}\Big( (\zeta_2\zeta_1^T-\zeta_1\zeta_2^T)\bM\Big)\\[2mm]
    &=& -\mu({\bf L})\Big( \langle \zeta_1,\bM\zeta_2\rangle -
    \langle\zeta_2,\bM\zeta_1\rangle\big)\\[2mm]
    &=& 2\mu({\bf L})\langle\zeta_2,\bM\zeta_1\rangle\\[4mm]
    &=& \displaystyle
    2\frac{\mu({\bf L})}{\|\g(\mathcal{U})\|^2} \frac{dI}{dc}\,,
  \end{array}
  \]
  with the last expression following from (\ref{a.3}) and
 \[
  \frac{dI}{dc} = \langle \nabla I(\mathcal{U}),\mathcal{U}_c\rangle
  = \langle \bM\g(\mathcal{U}),\mathcal{U}_c\rangle=
  \|\g(\mathcal{U})\|^2 \langle\bM\zeta_1,\zeta_2\rangle\,.
  \]
  Scaling $D(\lambda)$ by a positive constant then completes the proof.
  $\hfill\blacksquare$
  \vspace{.15cm}

\noindent{\bf Remark.} Since (\ref{A.10}) is ``the Evans function'' in
this example, scaling $D(\lambda)$ by a positive
constant should be interpreted as a trivial
re-definition of it to
\[
\widetilde D(\lambda) = \fr \|\g(\mathcal{U}\|^2\, D(\lambda)\,.
\]
and it is this new Evans function $\widetilde D(\lambda)$ that has
$\widetilde D''(0) = \mu({\bf L})dI/dc$.
\vspace{.15cm}

  \noindent{\bf Corollary 2.3.} {\it When} $(-1)^{\mbox{\footnotesize Morse}}\,\frac{dI}{dc}<0$ {\it the family of RE has an unstable eigenvalue.}
  \vspace{.15cm}
  
  \noindent{\bf Proof.} The condition assures that $D(\lambda)$ is negative
  for $\lambda$ near zero.  For large and real $\lambda$ the characteristic
  function has the asymptotic form
    \[
  D(\lambda) = (-1)^{2n}{\rm det}(\bM)\lambda^{2n} + \cdots\,,
  \]
  and so $D(\lambda)>0$ for $\lambda$ real, positive, and sufficiently large.
   By the intermediate value theorem $D(\lambda)$
  has at least one positive real root.$\hfill\blacksquare$
  \vspace{.15cm}

     Theorem 2.2 connects $dI/dc$ to the spectral problem and
     is a finite dimensional version of the derivative formula
     (\ref{Dpp-transversality}).
     The connection between $dI/dc$ and critical point type for relative
     equilibria appears in the literature from various perspectives
    (e.g. \textsc{Maddocks \& Sachs}~\cite{ms95} and references therein),
     and is an elementary example of the critical point theory in
     infinite-dimensional spaces in GSS \cite{gss87,gss90}.
        
     The plan is to extrapolate Theorem 2.2 to the context of solitary waves.
     Before proceeding with that proof, the
    next section shows that the canonical form for Hamiltonian PDEs
    (\ref{MKS-1}) is not only a useful representation of well known PDEs,
    it is also a universal class.

\section{A class of multisymplectic Hamiltonian PDEs}
\setcounter{equation}{0}
\label{sec-Diracoperators}

The class of multisymplectic Hamiltonian PDEs (\ref{MKS-1})
is a natural starting point for the theory.
Indeed all the theory in \cite{bd99,bd01,bd03} is based on this
class of PDEs.  In this section it is shown that this class of PDEs
can be obtained naturally and coordinate free
from a pseudo-Riemannian manifold, showing that
the class is universal in addition to being practical.

This approach is a generalization of the cotangent
bundle of a manifold as a natural and coordinate free generator of
symplectic structure.  Let $M$ be a smooth manifold,
which for simplicity is taken to be $\R^n$.
Let $(q_1,\ldots,q_n,p_1,\ldots,p_n)$ be local coordinates for
$T^*M\cong\R^{2n}$. The cotangent bundle hosts a canonical one
form $\bp\cdot\rd\bq$ with associated functional
\begin{equation}\label{pdq-integral}
\int_{t_1}^{t_2} \bp\cdot \bq_t\,\rd t\,.
\end{equation}
The first variation of this functional
\[
0=\delta\int_{t_1}^{t_2} \bp\cdot \bq_t\,\rd t=
\int_{t_1}^{t_2} \big(\delta\bp\cdot \bq_t+
 \bp\cdot \delta\bq_t\big)\,\rd t\,,
 \]
 with fixed endpoints on variations, $\delta{\bf q}(t_1)=\delta{\bf q}(t_2)=0$,
generates the operator
\[
  {\bf J}\frac{d\ }{dt}\quad\mbox{with}\quad {\bf J}= \left[ \begin{matrix}
      {\bf 0} & -{\bf I}_n\\
      {\bf I}_n & {\bf 0} \end{matrix}\right] \,.
  \]
  Two observations about this operator are of interest:
  firstly, it generates a Hamiltonian
  system by introducing the gradient
  of $H(\bq,\bp)$, a given smooth function,
  \[
    {\bf J}\frac{d\ }{dt}\begin{pmatrix} \bq\\ \bp\end{pmatrix} = \nabla H\,,
      \]
      and secondly it is a one-dimensional ``Dirac operator''
      \[
        {\bf J}\frac{d\ }{dt}\circ{\bf J}\frac{d\ }{dt} = - {\bf I}_{2n}\otimes
        \frac{d^2\ }{dt^2}\,.
        \]
        The strategy here is to generalize this construction to generate
        abstract multisymplectic Hamiltonian PDEs.  The main difference in
        the PDE case is
        that the manifold $M$ is the base manifold representing space-time,
        and the fiber is built on the total exterior algebra bundle rather
        than just the cotangent bundle.

        The starting point is a smooth
        pseudo-Riemannian manifold $M$ (which here is taken to be the
        flat space $M=\R^{q,p}$), with constant signature metric
        which can be represented by
\begin{equation}\label{metric-def}
\lbk{\bf u},{\bf v}\rbk_1 := \big\langle \Rg{\bf u},{\bf v}\big\rangle\,,
\end{equation}
with $\langle\cdot,\cdot\rangle$ a Euclidean inner product, and
\begin{equation}\label{G-def}
  \Rg = {\rm diag}(\underbrace{1,\cdots,1}_{\mbox{\footnotesize q times}},
  \underbrace{-1,\cdots,-1}_{\mbox{\footnotesize p times}})\,.
\end{equation}
The metric $\lbk\cdot,\cdot\rbk_1$
  induces a metric on each of the spaces
  $\bigwedge^k(T_x^*M)$.  The induced metrics are denoted by
  \begin{equation}\label{metric-def-k}
    \lprod{\bf u}^{(k)},{\bf v}^{(k)}\rprod_k\ \mbox{for}\
          {\bf u}^{(k)},{\bf v}^{(k)}\in\Dform^k(M)\,,
  \end{equation}
  where
  $\Dform^k(M)$ is the space of mappings
  from $M$ into $\bigwedge^k(T_x^*M)$. Concatenating these spaces gives
  the \emph{total exterior algebra} (TEA) \emph{bundle} denoted by
  $\Dform(M) := \bigcup_{k=0}^n\Dform^k(M)$.  Let
  $r={\rm dim}\left(\bigwedge^k(T_x^*M)\right)$,
  then we have the Euclidean space
  representation of the lifted metric
  \begin{equation}\label{Rr-def}
  \lprod{\bf u}^{(k)},{\bf v}^{(k)}\rprod_k = \langle
  \Rg^{(r)}{\bf u}^{{(k)}^r},{\bf v}^{{(k)}^r}\rangle\,,
  \end{equation}
  where ${\bf u}^{{(k)}^r}$ is the lift to $\R^r$ of ${\bf u}^{(k)}$.

  The generalization of ${\bf p}\cdot\rd{\bf q}$ on the cotangent bundle
  is the following form
  on the total exterior algebra bundle
  \begin{equation}\label{Theta-def}
  \Theta(Z) = \sum_{k=1}^m \mathbf{u}^{(k)}\wedge\bigstar\rd\mathbf{u}^{(k-1)}\,,
  \end{equation}
  where $\bigstar$ is Hodge star and $\rd$ is an
  exterior derivative.
  See \cite{tjb06} for the origin of this form in
  the case of a positive
  definite metric, which generates an elliptic partial differential
  operator (PDO), and see \cite{burchell-thesis} for the
  case of general indefinite metric, which
  generates a hyperbolic PDO.
  
  Define the Lagrangian density
          \begin{equation}\label{L-def-sec2}
          \mathcal{L}(Z) = \Theta(Z) - S(Z)\vol\,,
          \end{equation}
          where $S:\Dform(M)\to \R$ is a given generalized Hamilton function.
Taking the first variation of the functional
          \begin{equation}\label{L-def-sec2-1}
                      \delta\int \Theta(Z)-S(Z)\vol=0\,,
          \end{equation}
          generates a nonlinear Dirac operator
          \begin{equation}\label{JZS-abstract}
          \J_\partial Z = \Rg^{(r)}\nabla S(Z)\,.
          \end{equation}
          When written out in coordinates, and pre-multiplying by
          $\Rg^{(r)}$, the PDE becomes
          \begin{equation}\label{JZS-coordinates}
          \sum_{j=1}^n \Rg^{(r)}\J_j \partial_{x_j} Z = \nabla S(Z)\,.
          \end{equation}
          This PDE is now in standard form for a multisymplectic
          Hamiltonian PDE \cite{bhl10}.  Indeed it is a new class of
          multisymplectic Hamiltonian PDEs. It is not difficult to show
          that each $\Rg^{(r)}\J_j$ is symplectic (skew-symmetric and
          non-degenerate), acting on a space of dimension $2^n$,
          and so $n-$independent symplectic structures
          are generated.

\subsection{Multisymplectic Dirac operator based on $M=\R^{1,1}$}
\label{subsec-PDO_1-1}

The case of $M=\R^{1,1}$ with metric tensor $\Rg={\rm diag}(1,-1)$ is
the case of interest in this paper.  Take coordinates $(t,x)$ and volume form
$\vol=\rd t\wedge\rd x$.  Then differential forms in the TEA bundle
are of the form $Z=(\phi,\bu,v)$ with $\phi$ a
scalar-valued function,
\[
\bu = u_1\rd t + u_2\rd x \qand
v := v\rd t\wedge \rd x\,,
\]
where, to simplify notation, $v$ is both a form and a coordinate.
The PDO in this case acting on $Z\in\Dform(M)$ is
\[
\J_\partial Z =
\left(\!\begin{array}{c c c}
0 & \delta & 0 \\
\xd & 0 & \delta \\
0 & \xd & 0
\end{array}\!\right)\!\!
\left(\!\begin{array}{c}
\phi \\
\bu \\
v
\end{array}\!\right)\!.
\]
and it can be expressed in coordinates as
\[
\J_\partial = \J_1\delt + \J_2\delx
\]
with
\begin{equation}\label{J1J2-def_Cl11}
\J_1 =
\left(\!\begin{array}{c c c c}
0 & -1 & 0 & 0 \\
1 & 0 & 0 & 0 \\ 
0 & 0 & 0 & -1 \\
0 & 0 & 1 & 0
\end{array}\!\right)\!,
\qquad \J_2 =
\left(\!\begin{array}{c c c c}
0 & 0 & 1 & 0 \\
0 & 0 & 0 & -1 \\ 
1 & 0 & 0 & 0 \\
0 & -1 & 0 & 0
\end{array}\!\right)\,.
\end{equation}
The pair $\{\J_1,\J_2\}$ generates
the Clifford algebra $\mathscr{C}\ell_{1,1}$,
\[
\J_i\J_j+\J_j\J_i = - 2\Rg_{ij}{\bf I}_4\,.
\]
The Dirac property and the connection with the d'Alembertian is
\[
\J_\partial\circ\J_\partial = \left(\J_1\delt + \J_2\delx\right)^2 =
\J_1^2\partial_{tt} + (\J_1\J_2+\J_2\J_1)\partial_{tx} + \J_2^2\partial_{xx} =
-(\partial_{tt} - \partial_{xx}) \otimes{\bf I}_4\,.
\]
Introducing a scalar-valued function $S:\Dform(M)\to\R$, a nonlinear
Dirac equation is generated in the canonical form
\begin{equation} \label{system2}
  \J_\partial Z = \Rg^{(4)}\nabla S(Z)\,,\quad Z\in\Dform(M)\,.
\end{equation}
The induced metric in this case is
\begin{equation}\label{Rg-4-def}
\Rg^{(4)} = \mathrm{diag}(1,1,-1,-1) = \underbrace{[+1]}_{\wedge^0}\oplus\underbrace{\left[\begin{matrix}
    1 & 0 \\ 0 & -1 \end{matrix}\right]}_{\wedge^1}\oplus \underbrace{[-1]}_{\wedge^2}\,.
\end{equation}
This form follows by constructing the induced metric on each of the vector
spaces $\wedge^j$ and then concatenating.
The space $\bigwedge(\R^{1,1})$ is isomorphic to $\R^{2,2}$ with metric
$\langle\Rg^{(4)}\cdot,\cdot\rangle$.

The operator $\J_1$ is skew-symmetric and $\J_2$ is symmetric and
they are both invertible.  The induced skew symmetric operators are
\begin{equation}\label{J1J2-def_skew-sym_R11}
\bM :=\Rg^{(4)}\J_1 =
\left(\!\begin{array}{c c c c}
0 & -1 & 0 & 0 \\
1 & 0 & 0 & 0 \\ 
0 & 0 & 0 & 1 \\
0 & 0 & -1 & 0
\end{array}\!\right)\!,
\qquad \bK :=\Rg^{(4)}\J_2 =
\left(\!\begin{array}{c c c c}
0 & 0 & 1 & 0 \\
0 & 0 & 0 & -1 \\ 
-1 & 0 & 0 & 0 \\
0 & 1 & 0 & 0
\end{array}\!\right)\,.
\end{equation}
The system
\begin{equation}\label{R11_JS-eqn}
\bM Z_t + \bK Z_x = \nabla S(Z) \,,\quad Z\in\Dform(M)\,,
\end{equation}
is then in standard form for a multisymplectic Hamiltonian PDE,
and the two operators $\bM$ and $\bK$ define independent
symplectic vector spaces.
 
  When $\bM$ is invertible as above, the PDE can be written in evolution form
  \[
  Z_t + {\bf C}Z_x = \bM^{-1}\nabla S(Z)\,,\quad Z(x,0) = Z_0(x)\,,
  \]
  with ${\bf C}=\bM^{-1}\bK$.  When $[\bM,\bK]=0$, as is
  the case for $\bM$ and $\bK$ in (\ref{J1J2-def_skew-sym_R11}),
  the matrix ${\bf C}$
  is symmetric and the operator $Z_t+{\bf C}Z_x$ is hyperbolic.
  There are a range of results in the literature on existence and
  well-posedness of equations in this form in general, and
  Dirac equations in particular (e.g.\ \textsc{Pelinovsky}~\cite{dp11}
  and references therein).

\subsection{The coupled mode equation}
\label{subsec-coupledmodeeqn}

The \emph{coupled-mode equation} (CME) which appears in fluid mechanics
\cite{cgd11,grimshaw-2000,gc01,gs02} and optics \cite{splj09,cp06,bpz98} can
be characterised as a multisymplectic Dirac operator on
$\Dform(\R^{1,1})$ in the form
(\ref{system2}). In the literature, the CME is represented
in complex-amplitude form
\begin{equation}\label{cme-2}
\begin{array}{rcl}
\ri (A_t + A_x) + \alpha B + \tau|A|^2A + \nu|B|^2A + \mu
B^2\overline A &=& 0 \\[2mm]
\ri (B_t - B_x) + \alpha A + \tau|B|^2B + \nu|B|^2B + \mu
A^2\overline B &=& 0 \,.
\end{array}
\end{equation}
In this equation the coefficients $\alpha$, $\tau$, $\nu$ and $\mu$
are real-valued and $A(x,t)$ and $B(x,t)$ are complex valued functions.
Introduce coordinates
$(\phi,{\bf u},v)$ in
$\mathcal{A}^0\times \mathcal{A}^1\times\mathcal{A}^2$ and to link more closely
with the CME coordinates, take
\[
  {\bf w}= (w_1,w_2) := \big( \phi , u_1\big)\qand
  {\bf v}= (v_1,v_2) := \big( u_2 , v\big)\,.
  \]
  The system (\ref{cme-2}) is transformed using
\[
\begin{array}{rcl}
  A:= A_1 +\ri A_2 &=& w_1 - v_2 + \ri ( w_2 - v_1)\\[2mm]
  B:= B_1 +\ri B_2 &=& w_1 + v_2 + \ri ( w_2 + v_1)\,.
\end{array}
\]
In these coordinates the CME becomes
\begin{equation}\label{cme-3}
\left[\begin{matrix}
0 & -1 & 0 & 0 \\ 1 & 0 & 0 & 0 
\\ 0 & 0 & 0 & -1 \\ 0 & 0 & 1 & 0 \end{matrix}\right]
\begin{pmatrix} w_1\\ w_2\\ v_1\\ v_2\end{pmatrix}_t +
\left[\begin{matrix}
0 & 0 & 1 & 0 \\ 0 & 0 & 0 & -1 
\\ 1 & 0 & 0 & 0 \\ 0 & -1 & 0 & 0 \end{matrix}\right]
\begin{pmatrix} w_1\\ w_2\\ v_1\\ v_2\end{pmatrix}_x =
\Rg^{(4)}\left[\begin{matrix}
\partial S/\partial w_1\\ \partial S/\partial w_2\\
\partial S/\partial v_1\\ \partial S/\partial v_2\end{matrix}\right]\,,
\end{equation}
or, with $Z=(\bw,\bv)$ now identified with $\R^4$
\begin{equation}\label{J1J2-CME}
  \J_1 Z_t + \J_2 Z_x = \Rg^{(4)}\nabla S(Z)\,,
\end{equation}
using (\ref{J1J2-def_Cl11}) with $\Rg^{(4)}=
{\rm diag}(1,1,-1,-1)$.

A special case of (\ref{cme-2}) arises in optics with
$\tau=\gamma$, $\nu=2\gamma$, and $\mu=0$.  It is
the one-dimensional model that rules nonlinear wave propagation
around a forbidden frequency band gap (cf.\ \textsc{Sugny et al.}~\cite{splj09}).  An even more special case is the massive Thirring model (MTM) where
$\tau=\mu=0$,
\begin{equation}\label{mtm}
\begin{array}{rcl}
\ri (A_t + A_x) + \alpha B + \nu|B|^2A  &=& 0 \\[2mm]
\ri (B_t - B_x) + \alpha A + \nu|A|^2B &=& 0 \,.
\end{array}
\end{equation}
The transformed system for MTM is (\ref{cme-3}) with
\begin{equation}\label{S-mtm}
  S(Z) = -\fr\alpha({\bf w}\cdot{\bf w}-\bv\cdot\bv) -
  \frr\nu(\bw\cdot\bw+\bv\cdot\bv)^2
  +\nu(w_1v_2+w_2v_1)^2\,.
\end{equation}

\subsection{Coupled second-order nonlinear wave equations}
\label{subsec-coupled-nlws}

The pair of coupled second order nonlinear wave equations
\begin{equation}\label{nlws-2}
\phi_{tt} - \phi_{xx} + V_\phi =0 \qand
v_{tt}-v_{xx} - V_v = 0\,,
\end{equation}
where $V(\phi,v)$ is a given smooth function, can also be
transformed to the canonical form (\ref{R11_JS-eqn}).  Introduce
new coordinates $(\phi,u_1,u_2,v)$ via
\[
u_1 = \phi_t - v_x \qand u_2 = \phi_x-v_t\,.
\]
Then with
\begin{equation}\label{S-nlw-def}
S(Z) = \fr (u_1^2-u_2^2) + V(\phi,v)\,,
\end{equation}
the coupled equations (\ref{nlws-2}) are represented by
(\ref{R11_JS-eqn}).  An analysis of the Evans function for an
example in this class of nonlinear wave equations in given in
\S\ref{sec_coupledwave-eqn}.

\subsection{Reversibility}
\label{subsec-reversibility}

The canonical operators $\bM$ and $\bK$ in (\ref{J1J2-def_skew-sym_R11}) are reversible in
the following sense.  There exists an involution
\begin{equation}\label{reversor-def}
  {\bf R}={\rm diag}(1,-1,-1,1)\,,
\end{equation}
with
\begin{equation}\label{MK-reversor}
  {\bf RM} = -{\bf MR}\qand {\bf RK}=-{\bf KR}\,.
\end{equation}
  The implication of this symmetry is that if $S$ is reversible,
  that is $S({\bf R}Z)=S(Z)$, then ${\bf R}Z(-x,-t)$ is a solution whenever
  $Z(x,t)$ is a solution.

\section{Solitary waves and linearization}
\setcounter{equation}{0}
\label{sec-canonicalform}

The class of PDEs that we take as a starting point
for the development of the theory of linear stability of solitary waves is
\begin{equation}\label{MKS-def}
  {\bf M}Z_t + {\bf K}Z_x = \nabla S(Z)\,,\quad Z\in\R^4\,,
\end{equation}
  with $S:\R^4\to\R$ a smooth scalar-valued function, $\nabla$ the
  gradient on $\R^4$, and ${\bf M},{\bf K}$ are $4\times 4$
  skew-symmetric matrices.  
  The extension to higher dimension phase space is straightforward
  in principle, modulo a proliferation of indices, and is discussed
  in \S\ref{sec-cr}.
  
  The skew-symmetric matrices are required to satisfy
  \begin{equation}\label{Jc-def}
    {\rm det}\big(\J(c)\big)\neq 0\,,\quad \mbox{for all}\ c\in\mathscr{C}\,,
\quad \J(c) := {\bf K}+c{\bf M}\,.
        \end{equation}
  The structure (\ref{MKS-def}) and property (\ref{Jc-def})
  are hypothesis ({\bf H1}).
When $\J(c)$ has a non-trivial kernel
the theory goes through with the phase space restricted
to the complement of the kernel of $\J(c)$.
This is rare but does happen; an example is in \cite{bm95}. 

  In (\ref{Jc-def}), $\mathscr{C}$ is just
  an open subset of $\R$.  In examples, the existence of solitary waves
  will inform the definition of $\mathscr{C}$.
  With assumption ({\bf H1}), 
  $(\R^4,\bO)$ is a symplectic vector space with symplectic form
  \begin{equation}\label{Omega-def}
  \bO(\bu,\bv) = \langle\J(c)\bu,\bv\rangle\,,\quad \forall \bu,\bv\in\R^4\,.
  \end{equation}
         The operator $\J(c)$ is not canonical and $\J(c)^2\neq-{\bf I}$,
        so does not define a complex structure either.  However, neither of
      these properties are required for the theory here.
  
  In addition to the multisymplectic Dirac operators
  in  \S\ref{subsec-PDO_1-1}, the class of PDEs (\ref{MKS-def})
  also includes the case
where $\bM$ is of rank two and $\bK$ is of rank four.  The
KdV equation and NLS equation are multisymplectic Hamiltonian PDEs of
this latter type \cite{bd99,bd01}.

\subsection{Solitary wave solutions}
\label{sec-solitarywaves}

The abstract form (\ref{MKS-def})
is equivariant with respect to the translation group with action
\[
  {\bf T}_sZ(x,t) =  Z(x+s,t)\,,\quad \forall s\in\R\,,
  \]
  that is, ${\bf T}_s Z(x,t)$ is a solution of (\ref{MKS-def}) whenever
  $Z(x,t)$ is solution.  This latter property follows since
  $\bM$, $\bK$ and $S(Z)$ do not depend explicitly on $x$.
  A relative equilibrium associated with this group
  is a solution of the form
\begin{equation}\label{basic-sw}
Z(x,t) = {\bf T}_{ct}\Zh(x) := \Zh(\xi)\ \mbox{with}\ \xi=x+ct\,,
\end{equation}
where $c\in\R$.  This relative equilibrium solution
is called a solitary wave when the following asymptotic conditions
are operational
\begin{equation}\label{sw-asymptotics}
\lim_{\xi\to\pm\infty} \|\Zh(\xi)\| = 0\,,
\end{equation}
with the convergence exponential.  
As in \S\ref{Dprimeprime-odes}, the relative equilibrium $\Zh$ can
be characterized as
a constrained critical point problem.  Let
\begin{equation}\label{H-I-def}
H(\Zh) = \int_{-\infty}^{+\infty}\left[
  S(\Zh) - \fr \langle {\bf K} \Zh_\xi,\Zh\rangle \right]\,\rd\xi \qand
I(\Zh) = \int_{-\infty}^{+\infty}\fr \langle {\bf M} \Zh_\xi,\Zh\rangle\,\rd\xi\,.
\end{equation}
The functional $I(\Zh)$ is called the momentum of the solitary wave as
it is the conserved functional associated, via Noether's Theorem,
to the $x-$translation symmetry of (\ref{MKS-def}).
Indeed, the
symplectic Noether argument is the same as (\ref{a.3}):
the symplectic operator $\bM$
acting on the generator of the translation group, $\Zh_\xi$, generates
the gradient of a functional
\begin{equation}\label{a.3-sw}
  {\bf M}\Zh_\xi = \nabla I(\Zh)\,,
\end{equation}
where here $\nabla I$ is the gradient in $\R^4\otimes L^2(\R)$.
Analogous to (\ref{a.3}), the
left-hand side of (\ref{a.3-sw})
is the product of two commuting skew-symmetric
operators ${\bf M}$ and $\partial_\xi$.
The operator $\bM$
appears in both $I(\Zh)$ and the governing equation (\ref{MKS-def}) and
this will be useful for connecting $dI/dc$ to the Evans function.

Solitary wave solutions then correspond to critical points of $H(\Zh)$
restricted to level sets of the function $I(\Zh)$, with
Euler-Lagrange equation $\nabla H(\Zh)-c\nabla I(\Zh) = 0$, which
works out to
\[
0 = \nabla H(\Zh)-c\nabla I(\Zh) = \nabla S(\Zh) - \bK \Zh_\xi -c\bM\Zh_\xi =
\nabla S(\Zh) - \J(c)\Zh_\xi\,;
\]
that is, the solitary wave is a homoclinic orbit
of the Hamiltonian ODE
\begin{equation}\label{Zh-ode}
\J(c)\Zh_\xi = \nabla S(\Zh)\,.
\end{equation}
In the critical point argument, the speed $c$ is a Lagrange multiplier.
Solitary waves come in
one parameter families parameterized by $c$; that is, $\Zh(\xi,c)$,
and the family is non-degenerate when
\begin{equation}\label{dIdc-sw}
\frac{d\ }{dc}I\circ\Zh\neq 0\quad \mbox{for all}\ c\in\mathscr{C} \,. 
\end{equation}
This derivative will be represented by $\frac{dI}{dc}$.
The homoclinic orbit solution of (\ref{Zh-ode})
is assumed to be transversely constructed for
all $c\in\mathscr{C}$.  The concept of ``transversely constructed'' is
discussed below in (\ref{Pi-def-4}).

It is assumed that there exists a solitary wave solution of the
form (\ref{basic-sw})-(\ref{sw-asymptotics}) satisfying (\ref{Zh-ode}),
transversely constructed, with
(\ref{dIdc-sw}), and a smooth function of $\xi$ and $c$,
for all $c\in\mathscr{C}$. This is hypothesis ({\bf H2}).

\subsubsection{Reversible solitary waves}
\label{ssubsec-rev-sw}

The equation for solitary waves (\ref{Zh-ode}) is said to be reversible
if there exists an involution ${\bf R}$ with
\[
  {\bf R}\J(c)=-\J(c){\bf R} \qand S({\bf R}\Zh) = S(\Zh)\,.
  \]
  The first identity is satisfied with (\ref{reversor-def}) when
  $\J(c)=\bK+c\bM$ with $\bK$ and $\bM$ satisfying (\ref{MK-reversor}).
  \vspace{.15cm}

  \noindent{\bf Definition 4.1.} {\it
  A solitary wave solution $\Zh(\xi,c)$ is a ``reversible solitary
  wave'' when the equation (\ref{Zh-ode}) is reversible with involution
  ${\bf R}$ and the solitary wave satisfies}
  \begin{equation}\label{reversible-sw}
    {\bf R}\Zh(-\xi,c) = \Zh(\xi,c)\,.
\end{equation}

\subsection{Linearization and the operator ${\bf L}$}
\label{subsec-linear-and-L}

The second variation of $H-cI$ is a linear operator
\begin{equation}\label{second-variation}
  {\bf L}(\xi,c) := D^2 H(\Zh) - c D^2I(\Zh)\,.
\end{equation}
It is this operator to which the GSS spectral condition is applied.  Here
the operator ${\bf L}$ will play an important role in the Evans function
theory but the spectrum of ${\bf L}$, other than its zero eigenvalue,
will not enter the theory, being replaced
by transversality and the coefficient (\ref{Pi-def-a}).
Using (\ref{H-I-def}) another representation of ${\bf L}$ is
\begin{equation}\label{L-def-0}
  {\bf L}W = D^2S(\widehat Z)W - \bK W_\xi -
  c\bM W_\xi =
  {\bf B}(\xi,c)W -  \J(c) W_\xi\,,
\end{equation}
 with
  \begin{equation}\label{B-def}
    {\bf B}(\xi,c) = D^2S(\Zh(\xi,c))\,.
  \end{equation}
Differentiating (\ref{Zh-ode}) with respect to $\xi$ shows that
the tangent vector to the solitary wave is in the kernel of ${\bf L}$
\begin{equation}\label{ker-L}
  {\bf L}\widehat Z_\xi = 0 \,.
\end{equation}
It is assumed that ${\rm Ker}({\bf L})\cap L^2(\R) = {\rm span}\{
\widehat Z_\xi \}$.

The above properties of the linearization constitute hypothesis ({\bf H3}). 
  
In the analysis of the Evans function an equation for $\widehat Z_c$
will be needed.  Differentiate (\ref{Zh-ode}) with respect to $c$
  \[
  \big({\bf K}+c{\bf M}\big)(\Zh_c)_\xi + \bM\Zh_\xi = {\bf B}(\xi,c)\Zh_c\,,
  \]
   or
\begin{equation}\label{Zh-ode-diff-c}
{\bf L}\Zh_c = \bM\Zh_\xi \,.
\end{equation}
Note the similarity with the second equation in the ODE case (\ref{zeta12-def}).

The pair (\ref{ker-L}) and (\ref{Zh-ode-diff-c}) form a Jordan chain of
length exactly two.  This property is confirmed by noting
that length three would require the
existence of $W$ satisfying ${\bf L} W = \bM \Zh_c$.
The proof is the same as that in (\ref{jc-length-3}) and (\ref{jc-length-3a}).
By Hypotheses ({\bf H2}) and ({\bf H3})
and the symmetry of ${\bf L}$, the solvability
condition for ${\bf L} W = \bM \Zh_c$ is
  \begin{equation}\label{solvability}
    0= \int_{-\infty}^{+\infty}\langle \Zh_\xi,\bM \Zh_c\rangle\,\rd\xi
    = -\int_{-\infty}^{+\infty}\langle \nabla I(\Zh),\Zh_c\rangle\,\rd\xi =
    -\frac{dI}{dc}\,,
  \end{equation}
  which is non-zero by hypothesis ({\bf H2}).

\subsection{Eigenvalues and Eigenvectors in the system at infinity}
\label{subsec-SW-linearization}

We begin the set up of the Evans function by taking a dynamical
systems viewpoint of the linearization (\ref{L-def-0}).  The
linearized ODE is
\begin{equation}\label{MKS-steady-linear}
  \J(c) Z_\xi = {\bf B}(\xi,c)Z\,.
\end{equation}
Due to the boundary condition (\ref{sw-asymptotics}) the operator
${\bf B}$ is asymptotic to a constant matrix
\[
\lim_{\xi\to\pm\infty} {\bf B}(\xi,c) = {\bf B}^\infty(c)\,.
\]
The ``system at infinity'' for the steady problem is
$\J(c) Z_\xi = {\bf B}^\infty(c) Z$ which can be solved explicitly,
\[
Z(\xi) = \sum_{j=1}^4 z_j\zeta_j\re^{\mu_j\xi}\,,
\]
where $z_j$ are arbitrary complex constants, and $\mu_j(0,c)$ are eigenvalues
determined by
\[
\Delta(\mu,0;c) := {\rm det}\big[ {\bf B}^\infty(c) - \mu \J(c)\big] =0\,.
\]
The zero in one of the arguments anticipates the introduction of the stability
exponent $\lambda$ in the next section.  The vectors $\zeta_j$ are the
eigenvectors satisfying
\begin{equation}\label{zeta_j-def}
\big[{\bf B}^\infty(c) - \mu_j(c,0)\J(c)\big]
\zeta_j(c,0) = 0 \,, \quad j=1,\ldots,4\,.
\end{equation}

\begin{figure}[ht]
\begin{center}
\includegraphics[width=8.0cm]{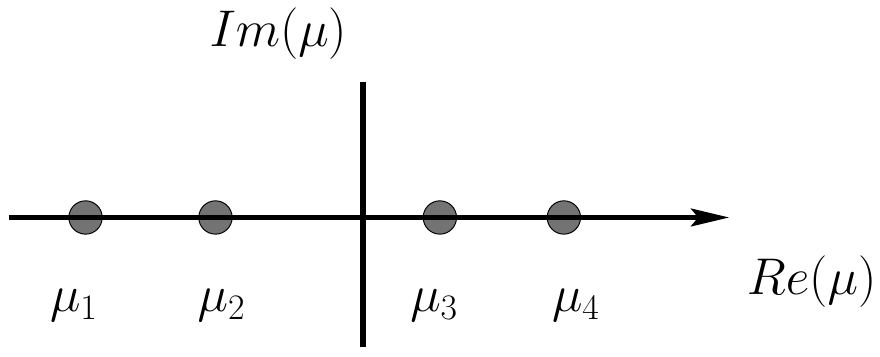}
\end{center}
\caption{Typical position of spatial exponents in the complex $\mu-$plane.}
\label{fig_mu-plane}
\end{figure}

The adjoint eigenvectors are of the form $\J(c)\overline{\eta_j(c,0)}$ with
$\eta_j(c,0)$ satisfying
\begin{equation}\label{adj.7}
\left[ {\bf B}^\infty(c) + \mu_j(c,0)\J(c)\right]\eta_j(c,0)= 0\,,\quad j=1,\ldots,4\,.
\end{equation}
The adjoint eigenvectors are normalized as
\begin{equation}\label{adj.8}
  \langle\J(c)\eta_i(c,0),\zeta_j(c,0)\rangle
  := {\bm\Omega}(\eta_i,\zeta_j) = \delta_{i,j}\,,\quad i,j=1,\ldots,4\,.
\end{equation}
  It is assumed that the spectrum of ${\bf B}^\infty(c)$ has a two-two splitting.
  Introducing a numbering the splitting is represented as
  $\mu_1(0,c)$ and $\mu_2(0,c)$ with negative real part, and
  $\mu_3(0,c)$ and $\mu_4(0,c)$ with positive real part.  For now,
  the eigenvalues are assumed to be real and are shown schematically in
  Figure \ref{fig_mu-plane}.
The case where the eigenvalues form a complex
  quartet is discussed below in \S\ref{sec-complex-eigs}.

  When $\lambda=0$ the sum of the eigenvalues is zero since
  \begin{equation}\label{eigenvalue-sum}
    \mu_1(c,0) + \mu_2(c,0) + \mu_3(c,0)+\mu_4(c,0) =
       {\rm Tr}\big(\J(c)^{-1}{\bf B}^\infty(c)\big)=0\,,
  \end{equation}
  using skew-symmetry of $\J(c)$ and symmetry of ${\bf B}^\infty(c)$.

  \subsection{Orientation}
  \label{subsec-orientation}
  
  Let $\{{\bf e}_1,\ldots,{\bf e}_4\}$ be any fixed (independent of $c$)
  basis, not necessarily the standard basis, for the
    phase space, $\R^4$ and fix an orientation by defining
    \begin{equation}\label{orientation-1}    
    \vol := {\bf e}_1\wedge{\bf e}_2\wedge{\bf e}_3\wedge{\bf e}_4\,.
\end{equation}
    The eigenvectors are related to the orientation by,
\begin{equation}\label{orientation-1a}
  \zeta_1(c,0)\wedge\zeta_2(c,0)\wedge\zeta_3(c,0)\wedge\zeta_4(c,0) =
  K(c,0) \vol\,,\quad\mbox{with}\ K\neq 0\,,\quad \forall\ c\in\mathscr{C}\,.
\end{equation}
  The simple and real $\mu-$eigenvalue assumption, the normalization
  (\ref{adj.8}), and the orientation property (\ref{orientation-1a})
  constitute hypothesis ({\bf H4}).

  \subsection{The $\xi-$dependent stable and unstable subspaces}
  \label{subsec-xi-spaces}
  
  Consistent with the two-two splitting are solutions of the $\xi-$dependent
  equation (\ref{MKS-steady-linear})
  \begin{equation}\label{Es-Eu-def}
  E^s(\xi,0) = {\rm span}\big\{ \Zh_\xi,{\bf a}^+\big\}
  \qand
  E^u(\xi,0) = {\rm span}\big\{ \Zh_\xi,{\bf a}^-\big\}\,,
  \end{equation}
  where $\Zh_\xi$ decays exponentially as $\xi\to\pm\infty$
  and ${\bf a}^{\pm}$ are the other solutions which satisfy
  \begin{equation}\label{a-pm-infty}
    \lim_{\xi\to +\infty}\re^{-\mu_1\xi}{\bf a}^{+}(\xi,c) = \zeta_1(c,0) \qand
     \lim_{\xi\to -\infty}\re^{-\mu_4\xi}{\bf a}^{-}(\xi,c) = \zeta_4(c,0)\,,
  \end{equation}
  with the convergence exponential.
  In general $\ba^{\pm}$ are \emph{not bounded} as $\xi\to\mp\infty$.
  \vspace{.15cm}

  \noindent{\bf Remark.} Here we have associated $\Zh_\xi$ with
  the exponents $\mu_2$ and $\mu_3$, and have associated ${\bf a}^{\pm}$
  with the exponents $\mu_1$ and $\mu_4$.  It could be the
  other way around.  However we will show in \S\ref{sec-Evansfunction} that
  the Evans function is independent of the permutation $3\leftrightarrow 4$.
\vspace{.15cm}

  Existence of ${\bf a}^{\pm}$ follows from the stable/unstable manifold
  theorem.  This strategy for constructing ${\bf a}^{\pm}$ is
  used by \textsc{Galvao \& Gelfreich}~\cite{gg11} and is effective
  in both the case of real eigenvalues and complex quartets.
  The existence of ${\bf a}^{\pm}$ can also be proved
  directly, by analyzing the linear system (\ref{MKS-steady-linear})
  using asymptotic theory for linear ODEs \cite{coppel,kp13}.
\vspace{.15cm}
  
  \noindent{\bf Proposition 4.1.} {\it $E^s(\xi,0)$ and $E^u(\xi,0)$ are
    Lagrangian subspaces with respect to the symplectic structure
    $\bO$ in (\ref{Omega-def}).}
  \vspace{.15cm}
  
  \noindent{\bf Proof.}  The proof is given for $E^s$.
  It is required to show that
  \begin{equation}\label{Omega-zero}
  \bO(\Zh_\xi,\ba^+)=0\quad\mbox{for all}\ \xi\in\R\,.
  \end{equation}
  Differentiate this expression, use (\ref{MKS-steady-linear}), and the
  skew symmetry of $\J(c)$,
  \[
  \begin{array}{rcl}
    \displaystyle\frac{d\ }{d\xi}\bO(\Zh_\xi,\ba^+)
    &=& \langle\J(c)\Zh_{\xi\xi},{\bf a}^{+}\rangle +
    \langle\J(c)\Zh_\xi, (\ba^+)_\xi\rangle \\[3mm]
&=& \langle {\bf B}(\xi,c)\Zh_{\xi},{\bf a}^{+}\rangle -
    \langle\Zh_\xi,\J(c) (\ba^+)_\xi\rangle \\[3mm]
    &=& \langle {\bf B}(\xi,c)\Zh_{\xi},{\bf a}^{+}\rangle -
    \langle\Zh_\xi, {\bf B}(\xi,c) \ba^+\rangle \\[3mm]
    &=& 0 \,,
    \end{array}
  \]
  using symmetry of ${\bf B}(\xi,c)$.
  This proves that $\bO(\Zh_\xi,\ba^+)$ is a constant for
  all $\xi\in\R$.  Now use the fact that $\Zh_\xi$
  and $\ba^+$ both go to zero as $\xi\to\infty$ to conclude (\ref{Omega-zero}).
  A similar proof confirms that $E^u$ is a Lagrangian subspace.
  $\hfill\blacksquare$

  \subsection{Transversality, the Lazutkin invariant, and orientation}
  \label{subsec-lazutkin}

  The Lazutkin invariant of a homoclinic orbit is defined in (\ref{Pi-def-a}).
    It is a property of the intersection between the stable and unstable manifolds which
  form the homoclinic orbit, and it is a symplectic invariant
  \cite{gl01,gg11,cb15}. The fact that $\Pi$ is independent of $\xi$ follows from
  \[
  \begin{array}{rcl}
    \frac{d\ }{d\xi}\bO({\bf a}^-,{\bf a}^+) &=&
    \bO\big( ({\bf a}^-)_\xi,{\bf a}^+\big) + \bO\big({\bf a}^-,({\bf a}^+)_\xi\big)\\[2mm]
&=&
    \langle {\bf B}(\xi,c){\bf a}^-,{\bf a}^+\rangle - \langle{\bf a}^-,{\bf B}(\xi,c){\bf a}^+\rangle\\[2mm]
    &=&
    \langle {\bf B}(\xi,c){\bf a}^-,{\bf a}^+\rangle - \langle{\bf B}(\xi,c){\bf a}^-,{\bf a}^+\rangle\\[2mm]
    &=& 0 \,.
  \end{array}
  \]
  A homoclinic orbit is said to be transversely constructed when
  \begin{equation}\label{Pi-def-4}
\Pi :=  \bO({\bf a}^-,{\bf a}^+)\neq 0\quad\mbox{for all}\ c\in\mathscr{C} \,.
  \end{equation}
 To fix the sign of $\Pi$ we have to synchronize
        the orientations of the stable and unstable subspace, relative
        to the ambient space.
              The unstable subspace satisfies
\begin{equation}\label{c.1}
{\rm span}\{\Zh_\xi,{\bf a}^-\} \to {\rm span}\{\zeta_3,\zeta_4\} \quad \mbox{as}\quad
\xi\to-\infty\,.
\end{equation}
Hence there exists a constant $C^{-}$ with the property that
\begin{equation}\label{Cminus}
\re^{-(\mu_3+\mu_4)\xi}\,\Zh_\xi(\xi)\wedge{\bf a}^-(\xi) \to
C^- \zeta_3\wedge\zeta_4\quad\mbox{as}\quad \xi\to -\infty \,.
\end{equation}
Similarly, the stable subspace satisfies
\begin{equation}\label{c.1-stable}
{\rm span}\{\Zh_\xi,{\bf a}^+\} \to {\rm span}\{\zeta_1,\zeta_2\} \quad \mbox{as}\quad
\xi\to+\infty\,.
\end{equation}
Hence there exists a constant $C^{+}$ with the property that
\begin{equation}\label{Cplus}
\re^{-(\mu_1+\mu_2)\xi}\,\Zh_\xi(\xi)\wedge{\bf a}^+(\xi) \to
C^+ \zeta_1\wedge\zeta_2\quad\mbox{as}\quad \xi\to +\infty \,.
\end{equation}
Take $\xi$ to $-\xi$ in this formula, and use the eigenvalue sum
formula (\ref{eigenvalue-sum}),
\begin{equation}\label{Cplus-1}
\re^{-(\mu_3+\mu_4)\xi}\,\Zh_\xi(-\xi)\wedge{\bf a}^+(-\xi) \to
C^+ \zeta_1\wedge\zeta_2\quad\mbox{as}\quad \xi\to -\infty \,.
\end{equation}
Take the wedge product between (\ref{Cminus}) and
(\ref{Cplus-1}) and take the limit $\xi\to-\infty$,
\[
\re^{-2(\mu_3+\mu_4)\xi}\,\Zh_\xi(-\xi)\wedge{\bf a}^+(-\xi)\wedge\Zh_\xi(\xi)\wedge
{\bf a}^-(\xi) \longrightarrow
C^+C^- \zeta_1(c,0)\wedge\zeta_2(c,0)\wedge\zeta_3(c,0)\wedge\zeta_4(c,0)\,,
\]
or
\begin{equation}\label{Pi-normalization}
\lim_{\xi\to-\infty}
\re^{-2(\mu_3+\mu_4)\xi}\,\Zh_\xi(-\xi)\wedge{\bf a}^+(-\xi)\wedge\Zh_\xi(\xi)\wedge
   {\bf a}^-(\xi) = C^+C^- K(c,0)\vol\,,
\end{equation}
using (\ref{orientation-1a}).
We take the sign of $C^+C^-K(c,0)$ to be fixed throughout the analysis.
This fixes the sign of $\Pi$ as the sign of $\ba^{-}$ can not be changed
without changing the sign of $\ba^{+}$ and vice versa.
These properties of $\Pi$ constitute Hypothesis ({\bf H5}).

      \subsection{Asymptotics of $\Zh_\xi$ and the factor $\chi$}
      \label{subsec-varrho}

      The asymptotic behaviour of $\Zh_\xi$ as $\xi\to-\infty$
      is associated with either the $\mu_3$ or the $\mu_4$ eigenvalue.
      We will see below in \S\ref{sec-Evansfunction}
      that the Evans function is independent of this choice.
      Hence, we will assume without loss of generality
      that $\Zh_\xi$ is associated with $\mu_3$ as $\xi\to-\infty$.

      There exists real numbers $\chi_{\pm}$ with
\begin{equation}\label{Zh-asymptotics}
\lim_{\xi\to-\infty}\re^{-\mu_3\xi}\Zh_\xi = \chi^-\zeta_3 \qand
\lim_{\xi\to+\infty}\re^{+\mu_3\xi}\Zh_\xi = \chi^{+}\eta_3\,.
\end{equation}
These expressions can be inverted to give formula for $\chi^{\pm}$
\begin{equation}\label{chi-plusminus-def}
\chi^- = \lim_{\xi\to-\infty}\re^{-\mu_3\xi}\bO(\eta_3,\Zh_\xi)
\qand
\chi^+ = \lim_{\xi\to+\infty}\re^{+\mu_3\xi}\bO(\Zh_\xi,\zeta_3)\,.
\end{equation}
The factor $\chi$ in the derivative of the Evans function is defined by
\begin{equation}\label{chi-def}
\chi = (\chi^+\chi^-)^{-1}\,.
\end{equation}
  
 \noindent{\bf Proposition 4.2.} {\it Suppose the basic state is reversible
 in the sense of Definition 4.1.  Then}
\begin{equation}\label{chi-rev}
{\rm sign}(\chi) = {\rm sign}\big[ \bO(\zeta_3,{\bf R}\zeta_3) \big]\,. 
\end{equation}
    \vspace{.1cm}

    \noindent{\bf Proof.}  Differentiating (\ref{reversible-sw}) with respect to $\xi$ gives
    \begin{equation}\label{RZh-xi}
    -{\bf R}\Zh_\xi(-\xi,c) = \Zh_\xi(\xi,c)\,.
    \end{equation}
    Now use the formula (\ref{Zh-asymptotics}) for the $\Zh_\xi$
    asymptotics as $\xi\to-\infty$, take $\xi\to-\xi$ in this equation,
    and act on it with ${\bf R}$,
    \[
    \lim_{\xi\to+\infty}\re^{+\mu_3\xi}{\bf R}\Zh_\xi(-\xi,c) = \chi^-{\bf R}\zeta_3\,.
    \]
    Substituting (\ref{RZh-xi}) then gives
    \begin{equation}\label{asymp-1}
    \lim_{\xi\to+\infty}\re^{+\mu_3\xi}\Zh_\xi(\xi,c) = -\chi^-{\bf R}\zeta_3\,.
    \end{equation}
    Now use the second asymptotic formula in (\ref{Zh-asymptotics})
    $\xi\to+\infty$
    \begin{equation}\label{asymp-2}
    \lim_{\xi\to+\infty}\re^{+\mu_3\xi}\Zh_\xi(\xi,c) = \chi^+\eta_3\,.
    \end{equation}
    Combining (\ref{asymp-1}) and (\ref{asymp-2}),
    \begin{equation}\label{zeta-eta-3}
    -\chi^-{\bf R}\zeta_3 = \chi^+\eta_3\,.
    \end{equation}
    On the other hand $\zeta_3$ and $\eta_3$ satisfy
    \[
      [{\bf B} - \mu_3\J(c)]\zeta_3=0\qand
      [{\bf B}+\mu_3\J(c)]\eta_3=0\,.
      \]
      Acting on the first equation with ${\bf R}$ gives
      \[
      \begin{array}{rcl}
        0 &=& {\bf R}[{\bf B} - \mu_3\J(c)]\zeta_3\\[2mm]
        &=& [{\bf R}{\bf B} - \mu_3{\bf R}\J(c)]\zeta_3\\[2mm]
        &=& [{\bf B}{\bf R} + \mu_3\J(c){\bf R}]\zeta_3\\[2mm]
        &=& [{\bf B} + \mu_3\J(c)]{\bf R}\zeta_3\,.
      \end{array}
      \]
The fact that ${\bf BR}={\bf RB}$ follows from differentiating
      $S({\bf R}Z)=S(Z)$ twice and setting $Z=0$.
Therefore there exists $\vartheta\in\R$ such that
\begin{equation}\label{vartheta-def}
{\bf R}\zeta_3 = \vartheta\eta_3\,.
\end{equation}
Act on both sides of (\ref{vartheta-def}) with ${\bf J}(c)$
and take the pairing with $\zeta_3$ to get
\[
\langle{\bf J}(c){\bf R}\zeta_3,\zeta_3\rangle = \vartheta\langle{\bf J}(c)\eta_3,\zeta_3\rangle\,,
\]
and so $\vartheta = -\bO(\zeta_3,{\bf R}\zeta_3)\rangle$.
Applying this formula and (\ref{vartheta-def}) to (\ref{zeta-eta-3}),
\[
    -\chi^-\vartheta\eta_3 = \chi^+\eta_3\quad
\Rightarrow\quad -\vartheta\chi^-=\chi^+\quad\Rightarrow\quad
\bO(\zeta_3,{\bf R}\zeta_3)\big(\chi^-\big)^2 = \chi^{-1}\,.
\]
Comparing signs on each side of this expression proves the proposition.
      $\hfill\blacksquare$

\subsection{Summary of the key geometric properties of a wave}
\label{subsec-summary}

      To summarize, there are three key properties of the family of
      solitary waves that will feed into the stability analysis:
      (a) the derivative of the momentum $dI/dc$ in (\ref{dIdc-sw}), (b) the transversality coefficient
      $\Pi$ in (\ref{Pi-def-a}), and
      (c) the asymptotic property (\ref{Zh-asymptotics}) which
      generates $\chi$ in (\ref{chi-def}).

\section{Linear stability and the Evans function}
\setcounter{equation}{0}
\label{sec-Evansfunction}

The linearization of the PDE
(\ref{MKS-def}) about the solitary wave solution (\ref{basic-sw}) is
\begin{equation}\label{MKS-linearisation}
  {\bf M} Z_t + \J(c) Z_\xi =  {\bf B}(\xi,c) Z\,, \quad Z\in \R^4\,.
\end{equation}
Introduce the spectral ansatz $Z(x,t) = \re^{\lambda t} \bu(\xi,\lambda)$. Then
the eigenvalue problem for $\lambda\in\C$ is
\begin{equation}\label{evans-eqn}
\bu_\xi = {\bf A}(\xi,\lambda) \bu\,,\quad \bu\in\C^4\,,\quad\lambda\in\Lambda\,,
\end{equation}
for some open set $\Lambda\subset\C$, with
\begin{equation}\label{A-def}
  {\bf A}(\xi,\lambda) := \J(c)^{-1}(
  {\bf B}(\xi,c)  - \lambda {\bf M})\,.
\end{equation}

\noindent{\bf Definition 5.1.} There
exists an unstable eigenvalue, in the linearization
about the solitary wave, if there exists a solution of
(\ref{evans-eqn}), for some $\lambda\in\C$ with ${\rm Re}(\lambda)>0$,
and $\bu(\cdot,\lambda)\in \R^4\otimes L^2(\R)$.
\vspace{.15cm}

The asymptotic condition (\ref{sw-asymptotics}) assures that
\begin{equation}\label{A-infty-integral}
\int_{-\infty}^{+\infty}\|{\bf A}(\xi,\lambda) - {\bf A}^\infty(\lambda)\|\,\rd\xi <+\infty\,,\quad \lambda\in\Lambda\,.
\end{equation}
In this integral the ``system at infinity'' is defined by
\begin{equation}\label{A-infty-def}
  {\bf A}^\infty(\lambda) := \lim_{\xi\to\pm\infty}{\bf A}(\xi,\lambda)\quad
  \forall\ \lambda\in\Lambda\,,
\end{equation}
with the dependence on $c$ suppressed for brevity.
This limit is assumed to exist for any fixed $c$ and
uniformly for $\lambda\in\Lambda$.
The spectrum of ${\bf A}^\infty(\lambda)$ consists of four eigenvalues
with the two-two splitting
\begin{equation}\label{mu-assumptions}
\mbox{Re}(\mu_1) \leq \mbox{Re}(\mu_2) < 0 < \mbox{Re}(\mu_3) \leq \mbox{Re}(\mu_4)\quad\forall\ \lambda\in\Lambda\,.
\end{equation}
The set $\Lambda$ is defined below.

The eigenvalue problem $[{\bf A}^{\infty}(\lambda)-\mu{\bf I}]\zeta=0$
is unwrapped into
\begin{equation}\label{zeta-def}
\big[{\bf B}^\infty(c) - \lambda\bM - \mu_j(c,\lambda)\J(c)\big]\zeta_j(c,\lambda) = 0\,,\quad j=1,\ldots,4\,,\ \forall\ \lambda\in\Lambda\,.
\end{equation}
The extension of (\ref{orientation-1a}) to include the
$\lambda-$dependent eigenvectors is
\begin{equation}\label{orientation-2}
\mathcal{V}(c,\lambda) := \zeta_1(c,\lambda)\wedge\zeta_2(c,\lambda)\wedge\zeta_3(c,\lambda)\wedge\zeta_4(c,\lambda) =
  K(c,\lambda) \vol\,,
\end{equation}
with $K(c,\lambda)\neq 0$, for all $c\in\mathscr{C}$ and $\lambda\in\Lambda$.

The adjoint eigenvalue problem
$[{\bf A}^{\infty}(\lambda)^H-\overline{\mu}{\bf I}]\psi=0$,
where the superscript $H$ denotes complex-conjugate transpose,
is unwrapped to
\begin{equation}\label{eta-def}
\big[{\bf B}^\infty(c) + \lambda\bM + \mu_j(c,\lambda)\J(c)\big]\eta_j(c,\lambda) = 0\,,\quad j=1,\ldots,4\,,\ \forall\ \lambda\in\Lambda\,,
\end{equation}
with $\eta_j=\J(c)^{-1}\overline{\psi_j}$.
The adjoint eigenvectors are normalized as
\begin{equation}\label{adj.8a}
  \langle\J(c)\eta_i(c,\lambda),\zeta_j(c,\lambda)\rangle
  := {\bm\Omega}(\eta_i,\zeta_j) = \delta_{i,j}\,,\quad i,j=1,\ldots,4\,.
\end{equation}
With this normalization, we can define a dual to $\mathcal{V}$ in
(\ref{orientation-2}) as
\begin{equation}\label{orientation-2-dual}
  \mathcal{V}^*(c,\lambda) := \J(c)\eta_1(c,\lambda)\wedge\J(c)\eta_2(c,\lambda)\wedge\J(c)\eta_3(c,\lambda)\wedge\J(c)\eta_4(c,\lambda)\,.
\end{equation}

The continuous spectrum is defined by
\begin{equation}\label{sigma-cont}
\sigma^{\footnotesize\mbox{cont}} = \Big\{\, \lambda\in\C\ :\ {\rm det}[ {\bf B}^\infty(c)-\lambda\bM -\ri \kappa\J(c)]=0\,,\quad \kappa\in\R\,\Big\}\,.
\end{equation}
Normally, it is assumed that the continuous spectrum
is purely imaginary, $\sigma^{\footnotesize\mbox{cont}}\subset\ri\R$, thereby limiting instability
considerations to unstable eigenvalues.  Here,
no special assumption is imposed on the continuous
spectrum since our main interest
is in the derivatives of the Evans function at $\lambda=0$, and
the existence of unstable point spectra. But
we will need $\lambda=0$ in the set $\Lambda$.

The set $\Lambda$ is an open set in the complex plane,
including the origin, such that for all $\lambda\in\Lambda$ the
four eigenvalues of ${\bf A}^\infty(\lambda)$ are simple and
satisfy the constraints (\ref{mu-assumptions}).
The above properties of the eigenvalues and eigenvectors at infinity
and the definition of the set $\Lambda$ constitutes hypothesis ({\bf H6}).
The assumption of simple eigenvalues can be relaxed by working
on exterior algebra spaces, or using maximally analytic eigenvectors
\cite{bd03}.

\subsection{Constructing the Evans function}
\label{subsec-evansfcn}

There are many equivalent ways of defining the Evans function
(e.g.\ Chapters 8--10 in \cite{kp13}).
The direct approach is to take the
wedge product of the individual vector-valued solutions of
(\ref{evans-eqn}). We will first define the Evans function
that way, and then introduce an equivalent definition which pairs
solutions of (\ref{evans-eqn}) with solutions of the adjoint equation,
and brings in the $\J(c)-$symplectic structure. 

Using standard asymptotic theory for
ODEs \cite{coppel,kp13} there are four $(\xi,\lambda)-$dependent vectors satisfying
\begin{equation}\label{bu-eqn-def}
({\bf u}_j)_\xi = {\bf A}(\xi,\lambda){\bf u}_j\,,\quad j=1,\ldots,4\,,
\end{equation}
with the asymptotic properties
\begin{equation}\label{bu-asymptotics}
\begin{array}{rclrcl}
&&\displaystyle\lim_{\xi\rightarrow+\infty} e^{-\mu_1(\lambda)\xi}\bu_1(\xi,\lambda) = \zeta_1(\lambda), \quad &&\displaystyle\lim_{\xi\rightarrow+\infty} e^{-\mu_2(\lambda)\xi}\bu_2(\xi,\lambda) = \zeta_2(\lambda), \\[2mm]
&&\displaystyle\lim_{\xi\rightarrow-\infty} e^{-\mu_3(\lambda)\xi}\bu_3(\xi,\lambda) = \zeta_3(\lambda), \quad &&\displaystyle
  \lim_{\xi\rightarrow-\infty} e^{-\mu_4(\lambda)\xi}\bu_4(\xi,\lambda) = \zeta_4(\lambda)\,,
\end{array}
\end{equation}
suppressing the dependence on $c$ as it is now secondary.  The stable
subspace is represented by $\bu_1\wedge\bu_2$, with
\[
\lim_{\xi\to+\infty} \bu_1(\xi,\lambda)\wedge\bu_2(\xi,\lambda) = 0\,,
\]
with the convergence exponential, and the unstable subspace
is represented by $\bu_3\wedge\bu_4$, with
\[
\lim_{\xi\to-\infty} \bu_3(\xi,\lambda)\wedge\bu_4(\xi,\lambda) = 0\,,
\]
with the convergence exponential. If these two spaces intersect then there is a solution which decays exponentially
as $\xi\to\pm\infty$, thereby generating an eigenfunction with eigenvalue
$\lambda$.  The Evans function captures these interesections and hence
eigenvalues.  It is defined by
\begin{equation}\label{evans-def-1}
D(\lambda)\vol = \re^{-\tau(c)\lambda\xi}
\bu_1(\xi,\lambda)\wedge\bu_2(\xi,\lambda)\wedge\bu_3(\xi,\lambda)\wedge\bu_4(\xi,\lambda)\,,
\end{equation}
with $\tau(c)\lambda$ related to the trace of ${\bf A}(\xi,\lambda)$ by
\begin{equation}\label{trace-A}
{\rm Tr}({\bf A}(\xi,\lambda)) =
{\rm Tr}\big(\J(c)^{-1}{\bf B}(\xi,\lambda)
-\lambda\J(c)^{-1}\bM\big)=
-\lambda{\rm Tr}\big(\J(c)^{-1}\bM\big):= \lambda\tau(c)\,,
\end{equation}
since $\J(c)$ is skew-symmetric and ${\bf B}(\xi,c)$ is symmetric.
The function $D(\lambda)$
is independent of $\xi$ and an analytic function of $\lambda$ for
all $\lambda\in\Lambda$ and $\lambda\in D^{-1}(0)$ is an
element of the point spectrum \cite{agj90,kp13}. 

Here an equivalent definition of the Evans function,
in terms of individual vectors of
(\ref{evans-eqn}) and its symplectic adjoint, are used.
\vspace{.15cm}

\noindent{\bf Proposition 5.1.} {\it The symplectic adjoint of $\bu_\xi=
  {\bf A}(\xi,\lambda)\bu$ is}
  \begin{equation}\label{bw-eqn}
    \bw_\xi = {\bf A}(\xi,-\lambda)\bw\,,\quad \bw\in\C^4\,,
    \end{equation}
      {\it and $\bw(\xi,\lambda)$ is an analytic functions of $\lambda$
        for $\lambda\in\Lambda$.}
  \vspace{.15cm}

  \noindent{\bf Proof.}
  Start with a construction of the adjoint of ${\bf A}(\xi,\lambda)$.
  The symplectic adjoint is that linear operator
  ${\bf A}^{\footnotesize\mbox{$\star$symp}}$ satisfying
  \[
  \bO({\bf A}^{\footnotesize\mbox{$\star$symp}}\bw,\bu) = \bO(\bw,{\bf A}\bu)\,.
  \]
  Calculating, using ${\bf A}(\xi,\lambda) = \J(c)^{-1}({\bf B}(\xi,c)-\lambda\bM)$,
  gives
  \[
  \begin{array}{rcl}
    \bO(\bw,{\bf A}\bu) &=& \langle\J\bw,{\bf A}\bu\rangle\\[2mm]
    &=& \langle{\bf A}^T\J\bw,\bu\rangle\\[2mm]
    &=& -\langle ({\bf B}(\xi,c)+\lambda\bM)\J^{-1}\J\bw,\bu\rangle\\[2mm]
    &=& -\langle\J\J^{-1} ({\bf B}(\xi,c)+\lambda\bM)\bw,\bu\rangle\\[2mm]
    &=& -\langle\J{\bf A}(\xi,-\lambda)\bw,\bu\rangle\\[2mm]
    &=& \bO({\bf A}^{\footnotesize\mbox{$\star$symp}}\bw,\bu)\,,
  \end{array}
  \]
  giving ${\bf A}^{\footnotesize\mbox{$\star$symp}}=-{\bf A}(\xi,-\lambda)$.
  Now bring in the derivative
  \[
   \bO(\bw_\xi,\bu) + \bO(\bw,\bu_\xi)  =  \frac{d\ }{d\xi}\bO(\bw,\bu)\,,
  \]
and so
\[
\begin{array}{rcl}
  \bO(\bw,\bu_\xi-{\bf A}(\xi,\lambda)\bu) &=&
  \bO(\bw,\bu_\xi) - \bO(\bw,{\bf A}(\xi,\lambda)\bu)\\[2mm]
  &=& -\bO(\bw_\xi,\bu)- \bO({\bf A}^{\footnotesize\mbox{$\star$symp}}\bw,\bu)
  +\frac{d\ }{d\xi}\bO(\bw,\bu)\\[2mm]
  &=& -\bO(\bw_\xi,\bu)+ \bO({\bf A}(\xi,-\lambda)\bw,\bu)
  +\frac{d\ }{d\xi}\bO(\bw,\bu)\\[2mm]
  &=& -\bO(\bw_\xi-{\bf A}(\xi,-\lambda)\bw,\bu)
  +\frac{d\ }{d\xi}\bO(\bw,\bu)\,,
\end{array}
\]
or
\[
\bO(\bw,\bu_\xi-{\bf A}(\xi,\lambda)\bu)
+\bO(\bw_\xi-{\bf A}(\xi,-\lambda)\bw,\bu) = \frac{d\ }{d\xi}\bO(\bw,\bu)\,.
\]
At this point,
one can either bring in integration, or note that if $\bu$ satisfies
  (\ref{evans-eqn}) and $\bw$ satisfies (\ref{bw-eqn}) then  $\bO(\bw,\bu)$ is
  independent of $\xi$.
  
  Another way to prove the form of the adjoint (\ref{bw-eqn}) is
  to start
  with the formal adjoint of (\ref{evans-eqn}),
\begin{equation}\label{W-adjoint-1}
W_\xi = -{\bf A}(\xi,\lambda)^H W \,,\quad W\in\C^4\,,\quad \lambda\in\Lambda\,,
\end{equation}
where the superscript $H$ signifies complex conjugate transpose, and $W$ is
not an analytic function of $\lambda$.
Pre-multiply by $\J(c)^{-1}$ and use the special form of
${\bf A}(\xi,\lambda)$ in (\ref{A-def}),
\[
(\J(c)^{-1}W)_\xi =-\J(c)^{-1}{\bf A}(\xi,\lambda)^H W 
= \J(c)^{-1}[{\bf B}(\xi,c)+\overline{\lambda}\bM]
  (\J(c)^{-1} W)\,,
\]
or
\[
(\J(c)^{-1}W)_\xi = {\bf A}(\xi,-\overline{\lambda})(\J(c)^{-1}W)\,.
\]
Define
\begin{equation}\label{w-conjugation}
  {\bf w}(\xi,\lambda) = \J(c)^{-1} \overline{W(\xi,\lambda)}\,,
\end{equation}
then the vector-valued functions $\bw(\xi,\lambda)$ are analytic and satisfy
(\ref{bw-eqn}).$\hfill\blacksquare$
\vspace{.15cm}
  
There are four solutions of (\ref{bw-eqn}) with the asymptotic properties
\begin{equation}\label{bw-asymptotics}
\begin{array}{rclrcl}
&&\displaystyle\lim_{\xi\rightarrow-\infty} e^{+\mu_1(\lambda)\xi}\bw_1(\xi,\lambda) = \eta_1(\lambda), \quad &&\displaystyle\lim_{\xi\rightarrow-\infty} e^{+\mu_2(\lambda)\xi}\bw_2(\xi,\lambda) = \eta_2(\lambda), \\[2mm]
&&\displaystyle\lim_{\xi\rightarrow+\infty} e^{+\mu_3(\lambda)\xi}\bw_3(\xi,\lambda) = \eta_3(\lambda), \quad &&\displaystyle
  \lim_{\xi\rightarrow+\infty} e^{+\mu_4(\lambda)\xi}\bw_4(\xi,\lambda) = \eta_4(\lambda)\,,
\end{array}
\end{equation}
where $\eta_j$, $j=1,\ldots,4$ are adjoint eigenvectors (\ref{eta-def}).
\vspace{.15cm}

\noindent{\bf Theorem 5.2.} {\it With the orientation (\ref{orientation-1}) and
  the normalizations (\ref{adj.8a}), the Evans function (\ref{evans-def-1})
  can be transformed to the representation} 
  \begin{equation}\label{Evans-uw-version}
  D(\lambda)\vol = {\rm det}\left[
    \begin{matrix} \bO(\bw_3,\bu_3) & \bO(\bw_3,\bu_4) \\
      \bO(\bw_4,\bu_3) & \bO(\bw_4,\bu_4) \end{matrix}\right]\mathcal{V}(c,\lambda)\,.
  \end{equation}
\vspace{.15cm}

\noindent The proof of the theorem is based on two remarkable formulae,
which bring the symplectic structure into the Evans function
(\ref{evans-def-1}),
    \begin{equation}\label{W-Upsilon-identity}
    \re^{-\tau(c)\lambda\xi}\bu_1(\xi,\lambda)\wedge
    \bu_2(\xi,\lambda)\intp \mathcal{V}^* =
    \J(c)\bw_3(\xi,\lambda)\wedge\J(c)\bw_4(\xi,\lambda)\,,
    \end{equation}
    and
        \begin{equation}\label{eta-identity}
          \J(c)\eta_3\wedge\J(c)\eta_4 = \zeta_1\wedge\zeta_2\intp
          \mathcal{V}^*\,.
        \end{equation}
        where $\intp$ is the interior product, and $\mathcal{V}^*$ is
        the dual four-form defined in (\ref{orientation-2-dual}).
  The proof of these
        two formulae and Theorem 5.2 requires
        an excursion into exterior algebra and they are given in Appendix
        \ref{app-a}.

\subsection{The equivalence class of Evans functions}
\label{subsec-eqiv-evans-fcns}

The two Evans functions (\ref{evans-eqn}) and (\ref{Evans-uw-version})
are equivalent.  Evans functions form an equivalence class.
Two representations $D^A(\lambda)$ and $D^B(\lambda)$, of an
  Evans function, are equivalent
  if there exists a non-vanishing analytic function $C(\lambda)$ such that
  $D^A(\lambda)=C(\lambda)D^B(\lambda)$ for all $\lambda\in\Lambda$.
  When orientation of $D^{A}(\lambda)$ and $D^{B}(\lambda)$
  along the real axis is of interest, then
  $C(\lambda)\big|_{\lambda\in\R}$ is required to be real and positive.
  
We will select a particular representation from this class.
Using $\mathcal{V}(c,\lambda)=
        K(c,\lambda)\vol$ from (\ref{orientation-2}),
 the explicit formula for the scalar-valued
        function $D(\lambda)$ in (\ref{Evans-uw-version}) is
        \[
        D(\lambda) = K(c,\lambda){\rm det}\left[
    \begin{matrix} \bO(\bw_3,\bu_3) & \bO(\bw_3,\bu_4) \\
      \bO(\bw_4,\bu_3) & \bO(\bw_4,\bu_4) \end{matrix}\right]\,.
        \]
        However, as noted in (\ref{orientation-2}),
        $K(c,\lambda)$ is non-zero for all $c\in\mathscr{C}$ and
        all $\lambda\in\Lambda$.  Hence it can be factored out giving
        the equivalent function
        \begin{equation}\label{D-lambda-def}
        D(\lambda) := {\rm det}\left[
    \begin{matrix} \bO(\bw_3,\bu_3) & \bO(\bw_3,\bu_4) \\
      \bO(\bw_4,\bu_3) & \bO(\bw_4,\bu_4) \end{matrix}\right]\,.
        \end{equation}
        Here and henceforth, it is this function that will be used as
``the Evans function''.
        The proof of the following property of the Evans function is
        evident from inspection.
        \vspace{.15cm}

        \noindent{\bf Propostion 5.3.} {\it The Evans function
          (\ref{D-lambda-def}) is invariant under the permutation
          $3\leftrightarrow 4$.}
        \vspace{.15cm}

        \noindent This property will be useful in the proof of the derivative formula
        as it shows that it does not matter whether $\Zh_\xi$ is asymptotic
        to the direction $\zeta_3$ or $\zeta_4$ as $\xi\to-\infty$.

        The factor $\chi$ can be eliminated by defining an equivalent
        Evans function with $\chi$ built in
  \begin{equation}\label{Evans-uw-version-R-pm}
  D_R(\lambda) = {\rm det}\left[
    \begin{matrix} \bO(\chi^+\bw_3,\chi^-\bu_3) & \bO(\chi^+\bw_3,\bu_4) \\
      \bO(\bw_4,\chi^-\bu_3) & \bO(\bw_4,\bu_4) \end{matrix}\right]\,.
  \end{equation}
  This Evans function is equivalent to (\ref{D-lambda-def}).  The advantage
  is that the derivative formula (\ref{Dpp-transversality}) simplifies to
  \[
  D_R''(0) = \Pi\,\frac{dI}{dc}\,.
  \]
  However, we prefer to use the definition (\ref{D-lambda-def})
  of the Evans function as
  a starting point and bring the factor $\chi$ out
  explicitly in the derivative formula.

If a class of PDEs of interest has solitary waves with $\Pi$ and
$\chi$ of
one sign, then an equivalent Evans function
can be defined without $\Pi$ or $\chi$, e.g.
\[
D_{new}(\lambda) = \frac{1}{\chi\Pi}\,D(\lambda)\,.
\]
In this case ${D_{new}}''(0) = dI/dc$. For example, this case is relevant for
the Evans function in PW (e.g.\ equation (\ref{Dpp-1})).

        \subsection{The $\lambda\to0$ limit of the stable and unstable spaces}
        \label{subsec-orientation-lambdazero}

        As above, we take, without loss of generality,
        $\bu_3(\xi,\lambda)\to{\rm span}\{\Zh_\xi\}$.  
        Hence we have the following limits
             \begin{equation}\label{bu-scaled}
             \lim_{\lambda\to0}{\bu}_3(\xi,\lambda) = \frac{1}{\chi_{-}} \Zh_\xi \qand
             \lim_{\lambda\to0}{\bu}_4(\xi,\lambda) = {\bf a}^-(\xi)\,.
             \end{equation}
             The first limit is due to the fact that
             \[
             \re^{-\mu_3(\lambda)\xi}{\bf u}_3 \to \zeta_3\quad\mbox{but}\quad
             \re^{-\mu_3(\lambda)\xi}\Zh_\xi \to \chi_{-}\zeta_3\,,
             \]
             using (\ref{bu-asymptotics}) and (\ref{Zh-asymptotics}).
             The $\lambda-$limit of the second
             term in (\ref{bu-scaled}) follows from (\ref{bu-asymptotics}) and (\ref{a-pm-infty}).
             
             Similarly, the adjoint eigenfunctions have $\lambda\to0$ limits
             \begin{equation}\label{bw-scaled}
             \lim_{\lambda\to0}{\bw}_3(\xi,\lambda) = \frac{1}{\chi^+}\Zh_\xi \qand
             \lim_{\lambda\to0}{\bw}_4(\xi,\lambda) = {\bf a}^+(\xi)\,.
             \end{equation}
  The form of the eigenfunctions $\bu_j$ and adjoint eigenfunctions
  $\bw_j$, the limits (\ref{bu-scaled}) and (\ref{bu-scaled}),
  and the form of the Evans function (\ref{D-lambda-def}) constitute 
  hypothesis ({\bf H7}).

\section{Derivatives of the Evans function}
\setcounter{equation}{0}
\label{sec-deriv-evans}

In this section the main result of the paper is proved, the connection
between $D''(0)$, transversality, $\chi$, and $dI/dc$ that was asserted in the
introduction in formula (\ref{Dpp-transversality}). 
\vspace{.15cm}

\noindent{\bf Theorem 6.1.} {\it The Evans function (\ref{D-lambda-def}),
  associated with
  the linearization of
  the class of PDEs (\ref{MKS-def}), about the solitary wave solutions (\ref{basic-sw}), under the hypotheses ({\bf H1})-({\bf H7}), has the following derivatives at $\lambda=0$,
  \begin{equation}\label{d-dp-dpp-theorem}
    D(0)=0\,,\quad D'(0)=0\,,\quad\mbox{and}\quad D''(0)=\chi\Pi\,\frac{dI}{dc}\,,
  \end{equation}
  where $\Pi$ is the transversality coefficient (\ref{Pi-def-a}),
  $\chi$ is the $\Zh_\xi$ coeffient (\ref{chi-def}), and
  $I(c)$ is the momentum (\ref{H-I-def})
  evaluated on the $c-$dependent family of solitary waves.}
\vspace{.15cm}

\noindent{\bf Proof.} 
The fact that $D(0)=0$ follows from the fact that
${\bf L}$ has a zero eigenvalue (\ref{ker-L}).  However, the proof
in the context of the Evans function is a bit more interesting,
as it brings in the Lagrangian subspace property of $E^s$ and $E^u$.
The proof proceeds with the evaluation of $D(\lambda)$ in
(\ref{D-lambda-def}) at $\lambda=0$,
\begin{equation}\label{D-zero}
D(0) = {\rm det}\left[ \begin{matrix} 0 & 0 \\ 0 & -\Pi
    \end{matrix}\right]\,.
\end{equation}
The zeros in the first column and row are confirmed by noting
that $\bu_3$ and $\bw_3$ satisfy (\ref{bu-scaled}) and (\ref{bw-scaled}).
Now use skew-symmetry of $\bO$ and the Lagrangian subspace property of
the stable and unstable subspaces (Proposition 4.1) to conclude
\[
\begin{array}{rcl}
  \bO(\bw_3,\bu_3)\big|_{\lambda=0} &=& \chi\bO (\widehat Z_\xi,\widehat Z_\xi)=0\\[2mm]
  \bO(\bw_3,\bu_4)\big|_{\lambda=0} &=& (\chi^+)^{-1}\bO(\widehat Z_\xi,\ba^-) = 0\\[2mm]
  \bO(\bw_4,\bu_3)\big|_{\lambda=0} &=&  (\chi^-)^{-1}\bO(\ba^+,\widehat Z_\xi) = 0\\[2mm]
  \bO(\bw_4,\bu_4)\big|_{\lambda=0} &=& \bO(\ba^+,\ba^-)=-\Pi\,.
\end{array}
\]
Substitution into (\ref{D-lambda-def}) then confirms the
zero structure in (\ref{D-zero}).

To prove the properties of the first and second derivatives
of $D(\lambda)$, define the entries of the matrix in $D(\lambda)$ as
  \begin{equation}
D(\lambda) = \left({d_1(\lambda)d_2(\lambda) - d_3(\lambda)d_4(\lambda)}_{\vphantom{0_0}}\right)
\end{equation}
  where
\[
\begin{array}{c c}
d_1(\lambda) = \bO(\w_3,\u_3)  &\qquad d_2(\lambda) = \bO(\w_4,\u_4)\,, \\[2mm]
d_3(\lambda) = \bO(\w_3,\u_4), &\qquad d_4(\lambda) = \bO(\w_4,\u_3)\,.
\end{array}
\]
It follows from (\ref{D-zero}) that
\begin{equation}\label{dn-zero}
d_1(0)=d_3(0)=d_4(0)=0\,,\qand d_2(0) = -\Pi\,.
\end{equation}
Computing the first derivative
\[
  D'(\lambda)
  = \big( d_1'(\lambda)d_2(\lambda) +
  d_1(\lambda)d_2'(\lambda) -
  d_3'(\lambda)d_4(\lambda)- d_3(\lambda)d_4'(\lambda)\big) \,.
\]
Evaluating at $\lambda=0$ and using (\ref{dn-zero})
\begin{equation}\label{Dp-zero}
D'(0) = - d_1'(0)\Pi\,.
\end{equation}
Now
\begin{equation}\label{d1p-lambda}
d_1'(\lambda) =
\bO(\partial_\lambda\bw_3,\bu_3) + \bO(\bw_3,\partial_\lambda\bu_3)\,.
\end{equation}
For $\partial_\lambda\bu_3$ and $\partial_\lambda\bw_3$,
start with their defining equation,
\[
\J(\u_3)_\xi = [\B - \lambda\M]\bu_3\qand
\J(\bw_3)_\xi = [\B + \lambda\M]\bw_3\,,
\]
and differentiate with respect to $\lambda$,
\begin{equation}\label{u3w3-lambda}
\J(\u_3)_{\xi\lambda}  =  [\B - \lambda\M](\bu_3)_\lambda - \bM\bu_3\qand
\J(\bw_3)_{\xi\lambda}  =  [\B + \lambda\M](\bw_3)_\lambda + \bM\bw_3\,.
\end{equation}
Set $\lambda=0$, 
\[
\begin{array}{rcl}
  {\bf L}(\bu_3)_\lambda &=&\displaystyle
 \bM \bu_3\Big|_{\lambda=0}=  (\chi^-)^{-1}\bM \Zh_\xi\\[4mm]
{\bf L}(\bw_3)_\lambda &=&\displaystyle
 -\bM \bw_3\Big|_{\lambda=0}= - (\chi^+)^{-1}\bM \Zh_\xi\,,
\end{array}
  \]
using (\ref{bu-scaled}) and (\ref{bw-scaled}).  Now use
equation (\ref{Zh-ode-diff-c}), giving
\begin{equation}\label{u3w3-lambda-soln}
  (\bu_3)_\lambda\big|_{\lambda=0} =  (\chi^-)^{-1}\Zh_c + C_a \Zh_\xi\qand
(\bw_3)_\lambda\big|_{\lambda=0} = - (\chi^+)^{-1}\Zh_c + C_b \Zh_\xi\,,
    \end{equation}
    with $C_a$ and $C_b$ arbitrary constants.

Substitute the expressions (\ref{u3w3-lambda-soln}) into (\ref{d1p-lambda})
evaluated at $\lambda=0$,
\[
\begin{array}{rcl}
  d_1'(0) &=& \displaystyle\left[\bO(\partial_\lambda\bw_3,\bu_3) + \bO(\bw_3,\partial_\lambda\bu_3)\right]\bigg|_{\lambda=0}\\[4mm]
  &=& \displaystyle\left[\bO(- (\chi^+)^{-1}\Zh_c+C_b\Zh_\xi, (\chi^-)^{-1}\Zh_\xi) + \bO((\chi^+)^{-1}\Zh_\xi, (\chi^-)^{-1}\Zh_c+C_a\Zh_\xi)\right]\\[4mm]
&=& 2\displaystyle\chi\bO(\Zh_\xi,\Zh_c)\,.
\end{array}
\]
This latter term is zero.  To see this, first show that it is independent
of $\xi$,
\begin{eqnarray*}
\dfrac{d}{d\xi}\bO(\Zh_c,\Zh_\xi) & = & \langle \J(\Zh_c)_\xi,\Zh_\xi \rangle - \langle \Zh_c,\J(\Zh_\xi)_\xi \rangle \\
& = & \langle \B\Zh_c,\Zh_\xi \rangle + \langle \M\Zh_\xi,\Zh_\xi \rangle - \langle \Zh_c,\B\Zh_\xi \rangle \\
& = & 0\,,
\end{eqnarray*}
and so $\bO(\Zh_\xi,\Zh_c)$ is a constant, but this constant
clearly
vanishes at $\xi=\pm\infty$ and so the form is zero for all $\xi$. This proves
that $d_1'(0)=0$ and so $D'(0)=0$ in (\ref{Dp-zero}).

The second derivative is
\begin{equation}\label{Dpp-formula}
\begin{array}{rcl}
  D''(\lambda) &=&
  \bigg( d_1''(\lambda)d_2(\lambda) + d_1'(\lambda)d_2'(\lambda) +
  d_1(\lambda)d_2''(\lambda) +  d_1'(\lambda)d_2'(\lambda)  \\[2mm]
  &&\quad -  d_3''(\lambda)d_4(\lambda) -d_3'(\lambda)d_4'(\lambda)
  - d_3(\lambda)d_4''(\lambda)-d_3'(\lambda)d_4'(\lambda) \bigg)\,.
\end{array}  
\end{equation}
Evaluating (\ref{Dpp-formula})
at $\lambda=0$ eliminates the second derivatives of
$d_j$ for $j=2,3,4$, leaving
\begin{equation}\label{dpp-0-calc}
D''(0) = \big( d_1''(0)d_2(0) +2d_1'(0)d_2'(0) - 2d_3'(0)d_4'(0) \big)\,.
\end{equation}
The second term is zero due to $d_1'(0)=0$ as was shown above.  A similar
argument can be used to show that $d_3'(0)=0$ as follows,
\begin{equation}\label{d3prime}
d_3'(0) =  (\chi^+)^{-1}\bO(\Zh_c,\a^-) +  (\chi^+)^{-1}\bO(\Zh_\xi,\left.(\u_4)_\lambda\right|_{\lambda=0} )\,.
\end{equation}
The sum of the two terms is constant (since $d_3(\lambda)$ and $d_3'(\lambda)$
are independent of $\xi$).  The first term goes to
zero as $\xi\to-\infty$ as both $\Zh_c$ and $\ba^-$ go to zero.
For the second term $\Zh_\xi$ also goes to zero as $\xi\to-\infty$, so
all that is needed is that $(\bu_4)_\lambda\big|_{\lambda=0}$ be bounded
as $\xi\to-\infty$.  But $\bu_4\big|_{\lambda=0}$ goes to zero exponentially
as $\xi\to-\infty$ and $\partial_\lambda\bu_4$ will only add a polynomial
in $\xi$ to the exponential decay, giving zero for the second term
in (\ref{d3prime}) as well.

Hence the second derivative (\ref{dpp-0-calc}) reduces to
\begin{equation}\label{dpp-0-calc-1}
D''(0) = d_2(0) d_1''(0) =  -\Pi d_1''(0)\,.
\end{equation}
To compute $d_1''(0)$ start with $d_1'(\lambda)$ in (\ref{d1p-lambda}).
Using (\ref{u3w3-lambda}), we can write,
\begin{equation}\label{d1pp-calc}
\begin{array}{rcl}
  \partial_\xi \langle \J(\w_3)_\lambda,\u_3 \rangle & = & \langle [\B + \lambda\M](\w_3)_\lambda,\u_3 \rangle + \langle \M\w_3,\u_3 \rangle\\[2mm]
  &&\hspace{2.0cm}
  - \langle (\w_3)_\lambda,[\B - \lambda\M]\u_3 \rangle\,, \\[2mm]
  \partial_\xi \langle \J\w_3,(\u_3)_\lambda \rangle & = &
  \langle [\B + \lambda\M]\w_3,(\u_3)_\lambda \rangle 
  + \langle \w_3,\M\u_3 \rangle\\[2mm]
  &&\hspace{2.0cm}
  - \langle \w_3,[\B - \lambda\M](\u_3)_\lambda \rangle\,,
\end{array}
\end{equation}
which combine to leave us with:
\begin{equation}
\partial_\xi \langle \J(\w_3)_\lambda,\u_3 \rangle = \langle \M\w_3,\u_3 \rangle = -\partial_\xi \langle \J\w_3,(\u_3)_\lambda \rangle.
\end{equation}
If we now take some $R>0$ then we can integrate the first part of this over the range $\xi\in[0,R]$ and the second part over $\xi\in[-R,0]$ to get:
\begin{eqnarray}
  \displaystyle\Big[\langle \J(\w_3)_\lambda,\u_3 \rangle\Big]_{\xi=0}^{\xi=R}
  & = & \displaystyle\int_0^R \!\langle \M\w_3,\u_3 \rangle\, \xd\xi, \\[3mm]
  \displaystyle\Big[-\langle \J\w_3,(\u_3)_\lambda \rangle\Big]_{\xi=-R}^{\xi=0} &=&
  \displaystyle \int_{-R}^0 \!\langle \M\w_3,\u_3 \rangle\, \xd\xi\,.
\end{eqnarray}
They can be combined to give:
\begin{equation}
\left. d_1'(\lambda)^{\vphantom{0^0}}_{\vphantom{0^0}}\right|_{\xi=0} = -\int_{-R}^R \!\langle \M\w_3,\u_3 \rangle\, \xd\xi + \left.\langle \J(\w_3)_\lambda,\u_3 \rangle^{\vphantom{0^0}}_{\vphantom{0^0}}\right|_{\xi=R} + \left.\langle \J\w_3,(\u_3)_\lambda \rangle^{\vphantom{0^0}}_{\vphantom{0^0}}\right|_{\xi=-R}.
\end{equation}
(Note that although this value of $d_1'(\lambda)$ is specifically evaluated at $\xi=0$, since $d_1$ is independent of $\xi$ it will take this value for all $\xi$.) Taking the limit $R\rightarrow\infty$ allows us to write this as
\begin{equation}
d_1'(\lambda) = -\int_{-\infty}^{+\infty} \!\langle \M\w_3,\u_3 \rangle\, \xd\xi + \ell(\lambda)\,,
\end{equation}
where the function $\ell(\lambda)$ is defined as
\[
\ell(\lambda) = \lim_{\xi\rightarrow+\infty} \langle \J(\w_3)_\lambda,\u_3 \rangle + \lim_{\xi\rightarrow-\infty} \langle \J\w_3,(\u_3)_\lambda \rangle.
\]
Differentiate this function with respect to $\lambda$ to get an expression for $d_1''(\lambda)$:
\begin{equation}
d_1''(\lambda) = -\int_{-\infty}^{+\infty} \!\langle \M(\w_3)_\lambda,\u_3 \rangle + \langle \M\w_3,(\u_3)_\lambda \rangle\, \xd\xi + \ell'(\lambda).
\end{equation}
Now
\begin{equation}\label{lprime}
\begin{array}{rcl}
  \ell'(\lambda) &=&\displaystyle
  \lim_{\xi\rightarrow+\infty} \left[\langle \J(\w_3)_{\lambda\lambda},\u_3 \rangle + \langle \J(\w_3)_\lambda,(\u_3)_\lambda \rangle\right] \\[2mm]
  &&\hspace{2.0cm} \displaystyle
  + \lim_{\xi\rightarrow-\infty} \left[\langle \J(\w_3)_\lambda,(\u_3)_\lambda \rangle + \langle \J\w_3,(\u_3)_{\lambda\lambda} \rangle\right]\,.
\end{array}
\end{equation}
However, since
\[
\lim_{\xi\rightarrow+\infty} e^{\mu_3(\lambda)\xi}\w_3 = \eta_3(\lambda)
\]
we can deduce that
\[
\w_3 e^{+\mu_3(\lambda)\xi} = \eta_3 + o(1) \quad \mbox{for}\ \xi\to+\infty\,.
\]
This in turn implies that
\[
(\w_3)_{\lambda\lambda} e^{+\mu_3(\lambda)\xi} = p(\xi,\lambda)\eta_3 + o(1) \ \mbox{for}\ \xi\to+\infty\,,
\]
where $p(\xi,\lambda)$ is a quadratic polynomial in $\xi$. Since the exponential term will dominate the quadratic polynomial this tells us that
\begin{equation}\label{w3-lambdalambda}
\lim_{\xi\rightarrow+\infty} (\w_3)_{\lambda\lambda} = 0.
\end{equation}
The same argument can be used to show that
\begin{equation}\label{u3-lambdalambda}
\lim_{\xi\rightarrow-\infty} (\u_3)_{\lambda\lambda} = 0\,.
\end{equation}
Now set $\lambda=0$ in (\ref{lprime}) and use these asymptotic
properties.  The first term becomes
\[
\lim_{\xi\rightarrow+\infty} \left[\langle
  \left.\J(\w_3)_{\lambda\lambda}\right|_{\lambda=0}, (\chi^-)^{-1}\Zh_\xi \rangle + \langle \J(- (\chi^+)^{-1}\Zh_c + C_b\Zh_\xi), (\chi^-)^{-1}\Zh_c + C_a\Zh_\xi \rangle \right]=0\,,
\]
using (\ref{w3-lambdalambda}) and the vanishing of $\bO(\Zh_\xi,\Zh_c)$.
Set $\lambda=0$ in the second term in (\ref{lprime}),
\[
 \lim_{\xi\rightarrow-\infty} \left[\langle \J(- (\chi^+)^{-1}\Zh_c + C_b\Zh_\xi), (\chi^-)^{-1}\Zh_c + C_a\Zh_\xi \rangle +  (\chi^+)^{-1}\langle \J \Zh_\xi,\left.(\u_3)_{\lambda\lambda}\right|_{\lambda=0} \rangle\right] = 0\,,
 \]
 using (\ref{u3-lambdalambda}) and $\bO(\Zh_\xi,\Zh_c)=0$.
Therefore $\ell'(0)=0$ and
\begin{eqnarray}
d_1''(0) & = & -\int_{-\infty}^{+\infty} \!\langle \M(- (\chi^+)^{-1}\Zh_c + C_b\Zh_\xi), (\chi^-)^{-1}\Zh_\xi \rangle +  (\chi^+)^{-1}\langle \M \Zh_\xi, (\chi^-)^{-1}\Zh_c + C_a\Zh_\xi \rangle\, \xd\xi \nonumber \\
& = & -2\chi\int_{-\infty}^{+\infty} \!\langle \M\Zh_\xi,\Zh_c \rangle\, \xd\xi \nonumber \\
& = & -2\chi\frac{\xd I}{\xd c}\,,
\end{eqnarray}
using the expression for $dI/dc$ in (\ref{solvability}).
Substituting this result into (\ref{dpp-0-calc-1}) then gives
\[
D''(0) = 2\,\chi\,\Pi\,\frac{dI}{dc}\,.
\]
Scaling $D$ to eliminate the $2$ then
completes the proof of the Theorem.$\hfill\blacksquare$
\vspace{.15cm}

\noindent{\bf Remark.} Since (\ref{D-lambda-def}) is ``the Evans function'', scaling
$D$ to eliminate the $2$ should be interpreted as a trivial re-definition 
$D\mapsto \fr D$.
\vspace{.15cm}

Using Theorem 6.1 above and Theorem 10.1 in \cite{cb15} an alternative
formula for the second derivative in terms of the Maslov index is obtained,
\vspace{.15cm}

\noindent{\bf Corollary 6.2.} {\it Under the above hypotheses, an alternative
formula for $D''(0)$ is},
\begin{equation}\label{D-Maslov}
D''(0) =  \chi (-1)^{\mbox{\footnotesize Maslov}}\, \frac{dI}{dc}\,.
\end{equation}
The Maslov index of an orbit, say a homoclinic or periodic orbit,
is a generalization of the Morse index, and equals the Morse index
in special cases (e.g.\ \textsc{Cox, et al.}~\cite{cjls16}, \textsc{Latushkin
\& Sukhtaiev}~\cite{ls18}).  The
use of the Maslov index in the stability of waves is now a burgeoning
subject (e.g. \textsc{Beck, et al.}~\cite{bcjlms18},
\textsc{Chardard et al.}~\cite{cdb09,cb15}, and references therein).
 To go into the details of the Maslov index would take us too far afield.
 The new observation here is just the connection
between the parity of the Maslov index and the Evans function.
\vspace{.15cm}

\noindent{\bf Corollary 6.3.} {\it Let $d_\infty={\rm sign}\big(D(\lambda_\infty)\big)$, for some real $\lambda_\infty\gg 1$. Then there exists an unstable
  real eigenvalue of the spectral problem (\ref{evans-eqn}) when
  $d_{\infty}D''(0)<0$, or}
\begin{equation}\label{d-infty-condition}
\chi \Pi \frac{dI}{dc} d_{\infty}<0\,.
\end{equation}
We will not study $d_\infty$ here as there are well-established strategies
for the determining the
large $\lambda$ behaviour of the Evans function in the Hamiltonian context
(e.g.\ \S1(g) of \cite{pw92}, \S8 of \cite{bd01},
and Appendix B of \cite{bd03}), and in the non-Hamiltonian context
(e.g. \S 5B of \cite{agj90} and Chapter 9 of \cite{kp13}).

\section{Example: a coupled wave equation}
\setcounter{equation}{0}
\label{sec_coupledwave-eqn}

To illustrate the theory it is applied to a nonlinear wave equation in
the class introduced in \S\ref{subsec-coupled-nlws},
\begin{subequations}
\begin{eqnarray}
\phi_{tt} - \phi_{xx} + \partial_\phi V(\phi,v) & = & 0, \\
v_{tt} - v_{xx} - \partial_v V(\phi,v) & = & 0\,,
\end{eqnarray}
\end{subequations}
with
\[
V(\phi,v)=2\phi^2 - 2\phi^3 - 2v^2 + v^3 + \frac{p}{2}(2\phi - v)^2\,.
\]
Hypothesis ({\bf H1}) is satisfied by rewriting the equations
in the form
\begin{equation}
\M Z_t + \K Z_x = \nabla S(Z), \quad Z\in\mathbb{R}^4
\end{equation}
with
\begin{equation}
Z =
\left(\!\begin{array}{c}
\phi \\
u_1 \\
u_2 \\
v
\end{array}\!\right), \quad
\M =
\left(\!\begin{array}{c c c c}
0 & -1 & 0 & 0 \\
1 & 0 & 0 & 0 \\
0 & 0 & 0 & 1 \\
0 & 0 & -1 & 0
\end{array}\!\right), \quad
\K =
\left(\!\begin{array}{c c c c}
0 & 0 & 1 & 0 \\
0 & 0 & 0 & -1 \\
-1 & 0 & 0 & 0 \\
0 & 1 & 0 & 0
\end{array}\!\right)
\end{equation}
where
\[
u_1 = \phi_t - v_x, \qquad u_2 = \phi_x - v_t\,,
\]
and
\begin{equation}\label{S-def-example}
S(Z) = \frac{1}{2}({u_1}^2 - {u_2}^2) + 2\phi^2 - 2\phi^3 - 2v^2 + v^3 + \frac{p}{2}(2\phi - v)^2.
\end{equation}
The symplectic operator $\J(c)$ is
\begin{equation}\label{JcM-example}
\J = \K + c\M =
\left(\!\begin{array}{c c c c}
0 & -c & 1 & 0 \\
c & 0 & 0 & -1 \\
-1 & 0 & 0 & c \\
0 & 1 & -c & 0
\end{array}\!\right)\,,
\end{equation}
with $\det(\J(c))=(1-c^2)^2$.  The other part of hypothesis ({\bf H1})
is satisfied by
taking $c\in\R\setminus{\pm1}$ but the set $\mathscr{C}$ will be
restricted further below for the existence of the
family of solitary waves.

An exact solitary wave solution is given by
\begin{equation}\label{sw-exact-soln}
\phi = \phat(\xi,c) \quad\mbox{and}\quad v = 2\phat(\xi,c), \quad \xi = x + ct
\end{equation}
where
\begin{equation}\label{sw-exact-soln-1}
\phat(\xi,c) = \sech^2(\alpha\xi), \quad \alpha = \dfrac{1}{\sqrt{1-c^2}}.
\end{equation}
Hence the set $\mathscr{C}$ is taken to be
\begin{equation}\label{C-set-def}
\mathscr{C}= \{c\,:\, -1 < c < +1\}\,.
\end{equation}
In terms of the $Z-$coordinates the solitary wave solution is
\begin{equation}\label{Zhat-example}
\Zh(\xi,c) =
\left(\!\begin{array}{c}
\phat \\
-(2-c)\phat_\xi \\
(1-2c)\phat_\xi \\
2\phat
\end{array}\!\right)\,,
\end{equation}
and satisfies
\begin{equation}\label{JS-eqn-example}
\J(c)\Zh_\xi = \nabla S(Z)\,.
\end{equation}
The hypothesis ({\bf H2}) is satisfied: the family of solitary waves
satisfies (\ref{JS-eqn-example}) and is
a smooth function of $\xi$ and $c$ for $c\in\mathscr{C}$.  Transversality,
and the non-degeneracy of
$dI/dc$ will be verified below.

Moreover this solitary wave is reversible.
With the reversor ${\bf R}$, defined in (\ref{reversor-def}),
the operator $\J(c)$ anti-commutes with ${\bf R}$ and $S(Z)$ in
(\ref{S-def-example}) is ${\bf R}-$invariant.  Hence the governing
equation (\ref{JS-eqn-example}) is reversible.  Acting on the solitary wave with
${\bf R}$ gives
\begin{equation}\label{sw-example-reversible}
{\bf R}\Zh(\xi,c) =
\left(\!\begin{array}{c}
\phat(\xi,c) \\
+(2-c)\phat(\xi,c)_\xi \\
-(1-2c)\phat(\xi,c)_\xi \\
2\phat(\xi,c) 
\end{array}\!\right)=
\left(\!\begin{array}{c}
\phat(-\xi,c) \\
-(2-c)\phat(-\xi,c)_\xi \\
(1-2c)\phat(-\xi,c)_\xi \\
2\phat(-\xi,c)
\end{array}\!\right) = \Zh(-\xi,c)\,,
\end{equation}
confirming Definition 4.1.

\subsection{Linearization about the solitary wave}

\begin{figure}[ht]
\begin{center}
\includegraphics[width=4.0cm]{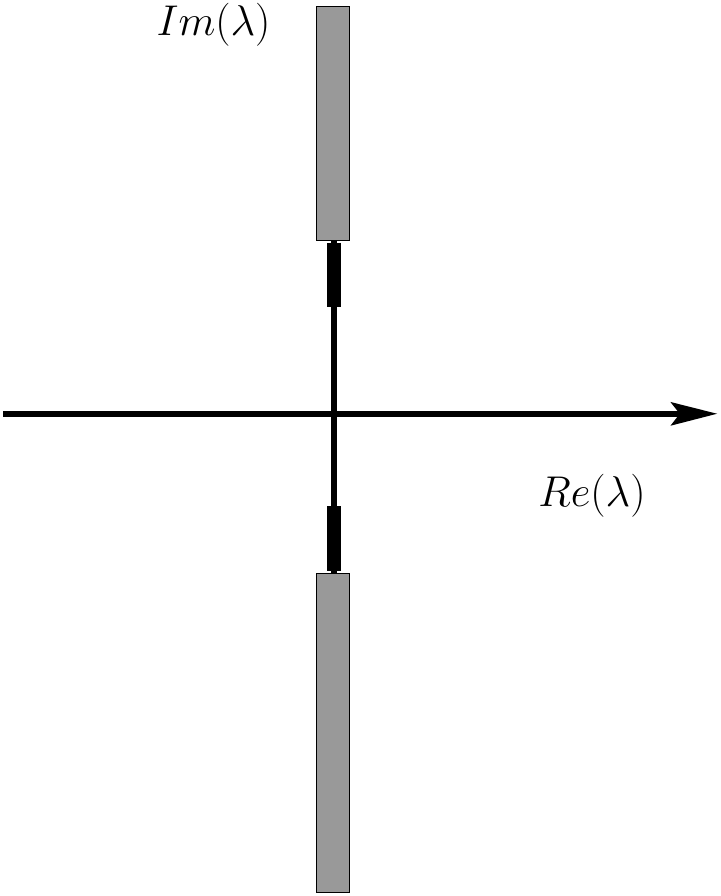}
\end{center}
\caption{Schematic of the four branches of continuous spectra on the imaginary
  axis.}
\label{fig_cont-spec}
\end{figure}

The linearization about the solitary wave solution is
\begin{equation}
\J Z_{\xi} = \left[\B(\xi,c) - \lambda\M\right]Z, \quad Z\in\mathbb{C}^4 \label{lineq}
\end{equation}
with
\begin{equation} \label{Bdef}
\left[\B(\xi,c) - \lambda\M\right] =
\left(\!\begin{array}{c c c c}
4p + 4 - 12\phat(\xi,c) & \lambda & 0 & -2p \\
-\lambda & 1 & 0 & 0 \\
0 & 0 & -1 & -\lambda \\
-2p & 0 & \lambda & p - 4 + 12\phat(\xi,c)
\end{array}\!\right).
\end{equation}
The system at infinity is
\begin{equation}
\J Z_{\xi} = \left[\B^\infty(c) - \lambda\M\right]Z, \quad Z\in\mathbb{C}^4\,,
\end{equation}
with $\left[\B^\infty(c) - \lambda\M\right]$ the same as \eqref{Bdef} but with $\phat$ set to zero. The eigenvalues of the system at infinity are defined by $\Delta(\mu,\lambda)=0$ with
\begin{eqnarray*}
\Delta(\mu,\lambda) & = & {\rm det}\left[\B^\infty(c) - \lambda\M - \mu\J(c)\right] \\
& = & (\mu^2 - \varrho^2)^2 - (8 + 3p)(\mu^2 - \varrho^2) + 16 + 12p
\end{eqnarray*}
where $\varrho=\lambda+c\mu$. The continuous spectrum is
\begin{eqnarray}
\sigma^{\footnotesize\mbox{cont}} & = & \left\lbrace \lambda\in\mathbb{C} : \Delta(\ri\kappa,\lambda) = 0, \quad \kappa\in\mathbb{R} \right\rbrace \nonumber\\
& = & \left\lbrace \lambda = -\ri c\kappa \pm \ri\sqrt{4 + \frac{3}{2}p + \kappa^2 \pm \frac{3}{2}p}, \quad \kappa\in\mathbb{R} \right\rbrace \nonumber \,.
\end{eqnarray}
There are four branches with a gap about the origin,
\[
\sigma^{\footnotesize\mbox{cont}} =
\left\lbrace \lambda = \ri\tilde{\kappa}\ :\ \left\lbrace\tilde{\kappa}\in\mathbb{R}\setminus\left(-2\alpha^{-1},2\alpha^{-1}\right) \right\rbrace\cup
\left\lbrace\tilde\kappa\in\R\setminus\big(-\alpha^{-1}\sqrt{4+3p},\alpha^{-1}\sqrt{4+3p}\big)\right\rbrace\right\rbrace\,.
\]
With $p>0$ the gap consists of $-2\alpha^{-1}<\tilde\kappa<2\alpha^{-1}$.
A schematic of the branches of continuous spectra, showing
the gap about the origin, is illustrated in Figure
\ref{fig_cont-spec}.
The set $\Lambda$ can be taken to be the right-half plane along with
a neighbourhood of the origin.

\begin{figure}[ht]
\begin{center}
\includegraphics[width=8.0cm]{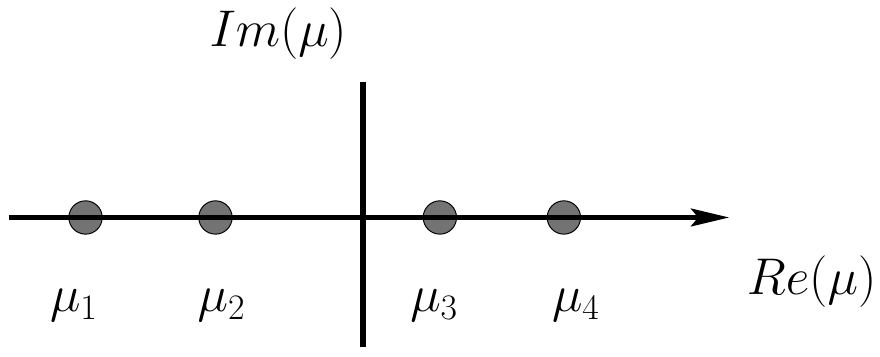}
\end{center}
\caption{Typical position of spatial exponents in the complex $\mu-$plane
  when $p>0$.}
\label{fig_mu-plane-example}
\end{figure}

Now set $\lambda=0$ and set $\Delta(\mu,0)=0$, giving the
four $\mu-$eigenvalues
\[
\mu_1 = -\sqrt{\dfrac{4+3p}{1-c^2}}, \quad \mu_2 = -\sqrt{\dfrac{4}{1-c^2}}, \quad \mu_3 = \sqrt{\dfrac{4}{1-c^2}}, \quad \mu_4 = \sqrt{\dfrac{4+3p}{1-c^2}}\,.
\]
With $p>0$ and $c\in\mathscr{C}$ they are simple and real, and a schematic of the eigenvalue
positions is shown in Figure \ref{fig_mu-plane-example}.
The tangent vector $\Zh_\xi$ is bi-asymptotic to the $\mu-$eigenvalues of
slowest decay ($\mu_2$ and $\mu_3$).  On the other hand,
when $p<0$ (but $p>-\frac{4}{3}$)
the tangent vector $\Zh_\xi$ is bi-asymptotic to the $\mu-$eigenvalues of maximal decay ($\mu_1$ and $\mu_4$).

The eigenvectors when $\lambda=0$ are
\[
\zeta_1 =
z_1\begin{pmatrix}
2 \\  (2c  - 1)\mu_1 \\ (2 -c )\mu_1 \\ 
1\end{pmatrix}\,,
\ \zeta_2 = z_2\begin{pmatrix}
1 \\ (c - 2)\mu_2 \\ (1 -2c)\mu_2 \\ 
2\end{pmatrix}\,,
\ \zeta_3 = z_3\begin{pmatrix}
1 \\ (c-2)\mu_3 \\ (1-2c)\mu_3  \\ 
2\end{pmatrix}\,,
\ \zeta_4 =
z_4\begin{pmatrix}
2 \\ (2c - 1)\mu_4 \\ (2 -c)\mu_4  \\ 
1\end{pmatrix}\,,
\]
with the $z_j$ arbitrary real numbers,
and the adjoint eigenvectors are
\[
\eta_1 =
h_1\begin{pmatrix}
2 \\ (1-2 c)\mu_1 \\ (c-2) \mu_1 \\ 
1\end{pmatrix}\,,\ 
\eta_2 = h_2\begin{pmatrix}
1 \\ (2-c)\mu_2 \\ (2c-1)\mu_2 \\ 
2\end{pmatrix}\,,\ 
\eta_3 = h_3\begin{pmatrix}
1 \\ (2-c)\mu_3 \\ (2c-1) \mu_3 \\ 
2\end{pmatrix}\,,\ 
\eta_4 =
h_4\begin{pmatrix}
2 \\ (1-2c)\mu_4 \\ (c-2) \mu_4 \\ 
1\end{pmatrix}\,,
\]
with $h_j$ arbitrary real numbers.  However,
the normalization $\bO(\eta_i,\zeta_j)=\delta_{i,j}$ gives
the constraints
\begin{equation}\label{hz-formulae}
h_1z_1 = \frac{\alpha}{6\sqrt{4+3p}}\,,\quad h_2z_2 = - \frac{\alpha}{12}\,,
\quad h_3z_3 = \frac{\alpha}{12}\,,\quad
h_4z_4 = -\frac{\alpha}{6\sqrt{4+3p}}\,.
\end{equation}
The full bevy of eigenvector expressions is given for information,
as this is more detail than is needed.
The function $K(c,0)$ is calculated by taking the wedge product of
the eigenvectors
\begin{equation}\label{zeta-wedge4}
\zeta_1\wedge
\zeta_2\wedge \zeta_3\wedge \zeta_4 = K(c,0)\,{\bf e}_1\wedge{\bf e}_2
\wedge{\bf e}_3\wedge{\bf e}_4\,,
\end{equation}
with ${\bf e}_1,\ldots,{\bf e}_4$ the standard basis.  Calculating gives
\[
\begin{array}{rcl}
  K(c,0) &=& \big[-9 (1-c^2)\big( (\mu_1-\mu_4)(\mu_2-\mu_3) \big)z_1z_2z_3z_4\big]\big|_{\lambda=0}\\[2mm]
  &=&     \big[ -36 z_1z_2z_3z_4\sqrt{\alpha^2\lambda^2+4}
  \sqrt{\alpha^2\lambda^2 + 4 + 3p}\,\big]\big|_{\lambda=0}\\[2mm]
   &=&      -72 z_1z_2z_3z_4\sqrt{4 + 3p}\,.
\end{array}
\]
This completes the verification of Hypothesis ({\bf H4}).

Since the solitary wave is reversible (\ref{sw-example-reversible}),
we have
that ${\rm sign}(\chi)$ is determined from
\begin{equation}\label{rev-correction}
\bO(\zeta_3,{\bf R}\zeta_3) = - 6\mu_3(z_3)^2(1-c^2)<0\,.
\end{equation}
By Proposition 4.2, $\chi<0$.  Although not necessary, an explicit formula for $\chi$ can be obtained as well.  Differentiating (\ref{Zhat-example})
with respect to $\xi$ and calculating gives
\[
\lim_{\xi\to-\infty}\re^{-\mu_3\xi}\Zh_\xi = \frac{8\alpha}{z_3}\zeta_3 \qand
\lim_{\xi\to+\infty}\re^{+\mu_3\xi}\Zh_\xi = -\frac{8\alpha}{h_3}\eta_3 \,.
\]
Hence
\begin{equation}\label{chi-formula-example}
\chi^-\chi^+ = \left(\frac{8\alpha}{z_3}\right)\left(-\frac{8\alpha}{h_3}\right) = - 64 \frac{\alpha^2}{z_3h_3} = -768\alpha\quad\Rightarrow\quad
\chi = - \frac{1}{768\alpha}\,,
\end{equation}
using (\ref{hz-formulae}).

\subsection{Derivative of the momentum $dI/dc$}
\label{subsec-dIdc}
       
It is now immediate from Theorem 6.1 that $D(0)=D'(0)=0$.
To compute $D''(0)$, the two factors $\Pi$ and $\frac{\xd I}{\xd c}$ need
to be computed. The easier of the two is $I(c)$ and its derivative:
\begin{eqnarray}
I & = & \frac{1}{2}\int^{+\infty}_{-\infty}\! \langle\M\Zhat_\xi,\Zhat\rangle \,\xd\xi \nonumber\\
& = & \frac{1}{2}\int^{+\infty}_{-\infty}\! (2-c)\phat\phat_{\xi\xi} - (2-c){\phat_\xi}^2 + 2(1-2c){\phat_\xi}^2 - 2(1-2c)\phat\phat_{\xi\xi} \,\xd\xi \nonumber\\
& = & \frac{3c}{2}\int^{+\infty}_{-\infty}\! \phat\phat_{\xi\xi} - {\phat_\xi}^2 \,\xd\xi \nonumber\\
& = & -3c\int^{+\infty}_{-\infty}\! {\phat_\xi}^2 \,\xd\xi \nonumber\\
& = &  -\frac{16}{5}\dfrac{c}{\sqrt{1-c^2}}\,,
\end{eqnarray}
after integrating by parts and using the exact solution (\ref{sw-exact-soln-1}).
Differentiating with respect to $c$,
\begin{equation}
\dfrac{\xd I}{\xd c} = -\frac{16}{5}(1 - c^2)^{-\frac{3}{2}} = -\frac{16}{5}\alpha^3
\end{equation}
which is strictly negative for $c\in\mathscr{C}$.  This completes the
verification of most of hypothesis ({\bf H2}).  Transversality will be
verified below.

\subsection{Computing the transversality coefficient $\Pi$}
\label{subsec-computing-Pi}

To compute the transversality coefficient $\Pi$, the tangent vectors
${\bf a}^\pm$ are needed.  Because of the nature of this example,
these tangent vectors can be computed explicitly.
Express the solutions of \eqref{lineq} by $Z=(\ptil,\util_1,\util_2,\vtil)$. Then $\tilde u_1$ and $\tilde u_2$ can be obtained from \eqref{lineq} as
\begin{equation}\label{u1u2-tilde-def}
\util_1 = \lambda\ptil + c\ptil_\xi - \vtil_\xi \quad\mbox{and}\quad \util_2 = -\lambda\vtil - c\vtil_\xi + \ptil_\xi.
\end{equation}
This allows us to rewrite the first and fourth components of \eqref{lineq} as the coupled pair of equations:
\begin{subequations}\label{tileq}
\begin{eqnarray}
(1 - c^2)\ptil_{\xi\xi} - 2c\lambda\ptil_\xi + (\widehat\varphi - 4p - \lambda^2)\ptil + 2p\vtil & = & 0, \\
(1 - c^2)\vtil_{\xi\xi} - 2c\lambda\vtil_\xi + (\widehat\varphi + p - \lambda^2)\vtil - 2p\ptil & = & 0
\end{eqnarray}
\end{subequations}
where $\widehat\varphi:=12\phat - 4$. If we now introduce the transformation
\[
\ptil = e^{\alpha^2c\lambda\xi}(\psi_1 + 2\psi_2), \quad \vtil = e^{\alpha^2c\lambda\xi}(2\psi_1 + \psi_2)\,,
\]
then equations \eqref{tileq} will become
\begin{subequations}
\begin{eqnarray}
\alpha^{-2}(\psi_1 + 2\psi_2)_{\xi\xi} + (\widehat\varphi - 4p - \alpha^2\lambda^2)(\psi_1 + 2\psi_2) + 2p(2\psi_1 + \psi_2) & = & 0, \\
\alpha^{-2}(2\psi_1 + \psi_2)_{\xi\xi} + (\widehat\varphi + p - \alpha^2\lambda^2)(2\psi_1 + \psi_2) - 2p(\psi_1 + 2\psi_2) & = & 0.
\end{eqnarray}
\end{subequations}
These equations decouple to give
\begin{subequations}
\begin{eqnarray}
\alpha^{-2}(\psi_1)_{\xi\xi} + (\widehat\varphi - \alpha^2\lambda^2)\psi_1 & = & 0, \\
\alpha^{-2}(\psi_2)_{\xi\xi} + (\widehat\varphi - 3p - \alpha^2\lambda^2)\psi_2 & = & 0
\end{eqnarray}
\end{subequations}
which can be rearranged as
\begin{equation}\label{psi1eq}
\begin{array}{rcl}
\alpha^{-2}(\psi_1)_{\xi\xi} + 12\sech^2(\alpha\xi)\psi_1 & = & (4 + \alpha^2\lambda^2)\psi_1, \\[2mm]
\alpha^{-2}(\psi_2)_{\xi\xi} + 12\sech^2(\alpha\xi)\psi_2 & = & (4 + 3p + \alpha^2\lambda^2)\psi_2. \end{array}
\end{equation}
The analysis of these equations for $\lambda\neq0$ and the construction
of the Evans function is given in Appendix \ref{app-b}.  Here we are interested
in the case $\lambda=0$ and the construction of $\ba^{\pm}$.

Now set $\lambda=0$ in the equations (\ref{psi1eq}) for $\psi_1$ and $\psi_2$.
These exact solutions can be used to confirm the kernel condition
(\ref{ker-L}), verifying Hypothesis ({\bf H3}).  $\a^\pm$ will be produced by the following two $\psi_2$ solutions (now denoted by $\psi^\pm$) of \eqref{psi1eq}
\[
\psi^\pm = e^{\mp\alpha\gamma\xi}\left[\pm\dfrac{p\gamma}{5} + \left(1+\frac{6p}{5}\right)\tanh(\alpha\xi) \pm \gamma\tanh^2(\alpha\xi) + \tanh^3(\alpha\xi)\right]
\]
where $\gamma=\sqrt{4+3p}$. The exponents are related to the $\mu-$eigenvalues
by
\[
\mu_1(c,0) = -\alpha\gamma\qand \mu_4(c,0) = + \alpha\gamma\,.
\]
By reversing the transformations to express $\ptil,\util_1,\util_2,\vtil$ in terms of $\psi^\pm$ we find that
\[
\a^+ =
\left(\!\begin{array}{c}
2\psi^+ \\
(2c-1)\psi^+_\xi \\
(2-c)\psi^+_\xi \\
\psi^+
\end{array}\!\right)
\quad\mbox{and}\quad
\a^- =
\left(\!\begin{array}{c}
2\psi^- \\
(2c-1)\psi^-_\xi \\
(2-c)\psi^-_\xi \\
\psi^-
\end{array}\!\right).
\]
%
Substitution into the formula for $\Pi$ then gives
\begin{eqnarray}
\bO(\a^-,\a^+) & = & 4(1-c^2)\psi^-_\xi\psi^+ + (2c-1)^2\psi^-\psi^+_\xi - (2-c)^2\psi^-\psi^+_\xi + (c^2-1)\psi^-_\xi\psi^+ \nonumber\\
& = & 3(1-c^2)\psi^-_\xi\psi^+ + (4c^2-4c+1-4+4c-c^2)\psi^-\psi^+_\xi \nonumber\\
& = & 3\alpha^{-2}\left(\psi^-_\xi\psi^+ - \psi^-\psi^+_\xi\right)\,.
\end{eqnarray}
It is easy to check that this expression is independent of $\xi$ so we can evaluate it at any value of $\xi$ we choose. If we take $\xi=0$ then since
\[
\psi^\pm_\xi = \mp\alpha\gamma\psi^\pm + \alpha e^{\mp\alpha\gamma\xi}\sech^2(\alpha\xi)\left[1 + \frac{6p}{5} \pm 2\gamma\tanh(\alpha\xi) + 3\tanh^2(\alpha\xi)\right]
\]
we get
\begin{eqnarray}
\Pi:=\bO(\a^-,\a^+) & = & 3\alpha^{-1}\left[\left(1 + \frac{6p}{5} - \frac{p\gamma^2}{5}\right)\left(+\frac{p\gamma}{5}\right) - \left(-\frac{p\gamma}{5}\right)\left(1 + \frac{6p}{5} - \frac{p\gamma^2}{5}\right)\right] \nonumber\\
& = & \frac{6p\gamma}{5\alpha}\left(1 + \frac{6p-p(4-3p)}{5}\right) \nonumber\\
& = & \frac{6p\gamma}{25\alpha}\left(5 + 2p - 3p^2\right) \nonumber\\
& = & \frac{6p}{25\alpha}\sqrt{4+3p}\,(5 - 3p)(1 + p)\,.
\end{eqnarray}
%
%
It remains to check the orientation condition (\ref{Pi-normalization}).
The asymptotics of each of the solutions in the
unstable direction, that is, as $\xi\to-\infty$ is
\[
\lim_{\xi\to-\infty}\re^{-\mu_3\xi}\Zh_\xi(\xi,c) = \frac{8\alpha}{z_3}\zeta_3(c,0)\,,
\quad
\lim_{\xi\to-\infty}\re^{-\mu_4\xi}{\bf a}^-(\xi) = -\frac{\mathcal{P}}{z_4}\zeta_4(c,0)\,,
\]
where $\mathcal{P}=\frac{1}{5}[(p+5)\sqrt{4+3p}+6p+10]$,
with $3p+4>0$, and as $\xi\to+\infty$,
\[
\lim_{\xi\to+\infty}\re^{-\mu_2\xi}\Zh_\xi(\xi,c) = -\frac{8\alpha}{z_2}\zeta_2(c,0)\,,
\quad
\lim_{\xi\to+\infty}\re^{-\mu_1\xi}{\bf a}^+(\xi) = -\frac{\mathcal{P}}{z_1}\zeta_1(c,0)\,.
\]
Combining gives
\[
\lim_{\xi\to-\infty}\re^{-(\mu_3+\mu_4)\xi}\Zh_\xi(\xi,c)\wedge{\bf a}^-(\xi) = -\frac{8\alpha\mathcal{P}}{z_3z_4}\zeta_3(c,0)\wedge\zeta_4(c,0)\,,
\]
and
\[
\lim_{\xi\to+\infty}\re^{-(\mu_1+\mu_2)\xi}\Zh_\xi(\xi,c)\wedge{\bf a}^+(\xi) = \frac{8\alpha\mathcal{P}}{z_1z_2}\zeta_1(c,0)\wedge\zeta_2(c,0)\,,
\]
and so, using $\mu_1+\mu_2+\mu_2 + \mu_4=0$,
\begin{equation}\label{Pi-normalization-example}
\lim_{\xi\to-\infty}
\re^{-2(\mu_3+\mu_4)\xi}\,\Zh_\xi(-\xi)\wedge{\bf a}^+(-\xi)\wedge\Zh_\xi(\xi)\wedge
   {\bf a}^-(\xi) = C^-C^+ K(c,0)\vol\,,
\end{equation}
with
\[
C^-C^+K(c,0) = -\frac{64\alpha^2\mathcal{P}^2}{z_1z_2z_3z_4}\,
\left(-288 z_1z_2z_3z_4\sqrt{4 + 3p}\right) = 18432 \alpha^2\mathcal{P}^2\sqrt{4+3p}>0\,.
\]
This quantity is considered fixed and it fixes the sign of $\Pi$,
confirming Hypothesis ({\bf H5}).  Hypotheses ({\bf H6}) and
({\bf H7}) are related to the $\lambda-$dependent linearization and
they are confirmed in Appendix \ref{app-b}.

\subsection{Summary and $D''(0)$}
\label{subsec-summary-and-Dpp}

The three key properties of the solitary wave that feed into $D''(0)$ are
\[
\begin{array}{rcl}
\chi &=&\displaystyle -\frac{1}{768 \alpha} <0 \\[4mm]
\displaystyle \frac{dI}{dc} &=& \displaystyle
-\frac{16}{5}(1-c^2)^{-3/2} <0\\[4mm]
\Pi &=&\displaystyle\frac{6p}{25\alpha}\sqrt{4+3p}(5-3p)(1+p)\,,
\end{array}
\]
with the sign of $\Pi$ dependent on ${\rm sign}(5-3p)$ when $p>0$.

\begin{figure}[ht]
\begin{center}
  \includegraphics[width=5.0cm]{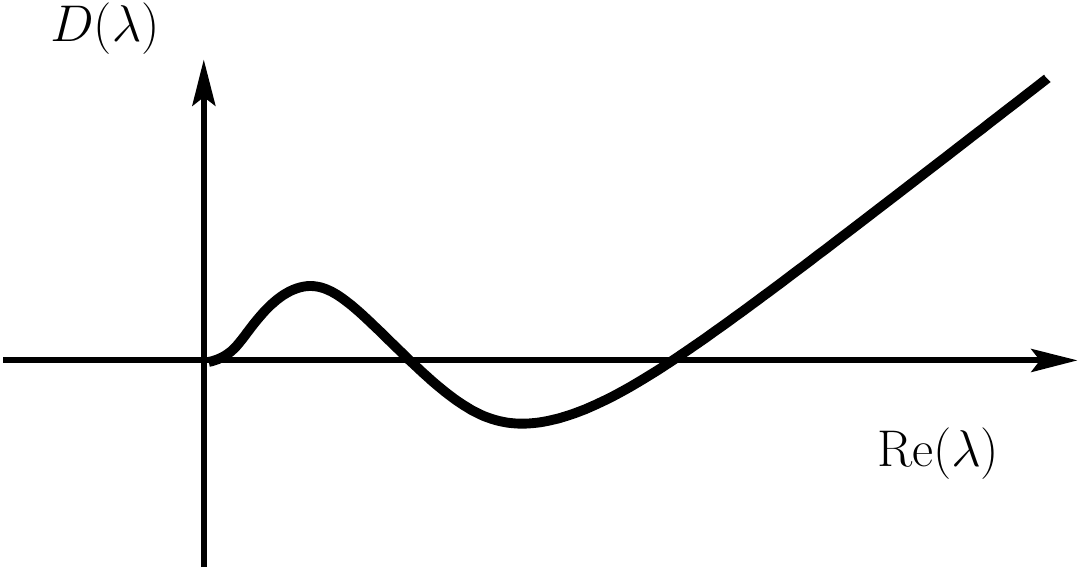}
  \includegraphics[width=5.0cm]{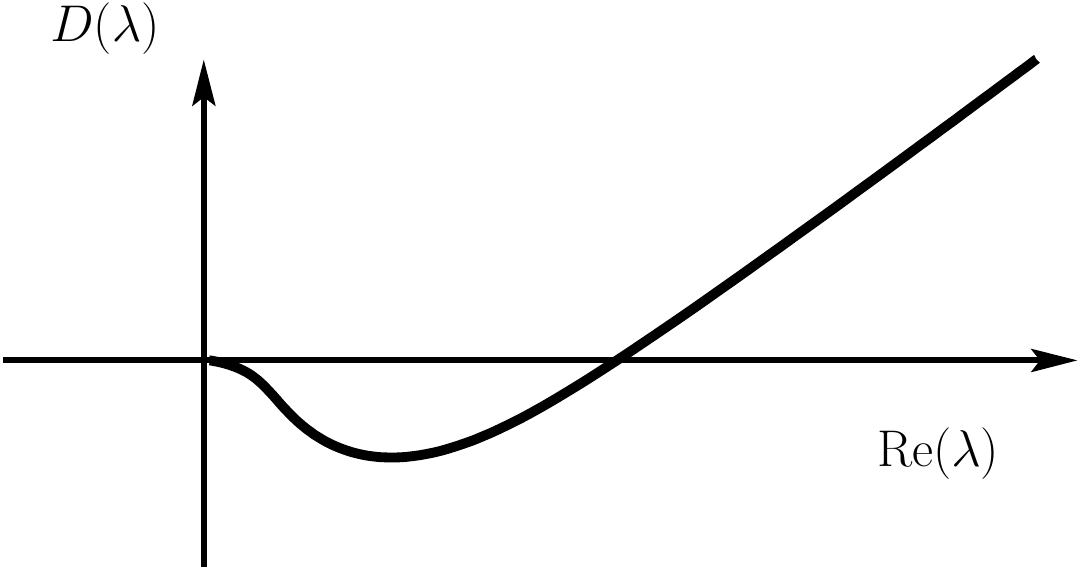}
\end{center}
\caption{Graph of $D(\lambda)$ along the real axis when
  $p<\frac{5}{3}$  (left graph) and $p>\frac{5}{3}$ (right graph).}
\label{fig_Dgraph}
\end{figure}

Combining these formulae and applying Theorem 6.1 gives
$D(0)=0$, $D'(0)=0$, and
\begin{equation}\label{Dpp-example}
D''(0) = \frac{\alpha p}{500}\sqrt{4+3p}\,(5 - 3p)(1 + p)\,.
\end{equation}
With the assumption of $p$ positive, the second derivative is positive for $0<p<5/3$ and negative for $p>5/3$.
Hence, there is at least one interval of $p$ with unstable eigenvalues.
To determine whether $p<\frac{5}{3}$ or $p>\frac{5}{3}$ is the unstable region
we need $d_\infty$.  The sign of $d_\infty$ can be obtained abstractly
(see comments below Corollary 6.3), but here
we have enough information to obtain it explicitly
(see \S\ref{subsec-exact-Evans-fcn-example} below), and we find that
\[
d_{\infty} = {\rm sign}\big( D(\lambda_\infty)\big) = +1\,,\quad \lambda_\infty \gg 1\,.
\]
Hence, applying Corollary 6.3 gives that the solitary wave is unstable when
$5-3p<0$ or $p>\frac{5}{3}$.  The theory is inconclusive in the case
$p<\frac{5}{3}$.  However, we will compute the exact Evans function below and
find that for $p<\frac{5}{3}$ there are two unstable eigenvalues.  Indeed,
we will find that the graph of the Evans function, along the real line,
has the form shown in Figure \ref{fig_Dgraph}.
This result shows both the strength and weakness of the theory.
The strength is that with limited information about the solitary wave,
and the correct sign, the existence of an unstable eigenvalue is proved.
The weakness is that it can not capture an even number of unstable
eigenvalues, and misses unstable eigenvalues that are off the real line.

\subsection{An exact Evans function}
\label{subsec-exact-Evans-fcn-example}

The Evans function for this example can be constructed explicitly, and the
details are given in Appendix \ref{app-b}.  The result is
  \begin{equation}\label{exact-evans-example}
           D_E(\lambda) = \frac{3 S\alpha\lambda^2 }{16(225)^2}(3+x^2)(5-x^2)(3+3p+x^2)(3p+x^2)(5-3p-x^2)\,,
  \end{equation}
           with $x^2=\alpha^2\lambda^2$ and
           \[
S = \sqrt{4+x^2}\sqrt{4+3p+x^2}\,,
\]
and the subscript $E$ is for exact. Taking the second derivative
\[
D_E''(0) = \frac{\alpha p}{500}\sqrt{4+3p}(5-3p)(p+1) \,,
\]
in agreement with $D''(0)$ in (\ref{Dpp-example}).

\begin{figure}[ht]
\begin{center}
\includegraphics[width=4.0cm]{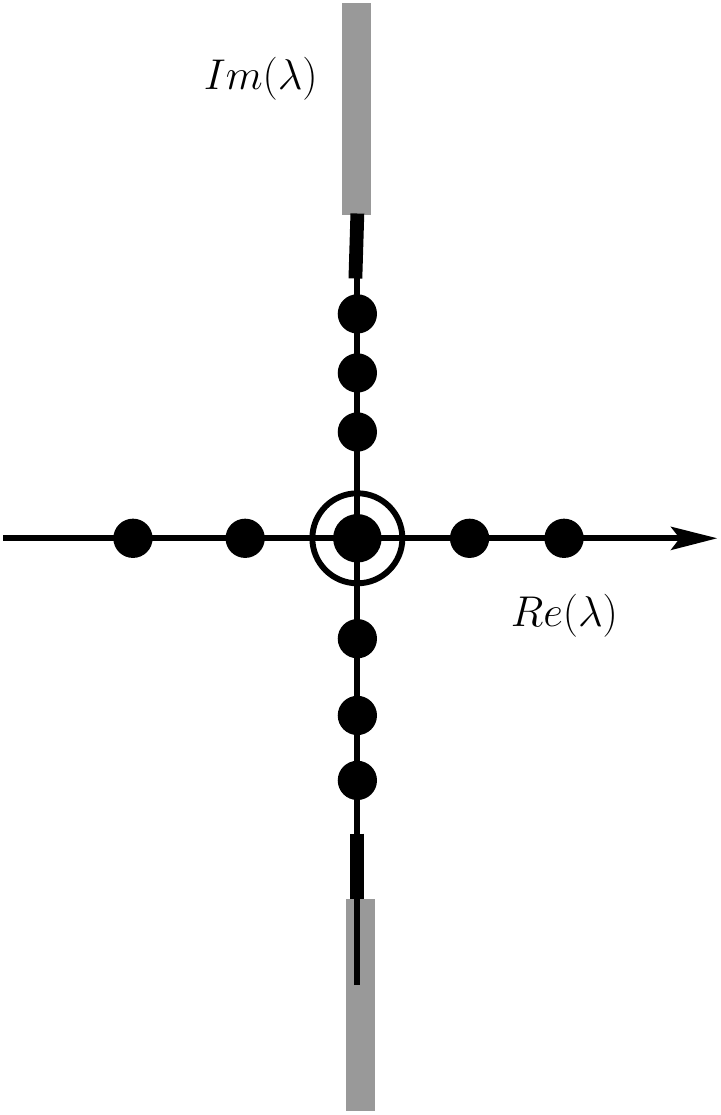}
\end{center}
\caption{Schematic of the branches of continuous spectra on the imaginary
  axis, along with the double zero eigenvalue at the origin, and
  the nonzero eigenvalues, for $0<p<1/3$.}
\label{fig_cont-spec-1}
\end{figure}

We see from the exact expression that there are many more eigenvalues that
go unnoticed, particularly on the imaginary axis in the gap of the continuous
spectrum. The set of eigenvalues for $5-3p>0$ are
           \[
           x=\pm\ri\sqrt{3}\,,\quad x=0^2\,,\quad x=\pm\sqrt{5}\,,\quad
           x=\pm\ri\sqrt{3+3p}\,,\quad x=\pm\ri\sqrt{3p}\,,\quad
           x=\pm\sqrt{5-3p}\,.
           \]
           When $0<p<1$
           the three positive purely imaginary eigenvalues satisfy
           \[
           0 < \ri\frac{\sqrt{3p}}{\alpha} < \ri\frac{\sqrt{3}}{\alpha} < \ri\frac{\sqrt{3+3p}}{\alpha}\,.
           \]
           However, when $p>1/3$, the largest imaginary eigenvalue is in the
           continuous spectrum.
           A schematic showing all the eigenvalues, as well as the
           branches of continuous spectrum, in the case $0<p<1/3$
is shown
           in Figure \ref{fig_cont-spec-1}.

\section{Complex $\mu-$eigenvalues in the system at infinity}
\setcounter{equation}{0}
\label{sec-complex-eigs}

The case where the $\mu-$eigenvalues, in the
system at infinity when $\lambda=0$, form a complex quartet,
has a number of interesting features
that are outside the scope of this paper.  In this section,
some of the issues
are highlighted and then encapsulated in an open question.
\begin{figure}[ht]
\begin{center}
  \includegraphics[width=4.0cm]{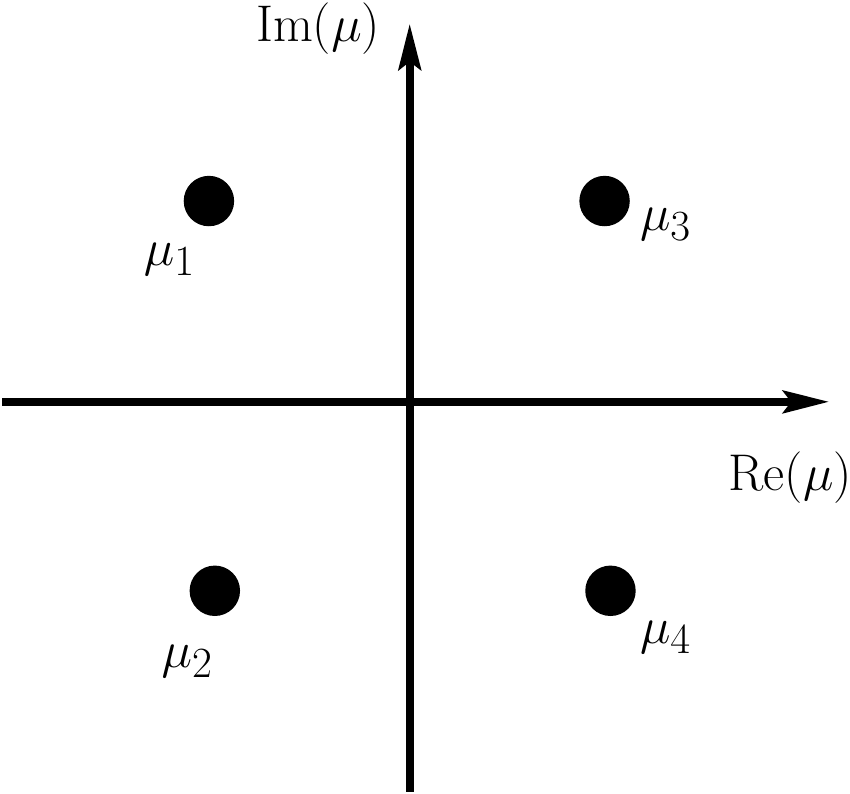}\hspace{1.0cm}
  \includegraphics[width=4.0cm]{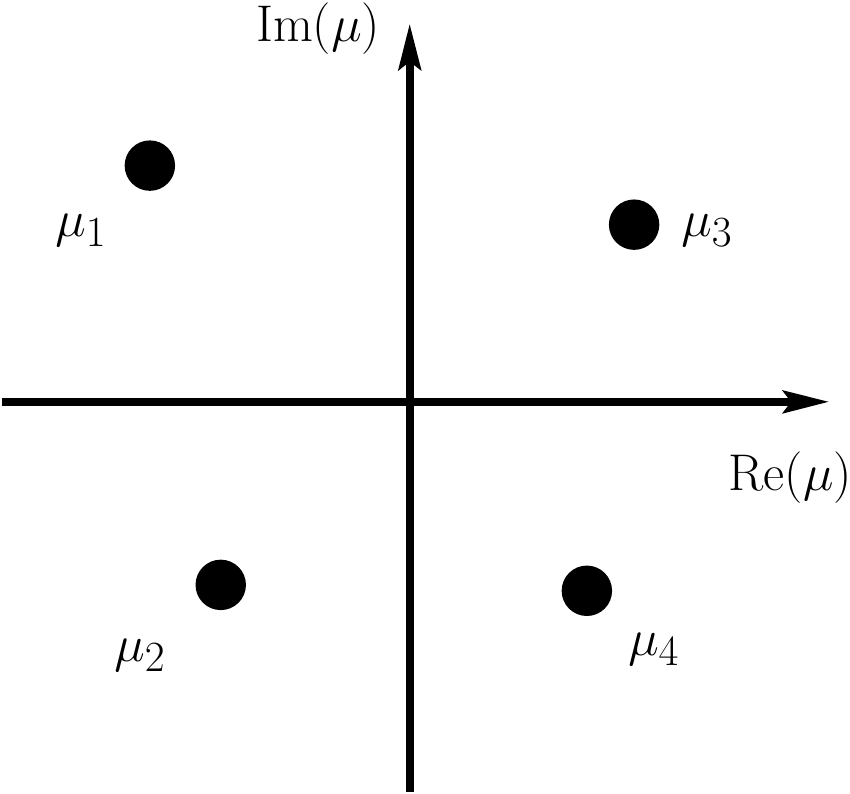}
\end{center}
\caption{Complex $\mu-$eigenvalues when $\lambda=0$ (left figure)
  and when $\lambda\neq0$ (right figure).}
\label{fig_mu-complex-eigs}
\end{figure}

A schematic of the $\mu-$eigenvalues in this
case when $\lambda=0$ and $\lambda\neq0$ is shown in Figure \ref{fig_mu-complex-eigs}.
When $\lambda=0$ the four eigenvalues form a complex quartet, symmetric
   about the origin numbered so that
\begin{equation}\label{mu-alpha-beta}
  \mu_1 = -\beta + \ri \alpha\,,\quad \mu_2=\overline{\mu_1}\,,\quad
  \mu_3 = \beta+\ri\alpha\,,\quad\mu_4=\overline{\mu_3}\,,
\end{equation}
with $\alpha>0$ and $\beta>0$.
The basic state $\Zh(\xi,c)$ will have the asymptotic form
\begin{equation}\label{mu-complex-Zh}
\lim_{\xi\to\infty}\re^{2\beta\xi} \Zh(\xi,c) = \big(c_1\cos(\alpha\xi+\nu)
+c_2\sin(\alpha\xi+\nu)\big)\zeta_1
+ \big(-c_1\sin(\alpha\xi+\nu)+c_2\cos(\alpha\xi+\nu)\big)\zeta_2\,,
\end{equation}
for some constants $c_1,c_2$ and phase shift $\nu$, with $\zeta_1$
and $\zeta_2$ real valued and satisfying
\[
  [{\bf B}^\infty(c) - (-\beta+\ri\alpha)\J(c)](\zeta_1+\ri\zeta_2)=0\,.
  \]
  The ratio of the constants $c_1/c_2$, $c_2\neq0$, is fixed on a given
  solitary wave.  Solitary waves of this type are called ``oscillatory
  solitary waves''.  The exponential decay assures that these solitary
  waves are still
  integrable on the real line, and so the momentum (\ref{H-I-def}) will
  be well defined. 

  The question is how to compute the three key ingredients in the
  derivative formula: $dI/dc$, $\Pi$, and $\chi$, and more importantly,
  how to adjust the proof of the derivative formula in Theorem 6.1.

  The calculation of $\chi$ is the least clear as the reality
  of the asymptotics is used in the definition of $\chi$
  in \S\ref{subsec-varrho}.
  On the other hand, with the momentum well defined, $dI/dc$ is in
  principal computable. In fact, a calculation of
  the invariants $H$ and $I$ for a class of
  oscillatory solitary waves is given
  in \textsc{Buryak \& Akhmediev}~\cite{ba95}.
  
  The tangent vector $\Zh_\xi$ will have a similar asymptotic
  form to (\ref{mu-complex-Zh}) with different constants,
  and $\ba^{+}$ is any tangent vector transverse to $\Zh_\xi$ with
  $\ba^{+}\to0$ as $\xi\to+\infty$.  There is a similar definition of $\ba^-$
  as $\xi\to-\infty$. Hence construction and normalization of $\Pi$ is
  similar to the real $\mu-$eigenvalue
  case in (\ref{Pi-normalization}) and \S\ref{subsec-lazutkin}.
Indeed, an example of the construction of ${\bf a}^\pm$ for the case of
  a homoclinic orbit with complex quartet eigenvalues in the system at
  infinity, and an
  explicit construction of the Lazutkin invariant for this
  case, is given in \textsc{Gaivao \& Gelfreich}~\cite{gg11}.
  There they
  use the Luzutkin invariant (called ``the homoclinic invariant'' in
  \cite{gg11}) to measure the splitting of separatices in
  the unfolding of the Hamiltonian Hopf bifurcation.

  When $\lambda\neq0$ the
 asymptotics (\ref{bu-asymptotics}) are defined exactly as before
\begin{equation}\label{bu-asymptotics-complex}
\begin{array}{rclrcl}
&&\displaystyle\lim_{\xi\rightarrow+\infty} e^{-\mu_1(\lambda)\xi}\bu_1(\xi,\lambda) = \zeta_1(\lambda), \quad &&\displaystyle\lim_{\xi\rightarrow+\infty} e^{-\mu_2(\lambda)\xi}\bu_2(\xi,\lambda) = \zeta_2(\lambda), \\[2mm]
&&\displaystyle\lim_{\xi\rightarrow-\infty} e^{-\mu_3(\lambda)\xi}\bu_3(\xi,\lambda) = \zeta_3(\lambda), \quad &&\displaystyle
  \lim_{\xi\rightarrow-\infty} e^{-\mu_4(\lambda)\xi}\bu_4(\xi,\lambda) = \zeta_4(\lambda)\,,
\end{array}
\end{equation}
with the complex vectors $\zeta_j$ defined as in (\ref{zeta-def}).
But the relation between $\bu_j(\xi,0)$ and $\Zh_\xi$ and $\ba^{\pm}$
will be more complex, something like
  \[
  \begin{array}{rcl}
    c_1 \bu_3(\xi,0) + \overline{c_1}\bu_4(\xi,0) &=& \Zh_\xi \\[2mm]
    c_2 \bu_3(\xi,0) + \overline{c_2}\bu_4(\xi,0) &=& \ba^- \\[2mm]
    c_3 \bw_3(\xi,0) + \overline{c_3}\bw_4(\xi,0) &=& \Zh_\xi \\[2mm]
    c_4 \bw_3(\xi,0) + \overline{c_4}\bw_4(\xi,0) &=& \ba^+\,.
  \end{array}
  \]
  The number of constants may be reducible, but the form
  of the limit still complicates the proof of the derivative formula
  in \S\ref{sec-deriv-evans}.

  The existence of oscillatory solitary waves is a difficult subject itself.
  Moreover when such a homoclinic orbit does exist, and is transversely
  constructed, then, by a theorem of
  \textsc{Devaney}~\cite{dev76}, there exists a countable set of multi-pulse
  homoclinic orbits nearby in parameter space.  This complicates the stability
  question.

  On the other hand, the Maslov index of multi-pulse solitary waves, when the system at
  infinity has complex $\mu-$eigenvalues, is computed in \cite{cdb09},
  and so $\Pi$ can be obtained in principle as the parity of the
  Maslov index (subject to orientation and normalization).
  It is clear from those results that the Maslov index
  can take any natural number on oscillatory solitary waves
  with more and more humps, suggesting a complicated stability picture.

  The open question, regarding the stability of oscillatory solitary waves, from the
  viewpoint of this paper, is: {\it does the derivative formula carry over to this case
  \[
  D''(0) = \chi \Pi \frac{dI}{dc}\,,
  \]
  and if so what is the definition of $\chi$?}

\section{Concluding remarks}
\setcounter{equation}{0}
\label{sec-cr}

There are two assumptions that are worthy of further comment:
the simple $\mu-$eigenvalue
assumption, and the restriction of four-dimensional phase space.

When the $\mu-$eigenvalues are simple at $\lambda=0$ they are simple for $\lambda$
near zero and so, with the emphasis of the paper on
        the derivatives of $D(\lambda)$ at $\lambda=0$, the assumption
        of simple eigenvalues is adequate.  For global (in $\lambda$)
        analysis, this assumption can be relaxed by working on
        exterior algebra spaces, where $\zeta_1\wedge\zeta_2$ and
        $\zeta_3\wedge\zeta_4$ can be constructed to be analytic for
        all $\lambda\in\Lambda$
        (cf.\ Appendix \ref{app-a} for the formulation in terms of exterior algebra),
        or maximally
        analytic individual eigenvectors can be constructed \cite{bd03}.

The assumption of a four dimensional phase space in the steady problem
(\ref{Zh-ode}) is sufficient to capture the essence of the theory, and
avoids unnecessary complexity.  When the phase space
has dimension $2n$ with $n>2$, the theory has a straightforward
generalization when all the $\mu-$eigenvalues are real and simple
when $\lambda=0$.
In this case the off diagonal terms in the matrix (\ref{D-lambda-def})
are of higher order in $\lambda$ and so the derivative formula is based
on the diagonal entries, and moreover it does not matter to which eigenvalue/eigenvector
$\Zh_\xi$ is asymptotic to as $\xi\to\pm\infty$.  The only change is that
$\Pi$ will have to be adapted to the new dimension of the stable and unstable
manifolds.  The $2-2$ splitting would be replaced by an $n-n$ splitting
and (\ref{Es-Eu-def}) would be replaced by
  \begin{equation}\label{Es-Eu-def-2n}
  E^s(\xi,0) = {\rm span}\big\{ \Zh_\xi,{\bf a}_1^+,\ldots,{\bf a}_{n-1}^+\big\}
  \qand
  E^u(\xi,0) = {\rm span}\big\{ \Zh_\xi,{\bf a}_1^-,\ldots,{\bf a}_{n-1}^-\big\}\,,
  \end{equation}
  with $\Pi$ replaced by the Lazutkin-Treschev invariant
  \begin{equation}\label{Pi-multidimen}
  \Pi = {\rm det}\left[\begin{matrix}
\bO({\bf a}_1^-,{\bf a}_1^+) & \cdots & \bO({\bf a}_1^-,{\bf a}_{n-1}^+) \\
\vdots & \ddots & \vdots \\
\bO({\bf a}_{n-1}^-,{\bf a}_1^+) & \cdots & \bO({\bf a}_{n-1}^-,{\bf a}_{n-1}^+)
    \end{matrix}\right]\,,
  \end{equation}
  (cf.\ \textsc{Treschev}~\cite{treschev} and \textsc{Chardard \& Bridges}~\cite{cb15}).  Hence, subject to the generalizations required in the Evans function
  construction, we expect that Theorem 6.1 generalizes with
  $\Pi$ replaced by (\ref{Pi-multidimen}).  
  The
  higher-dimensional case is better treated as it arises in applications.
  The case where the phase space has dimension $8$, with a four-four
  splitting of eigenvalues is studied by \textsc{Burchell}~\cite{burchell-thesis}. This case is relevant for the study of the transverse instability of solitary
  waves in the case of multisymplectic Dirac
  operators in two space dimensions and time.

  \section*{Acknowledgements}

  The first author acknowledges support from an EPSRC
  Doctoral Training Partnership,
  grant numbers EP/M508160/1 and EP/R513350/1.
The second author is grateful to Fr\'ed\'eric Chardard (Universit\'e Jean Monnet) for helpful discussions at an early stage of this project.  
\vspace{0.5cm}

\begin{center}
\hrule height.15 cm
\vspace{.2cm}
--- {\Large\bf Appendix} ---
\vspace{.2cm}
\hrule height.15cm
\end{center}
\vspace{.25cm}

\begin{appendix}

\renewcommand{\theequation}{A-\arabic{equation}}
\section{Matrix representation of the Evans function}
\label{app-a}
\setcounter{equation}{0}

For the purposes of this appendix $V$ is a four-dimensional complex vector
space with dual $V^*$, as the explicit use of a dual space makes
some of the calculations in exterior algebra clearer.
A basis for $V$ is $\{{\bf e}_1,{\bf e}_2,{\bf e}_3,{\bf e}_4\}$,
which is not necessarily the standard basis, and a basis for $V^*$ is
$\{{\bf e}_1^*,{\bf e}_2^*,{\bf e}_3^*,{\bf e}_4^*\}$,
normalized in the usual way $\delta_{ij} = \langle{\bf e}_i^*,{\bf e}_j\rangle$.

Exterior algebra spaces $\bwedge^k(V)$ and $\bwedge^k(V^*)$, $k=1,2,3,4$,
are defined
in the usual way with induced pairings, $\lbk\cdot,\cdot\rbk_k$.

The volume form
on $V$ is (\ref{orientation-1}) and on $V^*$ it is
\[
\vol^* = {\bf e}_1^*\wedge{\bf e}_2^*\wedge{\bf e}_3^*\wedge{\bf e}_4^*
    \,.
\]
The normalized basis then gives that $\lbk \vol^* , \vol \rbk_4 = 1$.

We will also find use for the equivalent bases
\begin{equation}\label{other-basis}
V={\rm span}\{\zeta_1,\ldots,\zeta_4\}\qand
V^*={\rm span}\{\J(c)\eta_1,\ldots,\J(c)\eta_4\}\,,
\end{equation}
where $\zeta_j$ and $\eta_j$ are the eigenvectors defined in (\ref{zeta-def})
and (\ref{eta-def}) respectively, with their associated four-forms
$\mathcal{V}(c,\lambda)$ in (\ref{orientation-2}) and
$\mathcal{V}^*(c,\lambda)$ in (\ref{orientation-2-dual}).

Although $V$ is a complex vector space, the pairing
between $V^*$ and $V$ does not involve conjugation (to be consistent
with the definition
of the dual vectors $\eta_j$, as in (\ref{eta-def}),
where conjugation is built in).  An example of the pairing on
$\bwedge^2$ is
\[
\lbk {\bf a}\wedge{\bf b} , {\bf c}\wedge{\bf d}\rbk_2 =
     {\rm det}\left(\begin{matrix} \lbk{\bf a},{\bf c}\rbk_1 & \lbk{\bf a},{\bf d}\rbk_1\\
     \lbk{\bf b},{\bf c}\rbk_1 & \lbk{\bf b},{\bf d}\rbk_1
     \end{matrix}\right)\,,\quad \mbox{for any}\ {\bf a},{\bf b}\in V^*\,,\quad
       {\bf c},{\bf d}\in V\,,
       \]
  for a decomposable element with obvious extension to arbitrary elements in
  $\bwedge^2$.  When reverting back to the notation in the body of the paper,
  $\lbk\cdot,\cdot\rbk_1$ is replaced by $\langle\cdot,\cdot\rangle$.

  The interior product $\intp$ is the adjoint operator associated with
  the of the wedge product,
  for example
  \[
  \lbk {\bf c}\intp ({\bf a}\wedge{\bf b}),{\bf d}\rbk_1 =
  \lbk {\bf a}\wedge{\bf b},{\bf c}\wedge {\bf d}\rbk_2\,,\quad
  {\bf a},{\bf b}\in V^*\,,\quad {\bf c},{\bf d}\in V\,.
  \]
  It is a contraction operator and is sometimes written
  ${\bf c}\intp ({\bf a}\wedge{\bf b}) = {\bf i}_{\bf c}({\bf a}\wedge{\bf b})$.
  
  The formula (\ref{Evans-uw-version}) for the Evans function
  can be proved using the Hodge
star operator as in \cite{bd99a,bd01}.  Here a new proof is
  given using the interior product. It is a little more general
  in that it does not require an inner product in its definition.

  The starting point is the Evans function in the form (\ref{evans-def-1}),
  relative to the orientation (\ref{orientation-1}). The interior product
  is used to transform it as follows,
  suppressing the arguments for brevity,
  \begin{equation}\label{D-transformation}
  \begin{array}{rcl}
D(\lambda)\vol &=& \re^{-\tau(c)\lambda\xi}
\bu_1\wedge\bu_2\wedge\bu_3\wedge\bu_4\\[2mm]
&=& \lbk \vol^*, \re^{-\tau(c)\lambda\xi}
\bu_1\wedge\bu_2\wedge\bu_3\wedge\bu_4\rbk_4\vol \\[2mm]
&=& \lbk \mathcal{V}^*(c,\lambda), \re^{-\tau(c)\lambda\xi}
\bu_1\wedge\bu_2\wedge\bu_3\wedge\bu_4\rbk_4\mathcal{V}(c,\lambda) \\[2mm]
&=& \lbk \re^{-\tau(c)\lambda\xi}(\bu_1\wedge\bu_2)\intp\mathcal{V}^*(c,\lambda), 
\bu_3\wedge\bu_4\rbk_2\mathcal{V}(c,\lambda) \,,
  \end{array}
  \end{equation}
  using $1 = \lbk\vol^*,\vol\rbk_4 = \lbk \mathcal{V}^*,\mathcal{V}\rbk_4$.
  Define the left argument in the pairing as
    \begin{equation}\label{WW-def}
    {\bf W}(\xi,\lambda) = \re^{-\tau(c)\lambda\xi}\bu_1(\xi,\lambda)\wedge
    \bu_2(\xi,\lambda)\intp \mathcal{V}^* \in \bwedge^2(V^*)\,.
    \end{equation}
    ${\bf W}$ is the left-hand side of the formula (\ref{W-Upsilon-identity}).
    Differentiate it with respect to $\xi$ 
    \begin{equation}\label{W-diff-eqn}
      {\bf W}_\xi = -\tau(c)\lambda
      \re^{-\tau(c)\lambda\xi}\bu_1(\xi,\lambda)\wedge
      \bu_2(\xi,\lambda)\intp \mathcal{V}^* +
      \re^{-\tau(c)\lambda\xi}{\bf A}^{(2)}\big(\bu_1(\xi,\lambda)\wedge
      \bu_2(\xi,\lambda)\big)\intp \mathcal{V}^*\,,
    \end{equation}
    using the definition ${\bf A}^{(2)}({\bf a}\wedge{\bf b}) :=
    {\bf Aa}\wedge {\bf b} + {\bf a}\wedge{\bf Ab}$ for any ${\bf a},{\bf b}\in V$, and the governing equation for ${\bf u}_j$ in (\ref{bu-eqn-def}).
    Now look at the second term.  It can be split using the identity
\[
        {\bf A}^{(2)}\big({\bf a}\wedge {\bf b}\big)\wedge{\bf c}\wedge{\bf d} +
        {\bf a}\wedge {\bf b}\wedge{\bf A}^{(2)}\big({\bf c}\wedge{\bf d}\big)
        ={\rm Trace}({\bf A}) 
        {\bf a}\wedge {\bf b}\wedge{\bf c}\wedge{\bf d}\,.
        \]
        Let
        ${\bf a},{\bf b}\in V$ be arbitrary and pair ${\bf a}\wedge{\bf b}\in\bwedge^2(V)$
        with the operator in the second term on the
        right-hand side of (\ref{W-diff-eqn}),
        \[
        \begin{array}{rcl}
        \lbk {\bf A}^{(2)}\big(\bu_1\wedge
        \bu_2\big)\intp \mathcal{V}^*, {\bf a}\wedge{\bf b}\rbk_2
        &=&       \lbk \mathcal{V}^*, {\bf A}^{(2)}\big(\bu_1\wedge
        \bu_2\big)\wedge{\bf a}\wedge{\bf b}\rbk_4\\[2mm]
        &=&    {\rm Trace}({\bf A})
        \lbk \mathcal{V}^*,\bu_1\wedge\bu_2\wedge{\bf a}\wedge{\bf b}\rbk_4 - \lbk \mathcal{V}^*, \bu_1\wedge
        \bu_2\wedge{\bf A}^{(2)}\big({\bf a}\wedge{\bf b}\big)\rbk_4\\[2mm]
        &=&    \tau(c)\lambda
        \lbk \big(\bu_1\wedge\bu_2\big)\intp\mathcal{V}^*,{\bf a}\wedge{\bf b}\rbk_2 - \lbk \big(\bu_1\wedge \bu_2\big)\intp \mathcal{V}^*, {\bf A}^{(2)}\big({\bf a}\wedge{\bf b}\big)\rbk_2\\[2mm]
        &=&    \tau(c)\lambda
        \lbk \big(\bu_1\wedge\bu_2\big)\intp\mathcal{V}^*,{\bf a}\wedge{\bf b}\rbk_2 \\[2mm]
        &&\hspace{2.0cm} - \lbk \big({\bf A}^{(2)}\big)^T\left[\big(\bu_1\wedge \bu_2\big)\intp \mathcal{V}^*\right], {\bf a}\wedge{\bf b}\rbk_2\,.
        \end{array}
        \]
        Since this holds for any ${\bf a},{\bf b}\in V$ we get that
        \[
        {\bf A}^{(2)}\big(\bu_1\wedge
        \bu_2\big)\intp \mathcal{V}^* = \tau(c)\lambda
        \big(\bu_1\wedge\bu_2\big)\intp\mathcal{V}^* -
        \big({\bf A}^{(2)}\big)^T\left[\big(\bu_1\wedge \bu_2\big)\intp \mathcal{V}^*\right]\,.
        \]
        Substitute into (\ref{W-diff-eqn})
        \[
        \begin{array}{rcl}
        {\bf W}_\xi &=& -\tau(c)\lambda
      \re^{-\tau(c)\lambda\xi}\bu_1(\xi,\lambda)\wedge
      \bu_2(\xi,\lambda)\intp \mathcal{V}^* +
      \re^{-\tau(c)\lambda\xi}{\bf A}^{(2)}\big(\bu_1(\xi,\lambda)\wedge
      \bu_2(\xi,\lambda)\big)\intp \mathcal{V}^*\\[2mm]
      &=& -\tau(c)\lambda {\bf W} + 
      \left( \tau(c)\lambda{\bf W} - \big({\bf A}^{(2)}\big)^T{\bf W}\right)\,,
        \end{array}
        \]
        or
        \begin{equation}\label{W-xi-eqn}
          {\bf W}_\xi = - \big({\bf A}^{(2)}(\xi,\lambda)\big)^T{\bf W}\,.
        \end{equation}
        This equation shows that ${\bf W}$ satisfies the adjoint equation on
        $\bwedge^2(V^*)$.

        Now write the right hand side of (\ref{W-Upsilon-identity}) as
    \begin{equation}\label{WW-w-eqn}
      \Upsilon(\xi,\lambda) = \J(c)\bw_3(\xi,\lambda)\wedge\J(c)\bw_4(\xi,\lambda)\,.
    \end{equation}
    Differentiate $\Upsilon$ and use the governing equation
    for $\bw$ in (\ref{bw-eqn}),
  \[
  \begin{array}{rcl}
    \Upsilon_\xi
    &=& \J(c)(\bw_3)_\xi\wedge\J(c)\bw_4 + \J(c)\bw_3\wedge\J(c)(\bw_4)_\xi \\[2mm]
    &=& \J(c){\bf A}(\xi,-\lambda)\bw_3\wedge\J(c)\bw_4 + \J(c)\bw_3\wedge\J(c){\bf A}(\xi,-\lambda)\bw_4 \\[2mm]
    &=& ({\bf B}(\xi,c)+\lambda\bM)\bw_3\wedge\J(c)\bw_4 + \J(c)\bw_3\wedge({\bf B}(\xi,c)+\lambda\bM)\bw_4 \\[2mm]
    &=& -{\bf A}(\xi,\lambda)^T\J(c)\bw_3\wedge\J(c)\bw_4 - \J(c)\bw_3\wedge {\bf A}(\xi,\lambda)^T\J(c)\bw_4 \,,
  \end{array}
  \]
  which simplifies to
  \begin{equation}\label{Upsilon-xi-eqn}
  \Upsilon_\xi =  -\big({\bf A}^{(2)}(\xi,\lambda)\big)^T \Upsilon\,.
  \end{equation}
  It is clear from (\ref{W-xi-eqn}) and (\ref{Upsilon-xi-eqn}) that
  ${\bf W}$ and $\Upsilon$ satisfy the same adjoint equation on
  $\bwedge^2(V^*)$.  But to show
  equality as in (\ref{W-Upsilon-identity})
  it is necessary to study their asymptotics.
  Using the asymptotic properties of $\bw_3$ and $\bw_4$ in (\ref{bw-asymptotics}),
  \[
   \lim_{\xi\to+\infty}\re^{(\mu_3+\mu_4)\xi}\Upsilon(\xi,\lambda)
   = \lim_{\xi\to+\infty}
  \left(
  \re^{\mu_3\xi}\J(c)\bw_3\right)\wedge
  \left( \re^{\mu_4\xi} \J(c)\bw_4\right) =
        \J(c)\eta_3\wedge\J(c)\eta_4 \,.
        \]
        On the other hand, using $\tau(c)\lambda=\mu_1+\mu_2+\mu_3+\mu_4$,
        \[
        \re^{(\mu_3+\mu_4)\xi}{\bf W}(\xi,\lambda) = \re^{-(\mu_1+\mu_2)\xi}\bu_1(\xi,\lambda)\wedge
        \bu_2(\xi,\lambda)\intp \mathcal{V}^*\,.
        \]
        Taking the limit $\xi\to+\infty$ and using the asymptotics of $\bu_1$ and
        $\bu_2$,
        \[
        \lim_{\xi\to+\infty}\re^{(\mu_3+\mu_4)\xi}{\bf W}(\xi,\lambda) = \lim_{\xi\to+\infty}\re^{-(\mu_1+\mu_2)\xi}\bu_1(\xi,\lambda)\wedge
        \bu_2(\xi,\lambda)\intp \mathcal{V}^* = \zeta_1\wedge\zeta_2\intp \mathcal{V}^*\,.
        \]
        This is where the second formula in (\ref{eta-identity}) comes into
        play.  It remains to prove that
        \begin{equation}\label{eta-identity-app}
        \J(c)\eta_3\wedge\J(c)\eta_4 = \zeta_1\wedge\zeta_2\intp \mathcal{V}^*\,.
        \end{equation}
        This identity is not immediately obvious, but a proof can be
        given using elementary linear algebra.
        Expand the right-hand side of (\ref{eta-identity-app}) in terms of
        a basis for $\bwedge^2(V^*)$ using (\ref{other-basis}), 
        \begin{equation}\label{cj-eqn}
        \begin{array}{rcl}
          \zeta_1\wedge\zeta_2\intp \mathcal{V}^*
          &=& c_1 \J(c)\eta_1\wedge\J(c)\eta_2 + c_2 \J(c)\eta_1\wedge\J(c)\eta_3 + c_3 \J(c)\eta_1\wedge\J(c)\eta_4 \\[2mm]
          &&\quad +  c_4 \J(c)\eta_2\wedge\J(c)\eta_3 + c_5 \J(c)\eta_2\wedge\J(c)\eta_4 + c_6 \J(c)\eta_3\wedge\J(c)\eta_4 \,.
        \end{array}
        \end{equation}
        A basis for $\bwedge^2(V)$ is the six-dimensional set
        $\{ \zeta_1\wedge\zeta_2, \ldots,\zeta_3\wedge\zeta_4\}$.
        Now pair both sides of (\ref{cj-eqn}) with each basis element from
        $\bwedge^2(V)$ in turn.  One finds that
        $c_1,\ldots,c_5$ are all zero, and $c_6=1$.  This
        proves (\ref{eta-identity}), and completes
        the proof of the formula (\ref{W-Upsilon-identity}).

        Now apply the identity ${\bf W}=\Upsilon$ to the Evans function
        in (\ref{D-transformation}),
  \begin{equation}\label{D-transformation-1}
  \begin{array}{rcl}
    D(\lambda)\vol &=&
 \re^{-\tau(c)\lambda\xi}
\bu_1\wedge\bu_2\wedge\bu_3\wedge\bu_4\\[2mm]
&=& \lbk \vol^*, \re^{-\tau(c)\lambda\xi}
\bu_1\wedge\bu_2\wedge\bu_3\wedge\bu_4\rbk_4\vol \\[2mm]
&=& \lbk \mathcal{V}^*(c,\lambda), \re^{-\tau(c)\lambda\xi}
\bu_1\wedge\bu_2\wedge\bu_3\wedge\bu_4\rbk_4\mathcal{V}(c,\lambda) \\[2mm]
&=&\lbk \re^{-\tau(c)\lambda\xi}(\bu_1\wedge\bu_2)\intp\mathcal{V}^*, 
\bu_3\wedge\bu_4\rbk_2\mathcal{V} \\[2mm]
&=& \lbk \J(c)\bw_3\wedge\J(c)\bw_4,\bu_3\wedge\bu_4\rbk_2\,\mathcal{V} \\[4mm]
&=&  {\rm det}\left[
    \begin{matrix} \bO(\bw_3,\bu_3) & \bO(\bw_3,\bu_4) \\
      \bO(\bw_4,\bu_3) & \bO(\bw_4,\bu_4) \end{matrix}\right]\,\mathcal{V}\,,
  \end{array}
  \end{equation}
  thereby completing the proof of Theorem 5.2.       
        $\hfill\blacksquare$.

\renewcommand{\theequation}{B-\arabic{equation}}
\section{Exact Evans function for the example}
\label{app-b}
\setcounter{equation}{0}

In this appendix the details of the construction of the exact
Evans function for the example in \S\ref{sec_coupledwave-eqn} are given.
When $\lambda\neq0$ the eigenvalues are
\[
\begin{array}{rcl}
  \mu_1 &=& c\alpha^2\lambda - \alpha\sqrt{ \alpha^2\lambda^2 + 4 + 3p}\\[2mm]
  \mu_2 &=& c\alpha^2\lambda - \alpha\sqrt{ \alpha^2\lambda^2 + 4}\\[2mm]
  \mu_3 &=& c\alpha^2\lambda + \alpha\sqrt{ \alpha^2\lambda^2 + 4}\\[2mm]
  \mu_4 &=& c\alpha^2\lambda + \alpha\sqrt{ \alpha^2\lambda^2 + 4 + 3p}\,.
\end{array}
\]
The $\mu-$eigenvalues $\mu(c,\lambda)$
are simple for all $\lambda\in\Lambda$, satisfying hypothesis ({\bf H6}).
The eigenvectors are
\[
\zeta_1 =
z_1\begin{pmatrix}
2 \\ 2 \varrho_1  - \mu_1 \\ 2 \mu_1 -  \varrho_1 \\ 
1\end{pmatrix}\,,
\ \zeta_2 = z_2\begin{pmatrix}
2 \\ 2 \varrho_2  - 4\mu_2 \\ 2 \mu_2 - 4 \varrho_2 \\ 
4\end{pmatrix}\,,
\ \zeta_3 = z_3\begin{pmatrix}
2 \\ 2 \varrho_3  - 4\mu_3 \\ 2 \mu_3 - 4 \varrho_3 \\ 
4\end{pmatrix}\,,
\ \zeta_4 =
z_4\begin{pmatrix}
2 \\ 2 \varrho_4  - \mu_4 \\ 2 \mu_4 -  \varrho_4 \\ 
1\end{pmatrix}\,,
\]
where $\varrho_j=\lambda+c\mu_j$.  The real numbers
$z_1,\ldots,z_4$ are arbitrary. 

Now construct the exact solutions, starting with (\ref{psi1eq}).
The two exact solutions with $\lambda\neq0$
for $\psi_1$ are
\begin{equation}\label{psi1-eqn}
\psi_1^{\pm} = \re^{\mp\alpha\sqrt{4+\alpha^2\lambda^2}\xi} A^{\pm}(\xi)\,,
\end{equation}
with $T={\rm tanh}(\alpha\xi)$ and
\[
A^{\pm}(\xi) = \pm a_0 + a_1 T \pm a_2 T^2 + T^3\,,
\]
where,
\[
a_0 = -\frac{\alpha^2\lambda^2}{15}\sqrt{4+\alpha^2\lambda^2}\,,\quad
a_1 = 1 + \frac{2}{5}\alpha^2\lambda^2\,,\quad
a_2 = -\sqrt{4+\alpha^2\lambda^2}\,.
\]
These two solutions are associated with the eigenvalues
\[
\psi_1^-\ :\ \mu_2 = c\alpha^2\lambda - \alpha\sqrt{4+\alpha^2\lambda^2}\qand
\psi_1^+\ :\ \mu_3 = c\alpha^2\lambda + \alpha\sqrt{4+\alpha^2\lambda^2}\,.
\]
Similarly, the two exact solutions with $\lambda\neq0$
for $\psi_2$ are
\begin{equation}\label{psi2-eqn}
\psi_2^{\pm} = \re^{\mp\alpha\sqrt{4+3p+\alpha^2\lambda^2}\xi} B^{\pm}(\xi)\,,
\end{equation}
with $T={\rm tanh}(\alpha\xi)$ and
\[
B^{\pm}(\xi) = \pm b_0 + b_1 T \pm b_2 T^2 + T^3\,,
\]
where
\[
b_0 = -\frac{(3p+\alpha^2\lambda^2)}{15}\sqrt{4+3p+\alpha^2\lambda^2}\,,\quad
b_1 = \frac{1}{5}(5+6p+2\alpha^2\lambda^2)\,,\quad
b_2 = -\sqrt{4+3p+\alpha^2\lambda^2}\,.
\]
These two solutions are associated with the eigenvalues
\[
\psi_2^-\ :\ \mu_1 = c\alpha^2\lambda - \alpha\sqrt{4+3p+\alpha^2\lambda^2}\qand
\psi_2^+\ :\ \mu_4 = c\alpha^2\lambda + \alpha\sqrt{4+3p+\alpha^2\lambda^2}\,.
\]

To simplify the construction of the Evans function, use the formula
(\ref{evans-def-1}) evaluated at $\xi=0$,
\begin{equation}\label{Dlambda-0}
    D_E(\lambda)\vol = {\bf u}_1(0,\lambda)\wedge
    {\bf u}_2(0,\lambda)\wedge{\bf u}_3(0,\lambda)\wedge{\bf u}_4(0,\lambda)\,.
\end{equation}
Hypothesis ({\bf H7}) is implicitly satisfied by using the original,
and equivalent, definition of the Evans function.

Using (\ref{u1u2-tilde-def}) and the above expressions for
    $\psi_j^\pm$ the exact solutions for $\bu_j(0,\lambda)$ are
\[
  {\bf u}_1(0,\lambda) = \begin{pmatrix}
    -2b_0 \\ -2\lambda b_0 + (2c-1)(\alpha b_1 - \mu_1 b_0) \\
    \lambda b_0 + (2-c) (\alpha b_1 - \mu_1 b_0) \\
    -b_0 \end{pmatrix}\,,\quad
  {\bf u}_2(0,\lambda) = \begin{pmatrix}
    -a_0 \\ -\lambda a_0 + (c-2)(\alpha a_1 - \mu_2 a_0) \\
    2\lambda a_0 + (1-2c) (\alpha a_1 - \mu_2 a_0) \\
    -2a_0 \end{pmatrix}\,,
  \]
  and
  \[
  {\bf u}_3(0,\lambda) = \begin{pmatrix}
    a_0 \\ \lambda a_0 + (c-2)(\alpha a_1 + \mu_3 a_0) \\
    -2\lambda a_0 + (1-2c) (\alpha a_1 + \mu_3 a_0) \\
    2a_0 \end{pmatrix}\,,\quad
  {\bf u}_4(0,\lambda) = \begin{pmatrix}
    2b_0 \\ 2\lambda b_0 + (2c-1)(\alpha b_1 + \mu_4 b_0) \\
    -\lambda b_0 + (2-c) (\alpha b_1 + \mu_4 b_0) \\
    b_0 \end{pmatrix}\,.
  \]
  Substituting these expressions into (\ref{Dlambda-0}) leads to
      \[
    D_E(\lambda))\vol = \Gamma(\lambda)
    \begin{pmatrix} 1\\ \lambda \\ -2\lambda \\ 2\end{pmatrix}\wedge
    \begin{pmatrix} 0 \\ c-2 \\ 1-2c \\ 0 \end{pmatrix}\wedge
         \begin{pmatrix} 2\\ 2\lambda \\ -\lambda \\ 1\end{pmatrix}\wedge
           \begin{pmatrix} 0 \\ 2c-1 \\ 2-c \\ 0 \end{pmatrix}\,,
           \]
           with
    \[
    \Gamma(\lambda) = a_0b_0
    (2\alpha a_1 +\mu_3a_0-\mu_2a_0)(2\alpha b_1 +\mu_4b_0-\mu_1b_0)\,.
    \]
    The latter wedge product evaluates to
    $9(c^2-1)\vol= -9\alpha^{-2}\vol$ where here
           $\vol$ is the standard volume form.  Therefore
           \[
           D_E(\lambda) = -\frac{9}{\alpha^2}a_0b_0(2\alpha a_1 +\mu_3a_0-\mu_2a_0)
           (2\alpha b_1 +\mu_4b_0-\mu_1b_0)\,.
           \]
           After some simplification we find
           \[
           D_E(\lambda) = -\frac{36 S}{(225)^2} (3+x^2)x^2(5-x^2)\,
(3+3p+x^2)(3p+x^2)(5-3p-x^2)\,,
           \]
           with $x^2=\alpha^2\lambda^2$ and
           \[
S = \sqrt{4+x^2}\sqrt{4+3p+x^2}\,.
\]
To synchronize with the Evans function used in \S\ref{sec_coupledwave-eqn}
we replace $D_E(\lambda))$ by its equivalent $\D_E(\lambda)\mapsto
4\chi D_E(\lambda)$,
\[
      D_E(\lambda) = -4\chi\frac{36 S}{(225)^2} (3+x^2)x^2(5-x^2)\,
      (3+3p+x^2)(3p+x^2)(5-3p-x^2)\,,
      \]
      which simplifies to
      \begin{equation}\label{Evans-exact-app}
      D_E(\lambda) = \frac{3 S}{16(225)^2}\alpha\lambda^2(3+x^2)(5-x^2)(3+3p+x^2)(3p+x^2)(5-3p-x^2)\,,
      \end{equation}
           This is the exact Evans function that is used in
           (\ref{exact-evans-example}).
           
      \end{appendix}

\end{document}